\documentclass{article}
\usepackage{graphicx}
\usepackage{color}
\usepackage{amsfonts}
\usepackage{amsthm}

\newtheorem{thm}{Theorem}[section]
\newtheorem{cor}[thm]{Corollary}
\newtheorem{lem}[thm]{Lemma}
\newtheorem{prop}[thm]{Proposition}

\newtheorem{defn}[thm]{Definition}

\newtheorem{rem}[thm]{Remark}

\newcommand{\p}{\varphi}
\newcommand{\s}{\psi}

\newcommand{\M}{\mathcal{M}}

\newcommand{\N}{\mathbb{N}}

\newcommand{\R}{{\mathbb{R}}}
\newcommand{\eps}{{\varepsilon}}

\newcommand{\sign}{{\rm{sign }}}

\begin{document}

\linespread{1.5}

\hoffset = -2cm \textwidth = 16.8cm

\title{Reparametrization invariant norms}
\author{P. Frosini$^{1,2}$\quad C.
Landi$^{1,3,}$\footnote{Corresponding author. E-mail Address:
\texttt{clandi@unimore.it}}\\\vspace{-0.3cm}
\small{$^{1}$ ARCES, Universit\`a di Bologna,}\\
\small{via Toffano $2/2$, I-$40135$ Bologna,
Italia}\\\vspace{-0.3cm} \small{$^{2}$Dipartimento di Matematica,
Universit\`a di
Bologna,}\\
\small{P.zza di Porta S. Donato 5, I-$40126$ Bologna,
Italia}\\\vspace{-0.3cm}
\small{$^3$DISMI, Universit\`a di Modena e Reggio Emilia,}\\
\small{via Amendola $2$, Pad. Morselli, I-$42100$ Reggio Emilia,
Italia}}
\maketitle
\abstract{This paper explores the concept of reparametrization invariant norm (RPI-norm) for $C^1_c$ functions, that is any norm invariant under composition with diffeomorphisms. The $L_\infty$-norm and the total variation norm are well known instances of RPI-norms. We prove the existence of an infinite family of RPI-norms, called standard RPI-norms, for which we exhibit both an integral and a discrete characterization. Our main result states that for every piecewise monotone function $\p$ in $C^1_c(\R)$ the standard RPI-norms of $\p$ allow us to compute the value of any other RPI-norm of $\p$. This is proved by using the standard RPI-norms in order to reconstruct the function $\p$ up to reparametrization and an arbitrarily small error with respect to the total variation norm.} 




\section{Introduction}
In recent papers the {\em natural pseudo-distance} $\sigma$
between manifolds endowed with regular real functions has been
studied as a tool for comparing the shape of manifolds (cf.
\cite{DoFr04a,DoFr04b, DoFrJEMS}). Each shape is represented by  pairs
$(\M,\p)$, where $\M$ is a connected manifold and $\p$ is
a real function defined on it (both $\M$ and $\p$ are supposed
to be sufficiently regular).  In this approach, the main
idea is to compare two diffeomorphic manifolds by measuring the
global change of the real functions they are endowed with, when
the manifolds  are deformed into each other:
$\sigma((\M,\p),(\mathcal N,\s))=\inf_h \sup_{p\in
\M}|\p(p)-\s\circ h(p)|$, where $h$ varies among all the
diffeomorphisms between $\M$ and $\mathcal N$.
We observe that $\sigma$ is a Fr\'echet-like
distance (cf., e. g., \cite{DuRo04}). Moreover, this line of research is strongly related
to the extensive study currently carried out about parametrization-independent shape comparison in Pattern Recognition
(cf, e. g., \cite{MiMu06}).

The definition of natural
pseudo-distance between the pairs $(\M,\p_1)$, $(\M,\p_2)$ can be reformulated as
 the value $\inf_{h\in D} F(\p_1-\p_2\circ h)$, where
$D$ denotes the set of all diffeomorphisms from $\M$ to $\M$  and
$F$ is the norm that takes each (sufficiently regular) function $\bar\p:\M\to \R$
to the number $\max_{P\in \M}|\bar \p(P)|$. In order that
$\inf_{h\in D} F(\p_1-\p_2\circ h)$ is a pseudo-distance, the key
property of the functional $F$ is that $F$ is a norm and $F(\bar
\p\circ h)=F(\bar\p)$ for every $\bar\p:\M\to \R$ and every $h\in
D$. In other words, the point is that $F$ is a {\em
reparametrization invariant norm}. Choosing a different
reparametrization invariant norm would allow us to obtain a
different pseudo-distance.

From this simple observation some natural questions arise. Are there
other reparametrization invariant norms on the vector space $V$ of
all (sufficiently regular) real functions defined on $\M$? What are their
properties? Are they induced by inner products? What kind of
information do they contain? Is it possible to reconstruct a
function $\p$ by using the values taken at $\p$ by the norms in the
set $RPI_V$ of all reparametrization invariant norms on $V$, or by
the norms in a suitable subset of $RPI_V$?

It is clear that the progress in this line of research requires
answers to these questions. This paper is a first step in this
direction, studying what happens in the simplest case,
i.e. $\M=\R$, when the considered reparametrizations are orientation-preserving.

We conclude this introduction by recalling that the invariance under
reparametrization appears to be relevant in several fields of research.
Just two examples are Statistics (cf., e. g., the Kolmogorov-Smirnov Test) and
the Theory of Interpolation Spaces (with reference, e. g., to the K-Method).

\subsection{The point and the main ideas in this paper}
This paper studies the reparametrization invariant norms that can
be defined  on a suitable set of  regular functions $\p$ from $\R$
to $\R$.
The norms $\max|\p|$, $\max\p-\min\p$, the total variation
$V_\p$ of $\p$, and the function $\sqrt{\max|\p|^2+V_\p^2}$ are
simple examples of reparametrization invariant norms, assuming
that the derivative of $\p$ has compact support and that $\p$ vanishes at $-\infty$. Actually,
there exist infinitely many examples of such norms, since each
linear combination with positive coefficients of reparametrization
invariant norms is obviously a reparametrization invariant norm.
However, in the set of all these norms, we have succeeded in
detecting a particular subset of norms, which we call {\em standard
reparametrization invariant (RPI-) norms}, such that
\begin{enumerate}
\item if the $C^1$-function $\p$ has compact support and is piecewise monotone, then the knowledge of all the standard RPI-norms of $\p$ allows us to reconstruct $\p$  up
to reparametrization, with an arbitrarily small error $\eps$ with
respect to the  total variation norm;
\item as a consequence of previous property, any other RPI-norm of such a function $\p$ is completely determined by the values taken on $\p$ by the standard RPI-norms.
\end{enumerate}
Therefore, we have focused our research on these norms.

 The main idea of this paper
originates from the following classical definition of the total
variation $V_\p$ for a regular function $\p:\R\to\R$ (see, e.g.,
\cite{AmFuPa00}):
$$V_\p=\sup_{\s\in \Psi}\left|\int_{-\infty}^{+\infty}\p(t)\cdot\frac{d\s}{dt}(t)\ dt\right|=
\sup_{\s\in
\Psi}\left|\int_{-\infty}^{+\infty}\p(-t)\cdot\frac{d\s}{dt}(t)\
dt\right|$$ where $\Psi$ is the set of all (sufficiently regular) functions $\s$
from $\R$ to $\R$ with $|\s|\le 1$. We observe that if we substitute
$\Psi$ with any subset $\hat \Psi$ of $\Psi$ \emph{that is closed
with respect to reparametrization}, then other reparametrization
invariant norms can be obtained (though, in this case, the two
suprema in the previous formula may be different). The closure with
respect to reparametrizations means that if $\hat \s\in \hat \Psi$
then $\hat \s\circ h\in \hat \Psi$ for every orientation-preserving
diffeomorphism $h:\R\to\R$.


In order to apply our idea, in first place  we  choose a functional
space. Many different choices are possible. As  a trade-off between
generality and simplicity we have chosen the space $AS^1(\R)$ of all
\emph{almost sigmoidal $C^1$ functions}. Roughly speaking, this
space could be defined as the space of all $C^1$ functions $\psi:\R\to \R$
that ``behave as a sigmoid outside a sufficiently large compact''
(see Section~\ref{genprop}, definition \ref{almsigm}). This choice is not very restrictive,
since $AS^1(\R)$ contains all the $C^1$ functions with compact
support.

By defining $[\s]$ as the set containing $\s\in AS^1(\R)$ and
all its re\-pa\-ra\-me\-tri\-za\-tions $\s\circ h$, and by setting
$$\|\p\|_{[\s]}=\sup_{\hat \s\in [\s]}\left|\int_{-\infty}^{+\infty}\p(-t)\cdot\frac{d\hat\s}{dt}(t)\ dt\right|$$
we obtain a reparametrizazion invariant norm on $AS^1(\R)$. The
norms obtained in this way are precisely the standard
reparametrizazion invariant norms, verifying the properties
described in previous statements 1 and 2
(Theorem~\ref{RT}). The reason for using $\p(-t)$
instead of $\p(t)$ inside the integral is that this choice allows
us to obtain the equality $\|\p\|_{[\s]}=\|\s\|_{[\p]}$, thanks to
the fact that our functions belong to $AS^1(\R)$. Incidentally,
this motivated also the choice of $AS^1(\R)$ as the functional
space to use.

We could proceed analogously for a general set $\hat \Psi$ closed with respect to reparametrization in
place of $[\s]$ and obtain a reparametrization invariant norm
$\|\p\|_{\hat\Psi}$, but the case $\hat \Psi=[\s]$ is the most
interesting one, since $\|\p\|_{\hat\Psi}$ can be easily expressed
as a supremum of standard reparametrization invariant norms.


In order to get the main results of this paper some technicalities
will be necessary. In particular, a key role will be played by the
Bounding Lemma~\ref{bounds}, asserting that, after normalization,
every reparametrization invariant norm $\|\p\|$ is upper bounded
by the total variation $V_\p$ and lower bounded by the value
$\lim_{t\to+\infty}|\p(t)|$. A stronger Bounding Lemma will be
proved for $C^1_c(\R)$. It asserts that, after normalization, every
reparametrization invariant norm $\|\p\|$ of a function $\p$
having compact support is upper bounded by half the total
variation $V_\p$ and lower bounded by the value $\max|\p(t)|$. The
proof of these key results will require some computations and a
preliminary study of the general properties of reparametrization
invariant norms, that will be carried out in Section~\ref{genprop}.
In particular, we shall examine the role played by two particular
functions, called $S$ and $\Lambda$. Moreover, in the same section
we shall prove the stability of RPI-norms with respect to small
perturbations in $C^1$ and the interesting fact that no inner
product can induce a RPI-norm.

The definition of standard reparametrization invariant norm
will be introduced in Section~\ref{SRPIN}, together with some examples and some basic properties.
However, in order to proceed further, we shall have to represent standard
reparametrizazion norms in a simpler way. We know that an
alternative definition exists for the total variation, saying that
$V_\p$ equals the value
$\sup_{n}\sup_{\tau_0\le\ldots\le\tau_i\le\ldots\le\tau_n}\sum_{i=0}^{n-1}
|\p(\tau_{i+1})-\p(\tau_i)|$.
In
Section~\ref{DROSRPIN} some computations will be necessary to
make available a similar representation also for standard reparametrization invariant norms (Theorem~\ref{teoremadirappresentazione}).
This new kind of representation will be used to prove the fundamental results in this paper, i.e. the possibility of
reconstructing piecewise monotone functions with compact support up to reparametrizations by means of the standard reparametrization invariant norms and the dependence of reparametrization invariant norms on standard reparametrization invariant norms (Section~\ref{DEP}).
Section~\ref{OP} will conclude this paper by illustrating some open problems.

\section{Re-pa\-ra\-me\-te\-ri\-za\-tion invariant norms: definition and general properties}
\label{genprop}

\subsection{Some notations and basic definitions}



In this paper the symbols $C^1(\R)$ and $C^1_c(\R)$ will represent
the set of all one time continuously differentiable functions
from $\R$ to $\R$ and the set of all functions in $C^1(\R)$ that
have compact support, respectively.
The symbol $D^1_+(\R)$ will represent
the set of all orientation-preserving $C^1$-diffeomorphisms
from $\R$ to $\R$.

We shall say that a function $f$ is {\em increasing} (\emph{strictly increasing}) if $t < t'$ implies $f(t) \le f(t')$ ($f(t) < f(t')$), and {\em decreasing}  (\emph{strictly decreasing}) if  $t< t'$ implies $f(t) \ge f(t')$ ($f(t) > f(t')$). A function will be called {\em monotone} if it is either decreasing or increasing, and \emph{strictly monotone} if it is either strictly increasing or strictly decreasing.

A number will be said to be {\em positive} when it is strictly greater than zero. The set of positive natural numbers will be denoted by $\mathbb{N}^+$.

First of all let us introduce the functional space we shall work
in.

\begin {figure}
 \begin {center}

  \includegraphics [height=5cm] {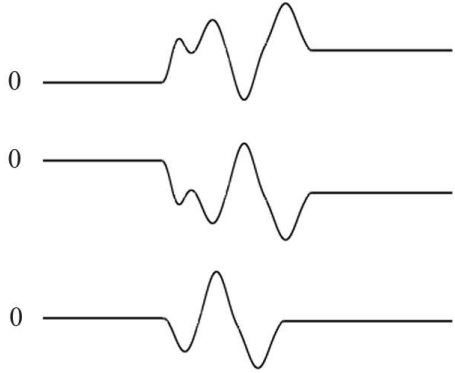}
  \caption {Three examples of almost sigmoidal functions.}

\label{sigfun}
 \end {center}
\end {figure}

\begin{defn}\label{almsigm}
Let us consider the set of all $C^1$-functions $\p:\R\to\R$ for
which two real values $a,b$ exist such that:

\begin{itemize}
  \item $\p(t)=0$ for every $t\in(-\infty,a]$;
  \item $\p(t)$ is constant in $[b,+\infty)$.
\end{itemize}

We shall denote this set by the symbol $AS^1(\R)$ and call each
function in $AS^1(\R)$ an \emph{almost sigmoidal function of class
$C^1$}.
\end{defn}
Examples of almost sigmoidal functions are shown in Figure \ref{sigfun}.

Obviously, $C^1_c(\R)\subset AS^1(\R)$ and every function in $AS^1(\R)$ has   bounded variation. We shall use the symbol $\textbf{0}$ to denote the almost sigmoidal function that vanishes everywhere.

The ideas described in this paper can be
extended to more general spaces, but we choose this setting in
order to simplify our proofs from the technical viewpoint.

For every $\s\in AS^1(\R)$ we shall denote by $V^+\langle\s\rangle(t)$ (resp.
$V^-\langle\s\rangle(t)$) the positive (resp. negative) variation of $\s$, and
by $V\langle\s\rangle(t)$ the variation of $\s$:
$$V^+\langle\s\rangle(t)=\int_{-\infty}^t\max\left\{\frac{d \s}{ds}(s),0\right\}\ ds,
\ V^-\langle\s\rangle(t)=\int_{-\infty}^t\max\left\{-\frac{d
\s}{ds}(s),0\right\}\ ds,$$
$$V\langle\s\rangle(t)=V^+\langle\s\rangle(t)+V^-\langle\s\rangle(t)=\int_{-\infty}^t\left|\frac{d \s}{ds}(s)\right|\ ds.$$
Since $\s\in C^1(\R)$, the functions $V^+\langle\s\rangle$, $V^-\langle\s\rangle$, and
$V\langle\s\rangle$ are $C^1$. Moreover we shall denote by $V^+_\s$,
$V^-_\s,V_\s$ the total positive variation, the total negative
variation and the total variation of $\s$, respectively:
$$V^+_\s=\int_{-\infty}^{+\infty}\max\left\{\frac{d \s}{ds}(s),0\right\}\ ds,
\ V^-_\s=\int_{-\infty}^{+\infty}\max\left\{-\frac{d
\s}{ds}(s),0\right\}\ ds,$$
$$V_\s=V^+_\s+V^-_\s=\int_{-\infty}^{+\infty}\left|\frac{d \s}{ds}(s)\right|\ ds.$$
We recall that $V^+\langle\s\rangle(t)$ and  $V^-\langle\s\rangle(t)$ are non-negative
increasing functions whose difference is exactly $\s$.
Observe that $V^+_\s-V^-_\s=\lim_{t\to+\infty}\s(t)$.
Obviously, if $\s\in C^1_c(\R)$ then
$V^+_\s=V^-_\s=\frac{1}{2}V_\s$.


\begin{defn}\label{rep}
For any two functions $\p_1,\p_2:\R\to \R$, we say that {\em
$\p_2$ is obtained from $\p_1$ by a reparametrization (of class
$C^1$)} if an orientation-preserving diffeomorphism  $h\in
D^1_+(\R)$ exists such that $\p_2=\p_1\circ h$. The diffeomorphism
$h$ will be called a \emph{reparametrization}. We denote by
$\sim$ the equivalence relation defined by setting $\p_2\sim \p_1$
if and only if $\p_2$ is obtained from $\p_1$ by a
reparametrization. The equivalence class of $\p_1$ in $AS^1(\R)$ will be denoted by $[\p_1]$.
\end{defn}

\begin {figure}
\begin{center}
\begin{picture}(0,0)%
\includegraphics[height = 5cm]{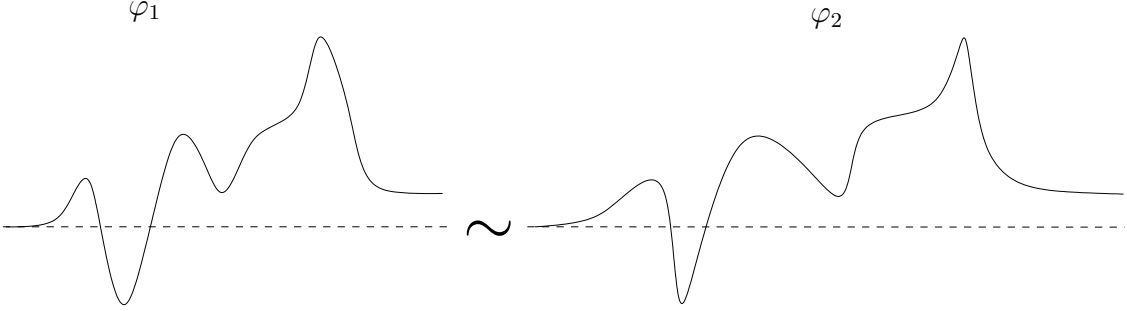}%
\end{picture}%
\setlength{\unitlength}{3947sp}%
\begingroup\makeatletter\ifx\SetFigFont\undefined%
\gdef\SetFigFont#1#2#3#4#5{%
  \reset@font\fontsize{#1}{#2pt}%
  \fontfamily{#3}\fontseries{#4}\fontshape{#5}%
  \selectfont}%
\fi\endgroup%
\begin{picture}(9272,3092)(1651,-3122)
\put(2449,-634){\makebox(0,0)[lb]{\smash{{\SetFigFont{12}{14.4}{\rmdefault}{\mddefault}{\updefault}{\color[rgb]{0,0,0}$\p_1$}}}}}
\put(6732,-681){\makebox(0,0)[lb]{\smash{{\SetFigFont{12}{14.4}{\rmdefault}{\mddefault}{\updefault}{\color[rgb]{0,0,0}$\p_2$}}}}}
\put(4550,-2080){\makebox(0,0)[lb]{\smash{{\SetFigFont{41}{49.2}{\rmdefault}{\mddefault}{\updefault}{\color[rgb]{0,0,0}$\sim$}}}}}
\end{picture}  
  \caption{We are interested in studying the norms that take both the functions
  $\p_1$ and $\p_2$ to the same value, since $\p_2$ is obtained by composing $\p_1$ with
  an orientation-preserving $C^1$-diffeomorphism of $\R$.}

\label{invarianza}
 \end{center}
\end {figure}

In this paper we study the norms that take equivalent functions to
the same value (see Figure~\ref{invarianza}).

%


\subsection{Reparametrization invariant norms}

Now we give the main definition in this paper.

\begin{defn}
Let us consider the real vector space $AS^1(\R)$. We say that a
norm $\|\cdot \|:AS^1(\R)\rightarrow \R$ is invariant under
reparametrization (or a {\em reparametrization invariant
norm})  if it is constant over each equivalence class of
$AS^1(\R)/\sim$.
\end{defn}

In the following the
reparametrization invariant norms will be often called
\emph{RPI-norms}.

The norms $\max|\p|$, $\max\p-\min\p$ and the total variation $V_\p$
of $\p$ are  simple examples of RPI-norms.

It is quite easy to see that infinitely many RPI-norms exist.
Indeed, it is trivial to prove that each linear combination with
positive coefficients of a finite number of RPI-norms  is still a
RPI-norm.

Another simple method to obtain a RPI-norm is to consider the $\sup$ of
a set of RPI-norms, under the assumption that such a $\sup$ is
finite at each point.

A third procedure consists in taking a norm $\|\cdot\|_*$ on $\R^k$
verifying the property that if $0\le x_i\le y_i$ for $1\le i\le k$
then $\|(x_1,\ldots,x_k)\|_*\le \|(y_1,\ldots,y_k)\|_*$. In this
case it is easy to verify that, if we have $k$ RPI-norms
$\|\cdot\|_1,\ldots,\|\cdot\|_k$, then  the function
$\|\p\|=\left\|\left(\|\p\|_1,\ldots,\|\p\|_k \right) \right\|_*$ is
a RPI-norm, too. E.g., the function $\sqrt{\max|\p|^2+V_\p^2}$ is a
RPI-norm. The hypothesis about $\|\cdot\|_*$ is necessary, as can be
seen by taking the norm $\|(x_1,x_2)\|_*=\sqrt{x_1^2+(x_1-x_2)^2}$
on $\R^2$, and setting $\|\p\|_1=\max|\p|$, $\|\p\|_2=V_\p$. Indeed
in this case $\|\cdot\|_*$ does not verify our hypothesis and
$F(\p)=\sqrt{\max|\p|^2+(\max|\p|-V_\p)^2}$ is \emph{not} a norm on
$AS^1(\R)$, since the triangular inequality does not hold (e.g.,
$F(\p_1+\p_2)>F(\p_1)+F(\p_2)$ if $\p_1(t)=\Lambda(t)$ and
$\p_2=\Lambda(t-4)$, where $\Lambda$ is the function defined in the
next Definition~\ref{SeLambda}).

Finally, another way to obtain a RPI-norm comes from the K-method in
the Theory of Interpolation Spaces (cf.~\cite{BeLo76,HoPe69}). Let
us consider two RPI-norms $\|\cdot\|_1$, $\|\cdot\|_2$ and the
function $K_p(t,\p)=\inf \left(\|\p_1\|_1^p+t^p\cdot
\|\p_2\|_2^p\right)^\frac{1}{p}$, where $t,p>0$ and the infimum is
computed for all possible decompositions $\p=\p_1+\p_2$ in
$AS^1(\R)$. Then $K_p(t,\cdot)$ is a RPI-norm.

\begin{rem}
Let $\|\cdot\|$ be a RPI-norm on $AS^1(\R)$. If $\p$ has  compact
support, also the composition $\p\circ h$ of $\p$ with an
\emph{orientation-reversing}
 $C^1$-diffeomorphism $h$ belongs to $AS^1(\R)$. Hence it makes sense to
 ask if  $\|\p\|$ equals $\|\p\circ h\|$ or not. In other words, the question is whether
 RPI-norms, that are invariant under  orientation-preserving reparametrizations by definition, are invariant   also under orientation-reversing reparametrizations.
   In general the answer is negative. As a counterexample, consider the RPI-norm $\|\p\|=\max_{t_1\le t_2}|2\p(t_1)-\p(t_2)|$.
\end{rem}

\begin {figure}
 \begin {center}

  \includegraphics [height=5cm] {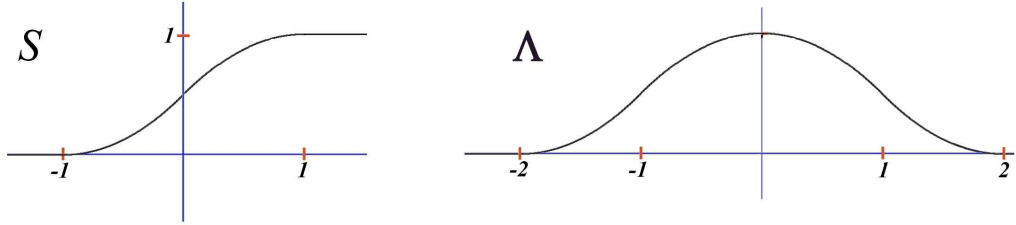}
  \caption {The graphs of the functions ${S}$ and $\Lambda$.}

\label{sigmadelta}
 \end {center}
\end {figure}

In order to proceed, we need to introduce two useful almost sigmoidal functions, represented in Figure \ref{sigmadelta}.

\begin{defn}
\label{SeLambda}
We shall denote by ${S}$  the  almost sigmoidal $C^1$-function
from $\R$ to $\R$ defined by setting
\[
{S}(t)=\left\{
\begin{array}{ll}
0 & \mbox{if $t< -1$} \cr \frac{(t+1)^2}{2} & \mbox{if $-1\le t\le
0$} \cr 1- \frac{(t-1)^2}{2}& \mbox{if $0<t\le 1$}\cr 1 & \mbox{if
$t>1$}
\end{array}
\right. .\]

We  define  $\Lambda:\mathbb{R}\to\mathbb{R}$ by setting
$\Lambda(t)=S(t+1)-S(t-1)$.
\end{defn}

In the following subsection we shall show that, in some sense,
each RPI-norm is controlled by the norm of the function $S$.

\subsubsection{The Bounding Lemma}\label{BLS}

The Bounding Lemma states that, after normalization, every
RPI-norm of $\p$ is bounded from above  by the total variation of
$\p$ and from below by $\lim_{t\to +\infty}|\p(t)|$.

This result will be proved as a consequence of the fact that  the
increasing functions of $AS^1(\R)$ can be approximated arbitrarily
well by functions  equivalent to  multiples  of the function $S$
(Prop.~\ref{sim}). From this, it follows that the RPI-norms of
monotone functions are multiples of the norm of $S$
(Prop.~\ref{simdef}).

\begin{prop}\label{sim}
Assume that a RPI-norm $\|\cdot \|$ is given. For any increasing
function $\p \in AS^1(\R)$ and any $\eps>0$, an increasing function
$\p_\eps \in AS^1(\R)$ exists such that $\|\p_\eps-\p\|=\eps$,
$\max\p_\eps=\max\p+\frac{\eps}{\|S\|}$ and
$\p_\eps\sim\max\p_\eps\cdot S$.
\end{prop}

\begin{proof}
Let $a$ and $b$ be two real numbers such that $a<b$, $\p$ equals
$0$ in the interval $(-\infty,a]$ and it is constant in the
interval $[b,+\infty)$. Let us define $\widetilde{S}(t)=
S\left(2\frac{t-a+\eps}{b-a+2\eps}-1\right)$ and
$\p_\eps=\p+\frac{\eps}{\|S\|}\cdot \widetilde{S}$. We point out
that $\|\widetilde{S}\|=\|S\|$, since $\widetilde{S}$ is obtained
by reparametrizing $S$. Hence $\|\p_\eps-\p\|=
\frac{\eps}{\|S\|}\cdot \|\widetilde{S}\|= \eps$. Moreover
$\p_\eps$ is clearly an increasing function belonging to
$AS^1(\R)$ and $\max\p_\eps=\max\p+\frac{\eps}{\|S\|}$.

Therefore,  we have only to prove that a reparametrization $h\in
D^1_+(\R)$ exists, such that $\p_\eps(t)= \max\p_\eps\cdot
S(h(t))$ for every $t\in \R$. In order to show this, we set
$$c=\frac{2}{b-a+2\eps}\cdot \sqrt{\frac{\eps}{\max\p_\eps\cdot
\|S\|}} $$ and define


\[
{h}(t)=\left\{
\begin{array}{ll}
c\cdot\left(t-a+\eps\right)-1 & \mbox{if $t\le a-\eps$} \cr
f_\eps\left(\p_\eps(t)\right) & \mbox{if $a-\eps < t < b+\eps$}
\cr c\cdot\left(t-b-\eps\right)+1 & \mbox{if $b+\eps\le t$}
\end{array}
\right. ,\] where $f_\eps:[0,\max\p_\eps]\to [-1,1]$ denotes the
inverse function of the restriction of $\max\p_\eps\cdot S$ to the
interval $[-1,1]$. Note that
$\max\p_\eps=\max\p+\frac{\eps}{\|S\|}$ is always positive. The
definition of $f_\eps$ implies that each $t\in (a-\eps,b+\eps)$ is
taken by $h$ to the unique point $h(t)$ for which

\begin{equation}
\label{eq1}
    \p_\eps(t)=\max\p_\eps \cdot S(h(t)).
\end{equation}

On the one hand, we observe that if $t\le a$ then the equality
$\p_\eps(t)=\frac{\eps}{\|S\|}\cdot \widetilde{S}(t)$ holds, and
hence for $a-\eps<t\le a$ the equality~(\ref{eq1}) becomes

\begin{equation}
\label{eq2}
\frac{\eps}{\|S\|}\cdot \widetilde{S}(t)=\max \p_\eps \cdot S(h(t)).
\end{equation}

If $t$ is also close enough to $a-\eps$ we have from (\ref{eq2}) that
$-1\le h(t)\le 0$ and in this case $\widetilde{S}(t)=
2\left(\frac{t-a+\eps}{b-a+2\eps}\right)^2$ and
$S(h(t))=\frac{\left(h(t)+1\right)^2}{2}$, because of the
definitions of $\widetilde{S}$ and $S$. Then, by a direct
computation we obtain from (\ref{eq2}) that if $a-\eps<t\le a$ and $t$ is close
enough to $a-\eps$ the equality $h(t)=c\cdot (t-a+\eps)-1$ holds.

On the other hand, if $t\ge b$ the equality
$\p_\eps(t)=\max\p+\frac{\eps}{\|S\|}\cdot \widetilde{S}(t)$
holds, and hence for $b\le t<b+\eps$ the equality~(\ref{eq1}) becomes

\begin{equation}
\label{eq3}
\max\p+\frac{\eps}{\|S\|}\cdot \widetilde{S}(t)=\max \p_\eps\cdot S(h(t)).
\end{equation}

If $t$ is also close enough to $b+\eps$ we have from (\ref{eq3})
that $0\le h(t)\le 1$ and in this case $\widetilde{S}(t)=1-
2\left(\frac{t-b-\eps}{b-a+2\eps}\right)^2$ and
$S(h(t))=1-\frac{\left(h(t)-1\right)^2}{2}$, once more because of
the definitions of $\widetilde{S}$ and $S$. Then, by a direct
computation (recalling that
$\max\p_\eps=\max\p+\frac{\eps}{\|S\|}$) we obtain from (\ref{eq3}) that if $b\le
t<b+\eps$ and $t$ is close enough to $b+\eps$ the equality
$h(t)=c\cdot (t-b-\eps)+1$ holds.

It follows that $h$ is differentiable at both the points $a-\eps$
and $b+\eps$, and that at both of them the derivative of $h$ takes
the positive value $c$.

Furthermore, we observe that the restriction of $h$ to the open
interval $(a-\eps,b+\eps)$ has positive
derivative, since both the derivative of $\p_\eps$ is positive in
this interval (due to the addend $\frac{\eps}{\|S\|}\cdot
\widetilde{S}$) and the derivative of $f_\eps$ is positive in the
open interval $(0,\max\p_\eps)$.
Also,  $h$ has  obviously  derivative equal to the positive value
$c$ outside the interval $[a-\eps,b+\eps]$, because of its definition, and at points $a$ and $b$, since it is $C^1$.
In conclusion, we have shown that  $h$ is an orientation-preserving $C^1$-diffeomorphism.

As a final step, it is easy to verify that $\p_\eps(t)=
\max\p_\eps\cdot S(h(t))$ for every $t\in \R$. Indeed we already
know that $\p_\eps(t)= \max\p_\eps\cdot S(h(t))$ for
$a-\eps<t<b+\eps$. For $t\le a-\eps$ we have $h(t)\le -1$ and
hence $\p_\eps(t)=0= \max\p_\eps\cdot S(h(t))$, while  for $t\ge
b+\eps$ we have $h(t)\ge 1$ and hence $\p_\eps(t)=\max\p_\eps=
\max\p_\eps\cdot S(h(t))$.

Therefore $\p_\eps$ is equivalent to the function
$\max\p_\eps\cdot S$ and our statement is proved.
\end{proof}


Now we can prove the following simple but crucial result,
underlining the importance of the function $S$.

\begin{prop}\label{simdef}
Assume that a RPI-norm $\|\cdot \|$ is given. For any
monotone function $\psi \in AS^1(\R)$ we have that
$\|\psi\|=\max|\psi|\cdot \|S\|$.
\end{prop}

\begin{proof}
Set $\p=|\psi|$. By applying the previous Proposition~\ref{sim} and the triangular
inequality, we obtain that $\big| \max \p_\eps\cdot\|S\| -
\|\p\|\big|=\big| \|\p_\eps\| - \|\p\|\big|\le \|\p_\eps -\p\|=
\eps$. By passing to the limit for $\eps$ tending to $0$, we get the
equality $\max \p\cdot\|S\| -\|\p\|=0$ and our statement is proved.
\end{proof}

\begin{rem}
We have to justify our line of proof of Proposition~\ref{simdef},
since the passage through Proposition~\ref{sim} could appear a
little cumbersome. The point is that the function
$\frac{\p}{\max\p\cdot S}$ may not tend towards a positive finite
constant for $t\to a^+$ or for $t\to b^-$. Furthermore, it may
happen that the derivative of $\p$ vanishes in the open interval
$(a,b)$. In these cases we cannot change directly $\p$ into
$\max\p\cdot S$ by a reparametrization, i.e. by composing $\p$
with an (orientation-preserving) $C^1$-diffeomorphism. Hence we have
to change $\p$ into an approximation $\p_\eps$ that does not present
the previous problems.
\end{rem}

Now we are ready to prove the Bounding Lemma. It gives a lower bound and an upper bound for
each RPI-norm, involving the norm of $S$.

\begin{lem}[\textbf{Bounding Lemma}]\label{bounds}
Let $\|\cdot \|:AS^1(\R)\rightarrow \R$ be a
reparametrization invariant norm. Then, for every $\p\in
AS^1(\R)$ the following inequalities hold:
$$\lim_{t\to+\infty}|\p(t)|\cdot \|S\| \le \|\p\|\le V_\p\cdot \|S\|.$$
\end{lem}

\begin{proof}


We can write $\p=V^+\langle\p\rangle-V^-\langle\p\rangle$, with $V^+\langle\p\rangle,V^-\langle\p\rangle\in AS^1(\R)$. Hence
$$\left|\| V^+\langle\p\rangle\|-\|V^-\langle\p\rangle\|\right|\le
\|\p\|\le \| V^+\langle\p\rangle\|+\|V^-\langle\p\rangle\|.$$

Since $V^+\langle\p\rangle,V^-\langle\p\rangle$ are increasing, Proposition~\ref{simdef} implies that  $\|V^+\langle\p\rangle\|= V^+_\p\cdot
\|S\|$ and $\|V^-\langle\p\rangle\|= V^-_\p\cdot \|S\|$, and hence our
statement is proved by recalling that $V^+_\p-V^-_\p=\lim_{t\to+\infty}\p(t)$.
\end{proof}

\begin{rem}
The inequalities in the Bounding Lemma are sharp, as we can easily see by setting $\p=S$.
\end{rem}

\begin{cor}\label{corbounds}
Let $\|\cdot \|:AS^1(\R)\rightarrow \R$ be a reparametrization
invariant norm. Then $\|\Lambda\|\le 2\cdot \|S\|$.
\end{cor}

\begin{proof}
Set $\p=\Lambda$ in the previous Lemma~\ref{bounds}.
\end{proof}

\begin{rem}
We observe that the inequality proved in Corollary~\ref{corbounds}
is sharp, since $\|\Lambda\|$ can equal $2\cdot\|S\|$. This happens,
e.g., when we consider as RPI-norm the total variation. Moreover, it
is interesting to note that no positive constant $c$ exists such
that the inequality $c\cdot \|S\|\le\|\Lambda\|$ holds for every
RPI-norm $\|\cdot\|$. To see this, it is sufficient to consider the
RPI-norm $\|\p\|_k=\max|\p|+k\cdot\lim_{t\to+\infty}|\p(t)|$, for
$k\ge 0$. Since $c\cdot \|S\|_k=c\cdot(1+k)$ and $\|\Lambda\|_k=1$,
if $k$ is large enough the inequality $c\cdot
\|S\|_k\le\|\Lambda\|_k$ does not hold. Also in this sense, the
lower bound in the Bounding Lemma cannot be improved. Incidentally,
we observe that the function $\lim_{t\to+\infty}|\p(t)|$ defines a
seminorm on $AS^1(\R)$ that is reparametrization invariant.
\end{rem}

\subsubsection{Derivatives and RPI-norms}
\label{Der}

The main consequence of the Bounding Lemma is that the closeness of
two almost sigmoidal functions with respect to the  total
variation norm implies their closeness with respect to any other
RPI-norm. From this we obtain the next proposition, showing that if
the derivatives of two functions $\p,\s\in AS^1(\R)$ are everywhere
close to each other, then $\p$ and $\psi$ are close to each other
with respect to any other reparametrization invariant norm.

\begin{prop}\label{c1continuity}
Let $\|\cdot\|$ be a
reparametrization invariant norm on $AS^1(\R)$. Assume that $\p,\s\in AS^1(\R)$ and that the compact
support of their derivative is contained in the interval $[a,b]$ with $a\neq b$. If $\max
\left|\frac{d\p}{dt}-\frac{d\s}{dt}\right|\le \frac{\eps}{(b-a)\cdot\|S\|}$, then
$\|\p-\s\|\le\eps$.
\end{prop}

\begin{proof}
If $\max \left|\frac{d\p}{dt}-\frac{d\s}{dt}\right|\le \frac{\eps}{(b-a)\cdot\|S\|}$, we
have that
$$V_{\p-\s}=\int_{-\infty}^{+\infty}\left|\frac{d(\p-\s)}{dt}(t)\right|\ dt\le  \frac{\eps}{(b-a)\cdot\|S\|}\cdot (b-a)=\frac{\eps}{\|S\|}.$$
By applying the right inequality in the
Bounding Lemma~\ref{bounds} we obtain
$$\|\p-\s\|\le V_{\p-\s}\cdot \|S\|\le \frac{\eps}{\|S\|}\cdot \|S\|=\eps.$$
\end{proof}


\subsubsection{A stronger Bounding Lemma for functions in $C^1_c(\R)$}\label{BLL}

The inequalities in the Bounding Lemma can be improved if $\p$
belongs to $C^1_c(\R)\subseteq AS^1(\R)$. In this case $\|\p\|$
stays somewhere between the norms $\max |\p|$ and
$\frac{1}{2}V_\p$, up to the multiplicative constant
$\|\Lambda\|$. This section is devoted to prove these stronger
bounds.

In order to prove these results we need the following technical lemma,where $V_\chi$ and $V_{\bar\chi-\chi}$ denotes the total variation of $\chi$ and $\bar\chi-\chi$, respectively, on $[\alpha,\beta]$.

\begin{lem}\label{approxlemma}
Let $\chi:[\alpha,\beta]\to\R$ be a monotone function of class $C^1$ with $\chi(\alpha)\ne \chi(\beta)$ and vanishing derivatives at $\alpha$ and $\beta$.
Then, for any  $\eps>0$, a $C^1$-function $\bar\chi :[\alpha,\beta]\to\R$ and an orientation-preserving $C^1$-diffeomorphism $\bar h:[\alpha,\beta]\to[\alpha,\beta]$ exist such that
\begin{enumerate}
\item  if $\chi$ is increasing then $\frac{d\bar\chi}{dt}> 0$ in the open interval $(\alpha,\beta)$; if $\chi$ is decreasing then $\frac{d\bar\chi}{dt}< 0$ in $(\alpha,\beta)$;
\item $V_{\bar\chi-\chi}\le\eps$.
\item $\bar\chi\left(\bar h(t)\right)=\chi(\alpha)+\left(\chi(\beta)-\chi(\alpha)\right)\cdot S\left(2\cdot \frac{t-\alpha}{\beta-\alpha}-1\right)$ for every $t\in [\alpha,\beta]$ and there exists $\eta>0$ such that  $\bar h$ is the identity on the intervals $[\alpha,\alpha+\eta]$ and $[\beta-\eta,\beta]$.
\end{enumerate}
\end{lem}

\begin{proof}
 Let us define $\widehat{S}(t)=
S\left(2\cdot \frac{t-\alpha}{\beta-\alpha}-1\right)$. Let us choose two small values $\eta_2>\eta_1>0$ and an $\hat\eps>0$, and set

 $$\bar\chi(t)=\left\{\begin{array}{ll} \chi(\alpha)+2\cdot (t-\alpha)^2\cdot \frac{\chi(\beta)-\chi(\alpha) }{(\beta-\alpha)^2} &  \mbox{if $t\in[\alpha,\alpha+\eta_1]$}\\
 \left(\chi+\sign(\chi(\beta)-\chi(\alpha))\cdot \hat\eps\cdot \widehat{S}\right)\cdot\frac{V_\chi}{ V_\chi +\hat\eps} & \mbox{if $t\in[\alpha+\eta_2,\beta-\eta_2]$}\\
 \chi(\beta)-2\cdot (\beta-t)^2\cdot \frac{\chi(\beta)-\chi(\alpha) }{(\beta-\alpha)^2} &  \mbox{if $t\in[\beta-\eta_1,\beta]$}
 \end{array}\right.$$

 Let us notice that if $\chi$ is increasing, then $\bar\chi$ is also increasing in the intervals where it is  defined, and viceversa  if $\chi$ is decreasing, then $\bar\chi$ is also decreasing. Moreover, in the intervals where it is  defined, the derivative of $\bar\chi$ does not vanish, except at $\alpha$ and $\beta$ (where $\chi$ and $\bar\chi$ coincide).

\begin{figure}
\begin{center}
\begin{picture}(0,0)%
\includegraphics[height=4.5cm]{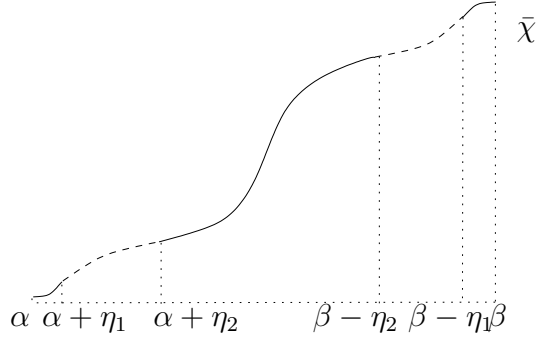}%
\end{picture}%
\setlength{\unitlength}{3947sp}%
\begingroup\makeatletter\ifx\SetFigFont\undefined%
\gdef\SetFigFont#1#2#3#4#5{%
  \reset@font\fontsize{#1}{#2pt}%
  \fontfamily{#3}\fontseries{#4}\fontshape{#5}%
  \selectfont}%
\fi\endgroup%
\begin{picture}(4342,1994)(1147,-3779)
\put(1162,-3702){\makebox(0,0)[lb]{\smash{{\SetFigFont{12}{14.4}{\rmdefault}
{\mddefault}{\updefault}{\color[rgb]{0,0,0}$\alpha$}%
}}}}
\put(1362,-3702){\makebox(0,0)[lb]{\smash{{\SetFigFont{12}{14.4}{\rmdefault}
{\mddefault}{\updefault}{\color[rgb]{0,0,0}$\alpha+\eta_1$}%
}}}}
\put(4344,-1881){\makebox(0,0)[lb]{\smash{{\SetFigFont{12}{14.4}{\rmdefault}
{\mddefault}{\updefault}{\color[rgb]{0,0,0}$\bar \chi$}%
}}}}
\put(2062,-3702){\makebox(0,0)[lb]{\smash{{\SetFigFont{12}{14.4}{\rmdefault}
{\mddefault}{\updefault}{\color[rgb]{0,0,0}$\alpha+\eta_2$}%
}}}}
\put(3062,-3702){\makebox(0,0)[lb]{\smash{{\SetFigFont{12}{14.4}{\rmdefault}
{\mddefault}{\updefault}{\color[rgb]{0,0,0}$\beta-\eta_2$}%
}}}}
\put(3662,-3702){\makebox(0,0)[lb]{\smash{{\SetFigFont{12}{14.4}{\rmdefault}
{\mddefault}{\updefault}{\color[rgb]{0,0,0}$\beta-\eta_1$}%
}}}}
\put(4162,-3702){\makebox(0,0)[lb]{\smash{{\SetFigFont{12}{14.4}{\rmdefault}
{\mddefault}{\updefault}{\color[rgb]{0,0,0}$\beta$}%
}}}}
\end{picture}%

 \caption{The function $\bar\chi$ used in Lemma \ref{approxlemma} (case
$\chi(\alpha)<\chi(\beta)$): it is quadratic near $\alpha$ and $\beta$, and without critical
points in $(\alpha,\beta)$. }

\label{approxchi}
 \end{center}
\end {figure}

 Now, if $\eta_1$ is small enough in comparison to $\eta_2$, we have that $\bar\chi(\alpha+\eta_2)-\bar\chi(\alpha+\eta_1)$ is positive or negative according to whether $\chi$, and hence $\bar\chi$, is increasing or decreasing. Analogously, $\bar\chi(\beta-\eta_1)-\bar\chi(\beta-\eta_2)$ is positive or negative according to whether $\chi$, and hence $\bar\chi$, is increasing or decreasing. Therefore, we can extend the definition of $\bar\chi$ to the open intervals $(\alpha+\eta_1,\alpha+\eta_2)$, $(\beta-\eta_2,\beta-\eta_1)$ in such a way that $\bar\chi$ is a $C^1$- function with non-vanishing derivative in the open interval $(\alpha,\beta)$.  Moreover, either $\bar\chi$ and $\chi$ are both increasing or both decreasing. So, $\bar\chi$ satisfies property {\em 1.}

Furthermore, if $\eta_1$ and $\eta_2$ have been chosen small enough, $\bar\chi$ satisfies also property {\em 2.}  Indeed,

\begin{eqnarray*}
&&V_{\bar\chi-\chi}=
\int_{\alpha}^{\alpha+\eta_2}\left|\frac{\ d\left(\bar\chi-\chi\right)}{dt}\right|\ dt+\int_{\alpha+\eta_2}^{\beta-\eta_2}\left|\frac{\ d\left(\bar\chi-\chi\right)}{dt}\right|\ dt+\int_{\beta-\eta_2}^{\beta}\left|\frac{\ d\left(\bar\chi-\chi\right)}{dt}\right|\ dt =\\
&=&
\int_{\alpha}^{\alpha+\eta_2}\left|\frac{\ d\left(\bar\chi-\chi\right)}{dt}\right|\ dt+\frac{\hat\eps}{V_\chi+\hat\eps}\cdot\int_{\alpha+\eta_2}^{\beta-\eta_2}\left|
\sign(\chi(\beta)-\chi(\alpha))\cdot V_\chi\cdot \frac{d\widehat{S}}{dt}-\frac{d\chi}{dt}
\right|\ dt+ \int_{\beta-\eta_2}^{\beta}\left|\frac{\ d\left(\bar\chi-\chi\right)}{dt}\right|\ dt\le\\
&\le&\int_{\alpha}^{\alpha+\eta_2}\left|\frac{\ d\bar\chi}{dt}\right|\ dt +
\int_{\alpha}^{\alpha+\eta_2}\left|\frac{\ d\chi}{dt}\right|\ dt+
\frac{\hat\eps}{V_\chi+\hat\eps}\cdot\left(\int_{\alpha+\eta_2}^{\beta-\eta_2}\left|
V_\chi\cdot \frac{d\widehat{S}}{dt}\right|\ dt+\int_{\alpha+\eta_2}^{\beta-\eta_2}\left|\frac{d\chi}{dt}\right|\ dt \right)+\\
& &+\int_{\beta-\eta_2}^{\beta}\left|\frac{\ d\bar\chi}{dt}\right|\ dt +
\int_{\beta-\eta_2}^{\beta}\left|\frac{\ d\chi}{dt}\right|\ dt\le\\
&\le& \left|\bar\chi(\alpha+\eta_2)-\bar\chi(\alpha)\right|+ \left|\chi(\alpha+\eta_2)-\chi(\alpha)\right|+\frac{\hat\eps}{V_\chi+\hat\eps}\cdot 2\cdot V_\chi +\left|\bar\chi(\beta)-\bar\chi(\beta-\eta_2)\right|+ \left|\chi(\beta)-\chi(\beta-\eta_2)\right|.
\end{eqnarray*}
Therefore, taking  $\eta_2$ (and hence $\eta_1$) small enough, by continuity we obtain that
$$V_{\bar\chi-\chi}\le 2\cdot \frac{\hat\eps}{V_\chi+\hat\eps}\cdot V_\chi+4\hat\eps.$$
It follows that if we choose $\hat\eps$, $\eta_1$ and $\eta_2$ small enough, then the inequality $V_{\bar\chi-\chi}\le \eps$ holds.

As for statement {\em 3.}, because of {\em 1.},  in the open interval $(\alpha,\beta)$, $\bar\chi$ admits the $C^1$ inverse
function $\bar\chi^{-1}$. We define $\hat h:(\alpha,\beta)\rightarrow (\alpha,\beta)$ by setting
$$\hat h=\bar\chi^{-1}\circ \left(\chi(\alpha)+\big(\chi(\beta)-\chi(\alpha)\big)\cdot \hat S\right).$$
Now, $\hat h$ is an orientation-preserving $C^1$-diffeomorphism because $\bar\chi$ and $\left(\chi(\alpha)+\big(\chi(\beta)-\chi(\alpha)\big)\cdot \hat S\right)$ are both increasing or both decreasing $C^1$-functions with non-vanishing derivatives.

Now,  by taking $\eta<\eta_1$, we obtain that $\hat h$ is the identity on $(\alpha,\alpha+\eta]\cup [\beta-\eta,\beta)$. This fact can be verified by  direct computations. Here the key point is that in the intervals $[\alpha,\alpha+\eta_1]$ and $ [\beta-\eta_1,\beta]$ the function $\bar\chi$ has been defined to be quadratic as $\hat S$. Therefore, if we  extend $\hat h$ to the closed interval $[\alpha,\beta]$ by taking $\bar h:[\alpha,\beta]\rightarrow [\alpha,\beta]$ with  $\bar h(\alpha)=\alpha$, $\bar h(\beta)=\beta$, and $\bar h(t)=\hat h(t)$ for $t\in (\alpha,\beta)$, then $\bar h$ is an  orientation-preserving $C^1$-diffeomorphism defined in $[\alpha,\beta]$ satisfying condition {\em 3.}
\end{proof}

Now we can prove the following result for $\Lambda$,
analogous to Proposition~\ref{simdef} proved for $S$.

\begin{prop}\label{simdefLambda}
Assume that a RPI-norm $\|\cdot \|$ is given. Assume also that $\p
\in C^1_c(\R)$, and that a value $\bar t\in \R$ exists such that
$\p$ is monotone both in $(-\infty,\bar t\  ]$ and in $[\bar
t,+\infty)$. Then $\|\p\|=\max|\p|\cdot \|\Lambda\|$.
\end{prop}

\begin{proof}
Let us assume that $\p\ne \mathbf 0$, otherwise the claim is trivial. Consider the smallest interval
$[a,b]$ containing   the compact support of
$\frac{d\p}{dt}$. Possibly by taking $-\p$ instead of $\p$, we can also assume that $\p$ is increasing in
 $[a,\bar t]$ and decreasing in $[\bar t, b]$ so that $\p\ge 0$. Furthermore, up to a reparametrization, we can assume that $a=-2$, $\bar t =0$ and $b=2$.

 Let $\chi_1$ denote  the restriction of $\p$ to the interval $[-2,0]$ and $\chi_2$ denote the restriction of $\p$ to the interval $[0, 2]$.  Let us apply Lemma \ref{approxlemma} for some  $\eps>0$ in order to obtain two functions $\bar \chi_1$ and $\bar \chi_2$ and the diffeomorphisms $\bar h_1$ and $\bar h_2$ such that $V_{\bar \chi_1-\chi_1}\le \frac{\eps}{2}$, $V_{\bar \chi_2-\chi_2}\le \frac{\eps}{2}$, $\bar\chi_1\left(\bar h_1(t)\right)=\p(0)\cdot S\left(2\cdot \frac{t+2}{2}-1\right)=\max|\p|\cdot S(t+1)$ and $\bar\chi_2\left(\bar h_2(t)\right)=\p(0)-\p(0)\cdot S\left(2\cdot \frac{t}{2}-1\right)=\max|\p|-\max|\p|\cdot S(t-1)$. Recall also that $\bar h_1$ is the identity in a neighbourhood of $-2$ and $0$, and $\bar h_2$ is the identity in a neighbourhood of $0$ and $2$.

 Consider the function $\p_\eps:\R\rightarrow \R$ in $AS^1(\R)$ defined by
 $$\p_\eps(t)=\left\{\begin{array}{ll}
 0 & \mbox{if $t \le -2$}\\
 \bar\chi_1(t) & \mbox{if $-2< t \le 0$}\\
  \bar\chi_2(t) & \mbox{ if $0< t\le 2$}\\
  0 & \mbox{if $t > 2$}\end{array}
 \right.$$
 We have that $\p_\eps$ is a function in $ C^1_c(\R)$,
with $V_{\p_\eps-\p}\le \eps$.
So, by applying the Bounding Lemma \ref{bounds}, we deduce that $\|\p_\eps-\p\|\le \eps\cdot \|S\|$.

Let us consider the orientation-preserving $C^1$-diffeomorphism $ h:\R\rightarrow \R$ defined by setting $ h(t)=\bar h_1(t)$ for $t\in [-2,0]$, $ h(t)=\bar h_2(t)$ for $t\in [0,2]$, $ h(t)=t$ otherwise. It holds that  $\p_\eps(h (t))=\max|\p|\cdot S(t+1)-\max|\p|\cdot S(t-1)=\max|\p|\cdot \Lambda$. Therefore,
 $\big| \max
\p\cdot\|\Lambda \|-\|\p\|\big|=\big| \|\p_\eps
\|-\|\p\|\big|\le \|\p_\eps-\p \|\le \eps\cdot \|S\|$. By passing to the
limit for $\eps$ tending to $0$, we get the equality $\max
\p\cdot\|\Lambda \|-\|\p\|=0$ and our statement is proved.
\end{proof}

The following result will be useful in the proof of the Reconstruction Theorem \ref{RT}. We omit its proof being quite similar to the ones used for Lemma \ref{approxlemma} and  Proposition~\ref{simdefLambda}.

\begin{prop}\label{generalizzazione}
 Let $\p\in AS^1(\R)$ admit $n$ points  $ t_0< t_1<\ldots < t_{n-1}$ such that
$\p$ is monotone  in each of the intervals $(-\infty,t_0]$, $[ t_0, t_1], \ldots , [ t_{n-2}, t_{n-1}]$, $[t_{n-1},+\infty)$. Then,  for any  $\eps>0$, an $AS^1$-function $\p_\eps :\R\to\R$ exists such that
\begin{enumerate}
\item  $V_{\p_\eps-\p}\le\eps$;
\item $\p_\eps\sim \p(t_0)\cdot S(t)+\sum_{i=1}^{n-1}\left(\p(t_{i})-\p(t_{i-1})\right)\cdot S(t-2i)$.
\end{enumerate}
In particular,  $\|\p\|= \left\|\p(t_0)\cdot S(t)+\sum_{i=1}^{n-1}\left(\p(t_{i})-\p(t_{i-1})\right)\cdot S(t-2i) \right\|$ for any  RPI-norm $\|\cdot \|$.
\end{prop}

Now we are ready to prove the stronger version of the Bounding Lemma
for functions with  compact support. It gives a
lower bound and an upper bound for each RPI-norm, involving the norm
of $\Lambda$.

\begin{lem}[\textbf{Bounding Lemma for $C^1_c(\R)$}]
\label{boundsC1}
Let $\|\cdot \|:AS^1(\R)\rightarrow \R$ be a
reparametrization invariant norm. Then, for every $\p\in
C^1_c(\R)$ the following inequalities hold:
$$\max|\p|\cdot \|\Lambda\| \le \|\p\|\le \frac{1}{2}V_\p\cdot \|\Lambda\|.$$
\end{lem}

\begin{proof}


First of all we prove the left inequality.  We take a point
$t_{\max}$ where $|\p|$ takes its maximum value and consider the
function

\[
{\psi}(\tau)=\left\{
\begin{array}{ll}
\int_{-\infty}^{\tau}\left|\frac{d\p}{dt}(t)\right|+\left|\frac{d\p}{dt}(2t_{\max}-t)\right|\
dt & \mbox{if $\tau\le t_{\max}$} \cr
\int^{+\infty}_{\tau}\left|\frac{d\p}{dt}(t)\right|+\left|\frac{d\p}{dt}(2t_{\max}-t)\right|\
dt & \mbox{if $\tau > t_{\max}$} \cr \end{array} \right. \] and set
$\hat \p=\p+\psi$.

We can easily verify that $\s$ is continuous also at $t_{\max}$, because of the two addends appearing in its definition. Then we observe that both $\hat\p$ and $\psi$ belong to $C^1_c(\R)$. In particular, the regularity of $\s$ follows from the fact that $t_{\max}$ is a critical point for $\p$.
Moreover, by computing their derivative we see that $\hat\p$ and $\psi$ are
increasing in $(-\infty,t_{\max}]$ and decreasing in
$[t_{\max},+\infty)$. Furthermore $\max \hat
\p=\hat \p(t_{\max})=\sign(\p(t_{\max}))\cdot\max|\p|+\max\s$.

Since $\p=\hat\p-\psi$, by applying Proposition~\ref{simdefLambda} with $\bar t=t_{\max}$
we get

\begin{equation}
\|\p\|=\|\hat\p-\psi\|\ge\big|\|\hat\p\|-\|\psi\|\big|=\big|\max\hat\p\cdot\|\Lambda\|-\max\psi\cdot\|\Lambda\|\big|=
\max|\p|\cdot\|\Lambda\|.
\end{equation}

As for the proof of the other inequality, we begin by  considering an interval
$[a,b]$ with $a\neq b$, such that the compact support of
$\frac{d\p}{dt}$ is contained in $[a,b]$. Let us define the function
$$F(\tau)=\int_{-\infty}^{\tau}\max\left\{\frac{d\p}{dt}(t),0\right\}\
dt- \int^{+\infty}_{\tau}\max\left\{-\frac{d\p}{dt}(t),0\right\}\ dt.$$
Since $F(a)=-V^-_\p\le 0$, $F(b)=V^+_\p\ge 0$ and $F$ is continuous,
a value $\bar t\in [a,b]$ exists such that $F(\bar t)=0$, i.e.
\[
\int_{-\infty}^{\bar t}\max\left\{\frac{d\p}{dt}(t),0\right\}\ dt=
\int^{+\infty}_{\bar t}\max\left\{-\frac{d\p}{dt}(t),0\right\}\ dt.
\]

Let us now assume that $\frac{d\p}{dt}(\bar t)=0$. In this case we
set

\[
{\p_1}(\tau)=\left\{
\begin{array}{ll}
\int_{-\infty}^{\tau}\max\left\{\frac{d\p}{dt}(t),0\right\}\ dt & \mbox{if $\tau\le \bar t$} \cr
\int^{+\infty}_{\tau}\max\left\{-\frac{d\p}{dt}(t),0\right\}\ dt & \mbox{if $\tau > \bar t$} \cr \end{array}
\right. \]
and
\[
{\p_2}(\tau)=\left\{
\begin{array}{ll}
\int_{-\infty}^{\tau}\max\left\{-\frac{d\p}{dt}(t),0\right\}\ dt & \mbox{if $\tau\le \bar t$} \cr
\int^{+\infty}_{\tau}\max\left\{\frac{d\p}{dt}(t),0\right\}\ dt & \mbox{if $\tau > \bar t$} \cr \end{array}
\right. .\]

Since $\int_{-\infty}^{+\infty}\frac{d\p}{dt}(t)\ dt=0$, it is immediate
to verify that
$\p_1(\tau)-\p_2(\tau)=\int_{-\infty}^{\tau}\frac{d\p}{dt}(t)\ dt=\p(\tau)$
for every $\tau\in\R$.
Because of the choice of $\bar t$, $\p_1$ and $\p_2$ are continuous also at $\bar t$. Moreover, we observe that both $\p_1$ and $\p_2$ are $C^1_c(\R)$ functions (here we are using the hypothesis $\frac{d\p}{dt}(\bar t)=0$). Furthermore they are
increasing in $(-\infty,\bar t\ ]$ and decreasing
in $[\bar t,+\infty)$.

By applying Proposition~\ref{simdefLambda} we get

\begin{eqnarray*}
&&\|\p\|=\|\p_1-\p_2\|\le\|\p_1\|+\|\p_2\|=\\
&=&\max\p_1\cdot\|\Lambda\|+\max\p_2\cdot\|\Lambda\|=\\
&=&(\max\p_1+\max\p_2)\cdot\|\Lambda\|=\\
&=&\left(\int_{-\infty}^{\bar t}\max\left\{\frac{d\p}{dt}(t),0\right\}\ dt+
\int^{+\infty}_{\bar t}\max\left\{\frac{d\p}{dt}(t),0\right\}\ dt\right)\cdot \|\Lambda\|=\\
&=&\left(\int_{-\infty}^{+\infty}\max\left\{\frac{d\p}{dt}(t),0\right\}\ dt\right)\cdot \|\Lambda\|=\\
&=&V^+_\p\cdot\|\Lambda\|=\frac{1}{2}V_\p\cdot \|\Lambda\|.
\end{eqnarray*}

Therefore the inequality $\|\p\|\le \frac{1}{2}V_\p\cdot
\|\Lambda\|$ is proved, in the case when $\frac{d\p}{dt}(\bar t)=0$.

Otherwise, if $\frac{d\p}{dt}(\bar t)\ne 0$, we observe that for
every $\eps>0$, $\p$ can be locally modified near $\bar t$ into a
function $\p_\eps\in C^1_c(\R)$ such that
\begin{itemize}
\item $V_{\p-\p_\eps}\le\eps$,
\item  $\frac{d\p_\eps}{dt}(\bar t)=0$,
\item   $ \int_{-\infty}^{\bar
t}\max\left\{\frac{d\p_\eps}{dt}(t),0\right\}\ dt= \int^{+\infty}_{\bar
t}\max\left\{-\frac{d\p_\eps}{dt}(t),0\right\}\ dt.$
\end{itemize}

The change we are using is represented in Figure \ref{modify}.




\begin {figure}
 \begin {center}

  \begin{picture}(0,0)%
\includegraphics{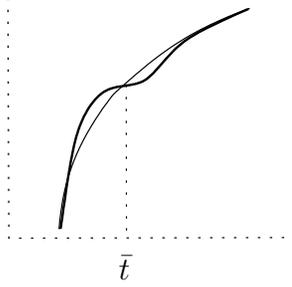}%
\end{picture}%
\setlength{\unitlength}{3947sp}%
\begingroup\makeatletter\ifx\SetFigFont\undefined%
\gdef\SetFigFont#1#2#3#4#5{%
  \reset@font\fontsize{#1}{#2pt}%
  \fontfamily{#3}\fontseries{#4}\fontshape{#5}%
  \selectfont}%
\fi\endgroup%
\begin{picture}(1827,1879)(2141,-3127)
\put(2849,-3050){\makebox(0,0)[lb]{\smash{{\SetFigFont{12}{14.4}{\rmdefault}{\mddefault}{\updefault}{\color[rgb]{0,0,0}$\bar t$}%
}}}}
\end{picture}%

  \caption {The change from $\p$ (thin) to $\p_\eps$ (thick) in order to get $\frac{d\p_\eps}{dt}(\bar t)= 0$ in the proof of the Bounding Lemma for $C^1_c(\R)$.}

\label{modify}
 \end {center}
\end {figure}

Because of what we have just proved in the case $\frac{d\p}{dt}(\bar t)=0$, it follows  that $\| \p_\eps\|\le
\frac{1}{2}V_{\p_\eps}\cdot \|\Lambda\|$ and hence $\| \p_\eps\|\le
\frac{1}{2}(V_{\p}+\eps)\cdot \|\Lambda\|$ (since $|V_{\p}-V_{\p_\eps}|\le V_{\p-\p_\eps}\le\eps$). Now, the Bounding
Lemma~\ref{bounds} assures that
$$\big|\|\p\|-\|\p_\eps\|\big|\le
\|\p-\p_\eps\|\le V_{\p-\p_\eps}\cdot\|S\|\le\eps\cdot\|S\|,$$ and hence
$$\|\p\|-\eps\cdot\|S\|\le\frac{1}{2}(V_{\p}+\eps)\cdot
\|\Lambda\|.$$

Then, the right inequality is proved for any $\p$ by passing to the
limit for $\eps$ tending to $0$.

\end{proof}

\begin{rem}
The double inequality that we have just proved shows that, if we confine ourselves to consider functions in $C^1_c(\R)$, half the total variation and $\max|\p|$
are the two extremal cases of RPI-norms. All other RPI-norms are somewhere between them, after normalization with respect to $\Lambda$.
We also observe that the two new inequalities are sharp, as we can immediately verify by setting $\|\p\|=\max|\p|$ and $\|\p\|=V_\p$.
\end{rem}

\begin{rem}
While the lower bound in the Bounding Lemma~\ref{bounds}
vanishes for all functions in $C^1_c(\R)$,
the lower bound in Lemma~\ref{boundsC1} never vanishes for non-zero functions in $C^1_c(\R)$.
This difference makes the study of functions with  compact support easier than the study of general almost sigmoidal functions.
\end{rem}

\subsection{Can a reparametrization invariant norm be induced by an inner product?}

We consider the question whether  a reparametrization invariant
norm can be associated with some inner
product. The next result shows that the answer to this question is negative.

\begin{prop}\label{scalprod}
 No inner product on $AS^1(\R)$
can induce a reparametrization invariant norm.
\end{prop}

\begin{proof}
Assume that an inner product
$\langle\cdot,\cdot\rangle$ exists on $AS^1(\R)$,  inducing a reparametrization invariant
norm.
 The
associated norm $\|\cdot\|$  satisfies the parallelogram identity:

\begin{eqnarray}
\label{parid}
\|\p_1+\p_2\|^2+\|\p_1-\p_2\|^2=2(\|\p_1\|^2+\|\p_2\|^2).
\end{eqnarray}

Let us take an almost sigmoidal function $\p$ with compact support, and set $\p_1=V^+\langle\p\rangle$ and $\p_2=V^-\langle\p\rangle$. Then,
$\p_1+\p_2=V\langle\p\rangle$, $\p_1-\p_2=\p$. By applying (\ref{parid}) and Proposition~\ref{simdef} about the norm of monotone almost sigmoidal functions we get
\begin{eqnarray*}
\|\p\|^2&=&2(\|V^+\langle\p\rangle\|^2+\|V^-\langle\p\rangle\|^2)-\|V\langle\p\rangle\|^2=\\
&=&
2(V^+_\p)^2\cdot\|S\|^2+2(V^-_\p)^2\cdot\|S\|^2-(V^+_\p+V^-_\p)^2\cdot\|S\|^2=\\
&=&
\left((V^+_\p)^2+(V^-_\p)^2-2(V^+_\p)(V^-_\p)\right)\cdot\|S\|^2=\\
&=&\left(V^+_\p-V^-_\p\right)^2\cdot\|S\|^2.
\end{eqnarray*}
Since for every $\p$ with compact support we have that
$V^+_\p-V^-_\p=0$, every such a function should have a vanishing norm. This contradicts the definition of norm.
\end{proof}

However, we remark that there exist degenerate symmetric bilinear
maps $\Phi$ inducing reparametrization invariant \emph{semi}-norms on
$AS^1(\R)$. An example is given by
$$\Phi(\p,\s)=\lim_{t\to +\infty}\p(t)\cdot
\lim_{t\to+\infty}\s(t).$$

\section{Standard reparametrization invariant norms}
\label{SRPIN}

In this section we  introduce  a class of reparametrization
invariant norms on $AS^1(\R)$. One well-known norm belonging to this
class is the $L_\infty$-norm. For the sake of conciseness and clearness in exposition, for every $\p\in AS^1(\R)$ we shall often use the symbol $\p^*$ to denote the function $\p^*(t)=\p(-t)$. Obviously, in general $\p^*$ is not an almost sigmoidal function, since it is obtained by composing $\p$ with an \emph{orientation-reversing} diffeomorphism of $\R$.

\subsection{The integral definition}

\begin{lem}\label{bound}
Let  $\p,\s\in AS^1(\R)$.  The following statements hold:
\begin{description}
  \item[i)] $$\int_{-\infty}^{+\infty}\p^*(t)\cdot\frac{d \s}{dt}(t)\ dt=-\int_{-\infty}^{+\infty}\frac{\ d \p^*}{dt}(t)\cdot\s(t)\ dt;$$
 \item[ii)] $$\left|\int_{-\infty}^{+\infty}\p^*(t)\cdot\frac{d \s}{dt}(t)\ dt\right|\le \max |\p|\cdot V_\s.$$
\end{description}
\end{lem}

\begin{proof}
\begin{description}
  \item[i)] Integrate by parts and observe that
  $\p^*(t)\cdot\s(t)|^{+\infty}_{-\infty}=0$.
 \item[ii)]  $$\left|\int_{-\infty}^{+\infty}\p^*(t)\cdot\frac{d \s}{dt}(t)\ dt\right|\le\max|\p|\cdot\int_{-\infty}^{+\infty}\left|\frac{d \s}{dt}(t)\right|\ dt = \max |\p|\cdot V_\s.$$
\end{description}
\end{proof}

\begin{thm} \label{thmRPInorms}
For every  function $\s\in AS^1(\R)-\{\mathbf{0}\}$ the setting
$$\|\p\|_{[\s]}=\sup_{ \hat \s\in [\s]}\left|\int_{-\infty}^{+\infty}\p(-t)\cdot\frac{d \hat \s}{dt}(t)\ dt\right|$$
defines a norm on the vector space $AS^1(\R)$, that  is invariant
under reparametrization. Moreover, if also
$\p\ne\mathbf{0}$, it holds that
\begin{enumerate}
\item $\|\p\|_{[\s]}=\|\p\|_{[-\s]}$;
\item $\|\p\|_{[\s]}=\|\s\|_{[\p]}$;
\item $\|\p\|_{[\s]}\le \min\left\{\max |\p|\cdot V_\s,\max |\s| \cdot
V_\p \right\}$.
\end{enumerate}
\end{thm}

\noindent \emph{Note.} In the rest of the paper the equality $\|\p\|_{[\s]}=\|\s\|_{[\p]}$ will be called \emph{exchange property}.
\begin{proof}

First of all, let us prove that $\|\cdot\|_{[\s]}$ is a norm.
Clearly $\|\p\|_{[{{\s}}]}$ is a non-negative real number. Indeed,
the finiteness of the $\sup$ follows from
Lemma~\ref{bound}\textbf{ii)} and the invariance of the total variation under reparametrizations.
Moreover it holds $\|\lambda
\p\|_{[{{\s}}]}=|\lambda|\cdot\|\p\|_{[{{\s}}]}$ for any
$\lambda\in \R$, and the triangle inequality is easily verified.
Also, if $\p\equiv 0$, then obviously $\|\p\|_{[{{\s}}]}=0$.
Therefore, the only thing we have to prove is that
$\|\p\|_{[{{\s}}]}=0$ implies  $\p\equiv 0$. We prove this
statement by contradiction, by assuming that $\|\p\|_{[{{\s}}]}=0$
and $\p^*(t_0)\ne 0$ for some $t_0$.
Let $[a,b]$ be a bounded closed interval $(a\neq b)$ containing
the compact support of $\frac{\ d\p^*}{dt}$ so that
$\p^*(t)=0$
for $t\ge b$. Analogously, define $[\alpha,\beta]$ to be a bounded
closed interval $(\alpha\neq \beta)$ containing the compact
support of $\frac{d\s}{dt}$, so that $\s(t)=0$ for $t\le
\alpha$. Since $\s$ is not constant, a point $t_1$ exists with
$\s(t_1)\neq 0$. It is easy to see that for every $\eps>0$ an
orientation-preserving $C^1$-diffeomorphism $ h_\eps:\R\to\R$ exists
such that
\begin{description}
  \item[i)] $ h_\eps([t_0-\eps,t_0+\eps])=[\alpha,t_1]$;
  \item[ii)] $\frac{d h_\eps}{dt}(t)=\eps$ for $t\ge t_0+\eps$.
\end{description}

\begin {figure}
 \begin {center}

  \includegraphics [height=5cm] {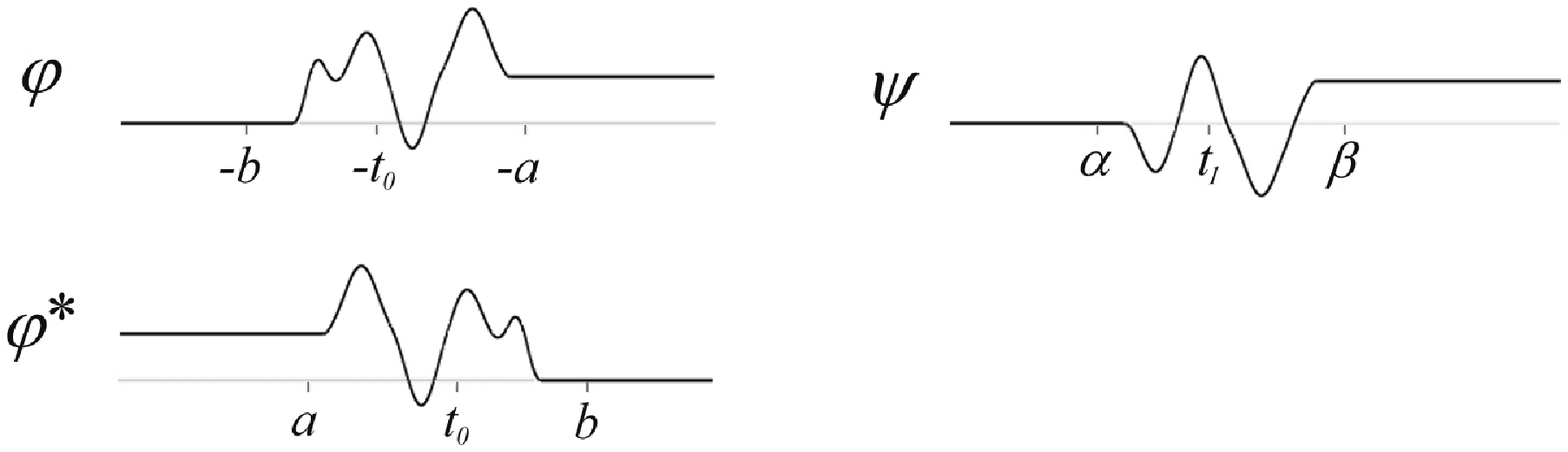}
  \caption {The functions $\p$, $\p^*$ and $\s$ (proof of Theorem \ref{thmRPInorms}).}

\label{dimstand}
 \end {center}
\end {figure}

We define $\hat \s=\s\circ h_\eps$. Obviously, $\hat \s\in AS^1(\R)$
and $\hat \s\sim \s$. Note that $\left|\frac{d \hat
\s}{dt}(t)\right|\le \eps\cdot\max\left|\frac{d \s}{dt}\right|$
for $t\ge t_0+\eps$, and that $\frac{d\hat\s}{dt}(t)=0$ for $t\le
t_0-\eps$ (since $\frac{d\s}{dt}(t)=0$ for $t\le
\alpha$). Taking $\eps$ small enough and remembering that
$\p^*(t)=0$ for $t\ge b$, we easily get

\begin{eqnarray*}
\label{diseq}
&&\left|\int_{-\infty}^{+\infty}\p^*(t)\cdot\frac{d \hat \s}{dt}(t)\
dt\right|=\left|\int_{t_0-\eps}^{b}\p^*(t)\cdot\frac{d \hat \s}{dt}(t)\
dt\right|\ge \\
&&\left|\int_{t_0-\eps}^{t_0+\eps}\p^*(t)\cdot\frac{d \hat
\s}{dt}(t)\ dt\right|-\max|\p|\cdot \eps\cdot\max\left|\frac{d
\s}{dt}\right|\cdot (b-(t_0+\eps)).\end{eqnarray*}

Now we observe that

$$\left|\int_{t_0-\eps}^{t_0+\eps}\p^*(t)\cdot\frac{d \hat
\s}{dt}(t)\ dt-\int_{t_0-\eps}^{t_0+\eps}\p^*(t_0)\cdot\frac{d \hat
\s}{dt}(t)\ dt\right|\le \max_{[t_0-\eps,t_0+\eps]}|\p^*(t)-\p^*(t_0)|\cdot V_\s.$$

Since we know that $\int_{t_0-\eps}^{t_0+\eps}\p^*(t_0)\cdot\frac{d \hat
\s}{dt}(t)\ dt=\p^*(t_0)\cdot\hat\s(t_0+\eps)=\p^*(t_0)\cdot\s(t_1)\neq 0$,
it follows that
$\lim_{\eps\to 0^+}\int_{t_0-\eps}^{t_0+\eps}\p^*(t)\cdot\frac{d \hat
\s}{dt}(t)\ dt=\p^*(t_0)\cdot\s(t_1)\neq 0$.
Therefore the value $\left|\int_{-\infty}^{+\infty}\p^*(t)\cdot\frac{d \hat
\s}{dt}(t)\ dt\right|$ is positive, if we have chosen a small enough $\eps$ in the costruction of $h_\eps$.
Hence
$\|\p\|_{[\s]}>0$, against our assumption. So, we have proved that
$\|\cdot\|_{[\s]}$ is a norm.

Now, let us consider a diffeomorphism $h\in D^1_+(\R)$ and prove
that $\|\p\circ h\|_{[\s]}=\|\p\|_{[\s]}$, i.e. $\|\cdot\|_{[\s]}$
is invariant under reparametrization. Setting $\tau=-h(-t)$, and
$\hat h(\tau)=-h^{-1}(-\tau)=t$, since $\hat h\in D^1_+(\R)$ we obtain that
$$\|\p\circ h\|_{[\s]}=\sup_{ \hat \s\in [\s]}\left|\int_{-\infty}^{+\infty}(\p\circ h)(-t)\cdot\frac{d \hat \s}{dt}(t)\ dt\right|=$$
$$=\sup_{ \hat \s\in [\s]}\left|\int_{-\infty}^{+\infty}\p^*(\tau)\cdot\frac{d \hat \s}{d\tau}(\hat h(\tau))
\cdot\frac{\ d \hat h}{d\tau}(\tau)\ d\tau\right|=$$
$$=\sup_{ \hat \s\in [\s]}\left|\int_{-\infty}^{+\infty}\p^*(\tau)\cdot\frac{d (\hat \s\circ \hat h)}{d\tau}(\tau)
\ d\tau\right|=$$ $$=\sup_{ \hat \s\in
[\s]}\left|\int_{-\infty}^{+\infty}\p^*(\tau)\cdot\frac{d \hat
\s}{d\tau}(\tau) \ d\tau\right|=\|\p\|_{[\s]}.$$

The equality $\|\p\|_{[\s]}=\|\p\|_{[-\s]}$ is trivial.

The equality $\|\p\|_{[\s]}=\|\s\|_{[\p]}$ follows from
Lemma~\ref{bound}\textbf{i)}, by observing that
$$\|\p\|_{[\s]}=\sup_{ h\in
D^1_+(\R)}\left|\int_{-\infty}^{+\infty}\p^*(t)\cdot\frac{d (\s\circ h)
}{dt}(t)\ dt\right|=$$
$$=\sup_{ h\in
D^1_+(\R)}\left|\int_{-\infty}^{+\infty} \frac{\ d\p^*}{dt}(t)
\cdot\s(h(t))\ dt\right|=\sup_{ h\in
D^1_+(\R)}\left|\int_{-\infty}^{+\infty} \frac{d\p}{dt}(t)
\cdot\s(h(-t))\ dt\right|=$$
$$=\sup_{ \hat h\in
D^1_+(\R)}\left|\int_{-\infty}^{+\infty} \frac{d\p}{d\tau}(\hat h (\tau))
\cdot\s^*(\tau)\cdot\frac{d\hat h}{d\tau}(\tau)\ d\tau\right|=$$
$$=\sup_{ \hat h\in
D^1_+(\R)}\left|\int_{-\infty}^{+\infty} \frac{d(\p\circ \hat h)}{dt}(\tau)
\cdot\s^*(\tau)\ d\tau\right|=\|\s\|_{[\p]},$$
where, once again, $\tau=-h(-t)$ and
$t=\hat h(\tau)=-h^{-1}(-\tau)$.

The inequality $\|\p\|_{[\s]}\le \min\left\{\max |\p|\cdot V_\s,\max |\s| \cdot V_\p \right\}$
follows from Lem\-ma~\ref{bound}\textbf{ii)} and the equality
$\|\p\|_{[\s]}=\|\s\|_{[\p]}$.
\end{proof}

In what follows, the norms $\|\cdot\|_{[\s]}$ will be called
\emph{standard reparametrization invariant norms} (or
\emph{standard RPI-norms}).

\subsection{Two examples of standard RPI-norms}
\label{twoex}

A simple instance of standard RPI-norm is given by the
$L_\infty$-norm, as the following proposition states.

\begin{prop}
\label{propS}
$\|\p\|_{[S]}=\max|\p|$.
\end{prop}

\begin{proof}
When $\p=\mathbf{0}$ the claim is trivial, so let us assume
$\max|\p|\ne 0$. By Theorem~\ref{thmRPInorms},
$\|\p\|_{[S]}\le\max|\p|$. Let now $t_{\max}$ be a value for which
$|\p^*(t_{\max})|=\max|\p|$ and consider, for $\eps>0$, the function
${S}_\eps(t)={S}\left(\frac{t-t_{\max}}{\eps }\right)$. Obviously,
${S}_\eps\sim{S}$. We observe that $\frac{d {S}_\eps }{dt}$ vanishes
outside the interval $[t_{\max}-\eps,t_{\max}+\eps]$ and
$\int_{t_{\max}-\eps}^{t_{\max}+\eps}\frac{d {S}_\eps }{dt}(t)\ dt=
S_\eps(t_{\max}+\eps)-S_\eps(t_{\max}-\eps)=1$. For $\eps$
sufficiently small, we have that $\p^*(t)$ has constant non-zero sign in
$[t_{\max}-\eps,t_{\max}+\eps]$. So, it follows that

\begin{eqnarray*}
& &\max|\p|-\left|\int_{-\infty}^{+\infty}\p^*(t)\cdot\frac{d {S}_\eps }{dt}(t)\ dt\right|=\\
&=&\max|\p|-\left|\int_{t_{\max}-\eps}^{t_{\max}+\eps}\p^*(t)\cdot\frac{d {S}_\eps }{dt}(t)\ dt\right|=\\
&=&\int_{t_{\max}-\eps}^{t_{\max}+\eps}\max|\p|\cdot \frac{d {S}_\eps }{dt}(t)\ dt-\int_{t_{\max}-\eps}^{t_{\max}+\eps}|\p^*(t)|\cdot\frac{d {S}_\eps }{dt}(t)\ dt=\\
&=&\int_{t_{\max}-\eps}^{t_{\max}+\eps}(\max|\p|-|\p^*(t)|)\cdot\frac{d {S}_\eps }{dt}(t)\ dt\le\\
&\le&\int_{t_{\max}-\eps}^{t_{\max}+\eps}(\max|\p|-\min_{|t-t_{\max}|\le \eps}|\p^*(t)|)\cdot\frac{d {S}_\eps }{dt}(t)\ dt=\\
&=&(\max|\p|-\min_{|t-t_{\max}|\le \eps}|\p^*(t)|)\cdot \int_{t_{\max}-\eps}^{t_{\max}+\eps}\frac{d {S}_\eps }{dt}(t)\ dt=\\
&=&\max|\p|-\min_{|t-t_{\max}|\le \eps}|\p^*(t)|.
\end{eqnarray*}
The continuity of $\p$ implies that $\lim_{\eps\to 0^+}\min_{|t-t_{\max}|\le \eps}|\p^*(t)|=\p^*(t_{\max})=\max|\p|$.
Hence  the equality $\lim_{\eps\to 0^+}\left|\int_{-\infty}^{+\infty}\p^*(t)\cdot\frac{d
{S}_\eps }{dt}(t)\
dt\right|=\max|\p|$ holds. Since we have already seen that $\|\p\|_{[S]}\le\max|\p|$, this proves that $\|\p\|_{[S]}=\max|\p|$.
\end{proof}

Another simple standard RPI-norm on $AS^1(\R)$ is given by $\max\p-\min\p$, as the
following proposition states.

\begin {figure}
 \begin {center}

  \includegraphics [height=10cm] {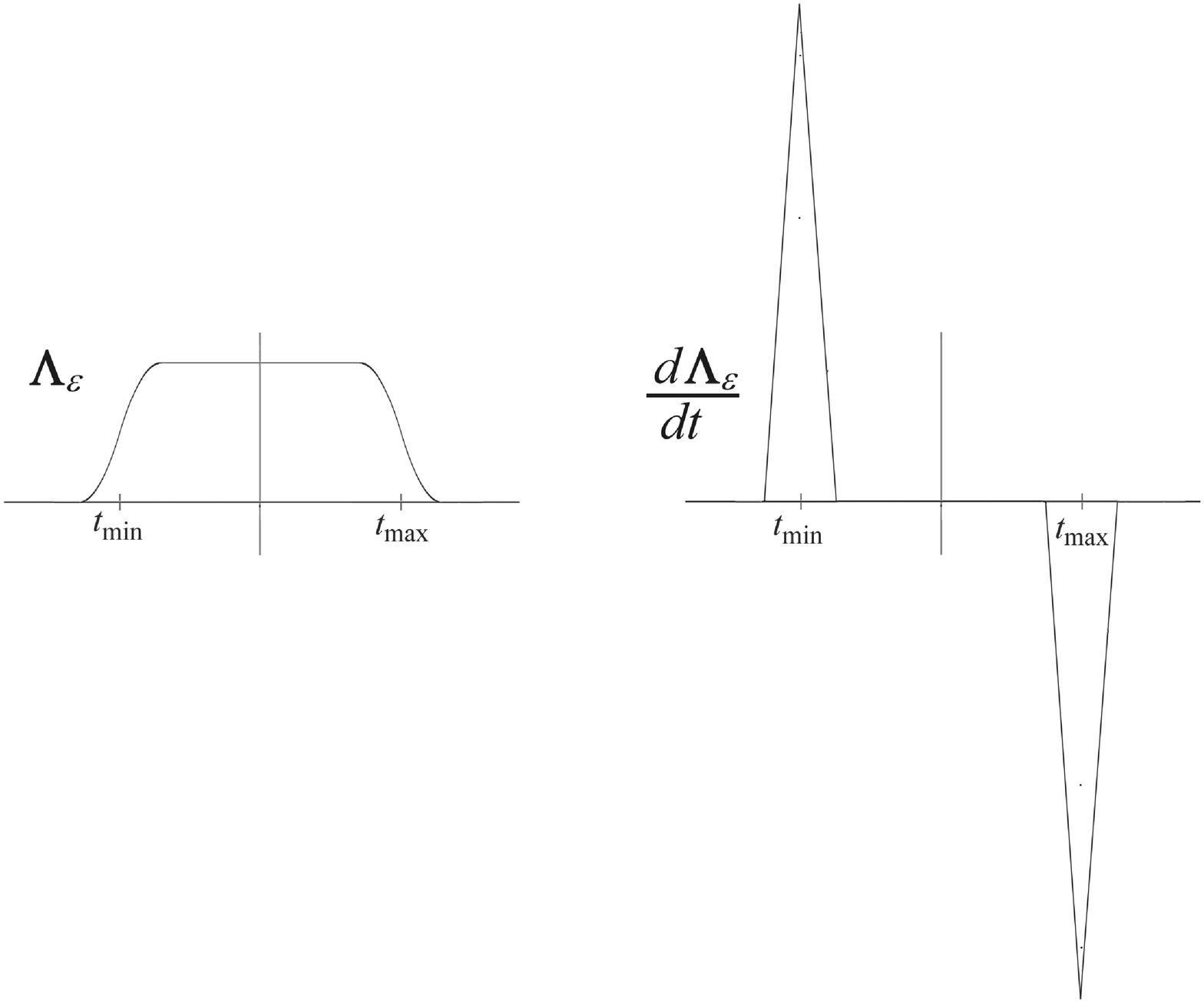}
  \caption {The function $\Lambda_\eps$ used in the proof of Proposition~\ref{maxmin} and its derivative (case $t_{\min}< t_{\max}$).}

\label{Leps}
 \end {center}
\end {figure}

\begin{prop}
\label{maxmin}
$\|\p\|_{[\Lambda]}=\max\p-\min\p$.
\end{prop}

\begin{proof}
Let us take a $C^1$-diffeomorphism $h\in D^1_+(\R)$. Possibly by
substituting $\p$ with $-\p$ we can assume that
$\int_{-\infty}^{+\infty}\p^*(t)\cdot\frac{d \left({\Lambda}\circ
h\right)}{dt}(t)\ dt\ge 0$. Hence, by recalling that $\Lambda$ is
increasing in $(-\infty,0]$ and decreasing in $[0,-\infty)$,
we obtain that

\begin{eqnarray*}
& &\left|\int_{-\infty}^{+\infty}\p^*(t)\cdot\frac{d \left({\Lambda} \circ h\right)}{dt}(t)\ dt\right|=\\
&=&\int_{-\infty}^{+\infty}\p^*(t)\cdot\frac{d \left({\Lambda} \circ h\right)}{dt}(t)\ dt=\\
&=&\int_{-\infty}^{h^{-1}(0)}\p^*(t)\cdot\frac{d \left({\Lambda} \circ h\right)}{dt}(t)\ dt+\int_{{h^{-1}(0)}}^{+\infty}\p^*(t)\cdot\frac{d \left({\Lambda} \circ h\right)}{dt}(t)\ dt\le\\
&\le&\int_{-\infty}^{{h^{-1}(0)}}\max\p\cdot\frac{d \left({\Lambda} \circ h\right)}{dt}(t)\ dt+\int_{{h^{-1}(0)}}^{+\infty}\min\p\cdot\frac{d \left({\Lambda} \circ h\right)}{dt}(t)\ dt=\\
&=&\max\p\cdot\left({\Lambda} (0)-\lim_{t\to -\infty}{\Lambda} (t)\right) +\min\p\cdot\left(\lim_{t\to +\infty}{\Lambda} (t)-{\Lambda} (0)\right) =\\
&=&\max\p-\min\p.
\end{eqnarray*}

Since $\|\p\|_{[\Lambda]}=\|-\p\|_{[\Lambda]}$, it follows that $\|\p\|_{[\Lambda]}\le\max\p-\min\p$.

Let $t_{\min}$ and $t_{\max}$ be a minimum point and a maximum point for $\p$, respectively.
 If $t_{\min}=t_{\max}$ then $\p\equiv 0$ and our statement is trivial. So, let us assume that $t_{\min}\neq t_{\max}$ .
Let us define $t_0=\min \left\{t_{\min},t_{\max}\right\}$ and $t_1=\max \left\{t_{\min},t_{\max}\right\}$. We consider the function $\Lambda_\eps=S\left(\frac{t-t_{0}}{\eps}\right)-
S\left(\frac{t-t_{1}}{\eps}\right)$ (see Figure~\ref{Leps}). Even if $\Lambda_\eps$ is not equivalent to $\Lambda$ we have that $\|\p\|_{[\Lambda]}=\|\p\|_{[\Lambda_\eps]}$ for $\eps$ small enough. Indeed,  from Theorem~\ref{thmRPInorms} (exchange property) and Proposition~\ref{simdefLambda} it follows that $\|\p\|_{[\Lambda]}=\|\Lambda\|_{[\p]}=\|\Lambda_\eps\|_{[\p]}=\|\p\|_{[\Lambda_\eps]}$.
Since it is easy to verify  that the equality $\lim_{\eps\to 0^+}\left|\int_{-\infty}^{+\infty}\p^*(t)\frac{d
{\Lambda}_\eps }{dt}(t)\
dt\right|=\max\p-\min \p$ holds, this implies that
$\|\p\|_{[\Lambda]}\ge\max\p-\min \p$. Therefore
our statement follows.
\end{proof}

\subsubsection{The key idea in using standard RPI-norms}

The two examples seen in the previous section show that, in some sense, computing standard RPI-norms is equivalent to computing the absolute value of a linear combination of Dirac deltas, applied to the function $\p^*(t)$. Indeed, it is easy to verify that in order to get $\|\p\|_{[S]}$ and $\|\p\|_{[\Lambda]}$ we have to compute
the values $\sup_{t}\left|\delta_{t}(\p^*)\right|$ and $\sup_{t_0\le t_1}\left|\delta_{t_0}(\p^*)-\delta_{t_1}(\p^*)\right|$, where
$\delta_{t}$ is the usual Dirac delta at point $t$. The ``weights'' of the Dirac deltas are determined by the integral
$\int_{-1}^{1}\frac{dS}{dt}(t)\ dt=1$ in the first case, and by
the integrals $\int_{-2}^{0}\frac{d\Lambda}{dt}(t)\ dt=1$, $\int_{0}^{2}\frac{d\Lambda}{dt}(t)\ dt=-1$ in the latter, i.e. the integrals of $\frac{d\psi}{dt}$ on the maximal intervals where the derivative of the function $\psi$ defining the norm $\|\cdot\|_{[\psi]}$ does not vanish. In order to compute $\|\p\|_{[S]}$ we place $\delta_{t}$ at a point where $|\p|$ takes its maximum value, while when we compute $\|\p\|_{[\Lambda]}$ we place $\delta_{t_0}$ and $\delta_{t_1}$ at the points  where $\p$ takes its maximum value and its minimum value (not necessarily in this order).
We shall carefully analyze and generalize this approach in Section~\ref{DROSRPIN}.

\subsection{Not every RPI-norm is a standard RPI-norm}
It is important to observe that some RPI-norms are not
standard RPI-norms.

In order to show this, now we give a useful property of standard RPI-norms.

\begin{prop}\label{comblin}
Let $\|\cdot \|$ be a RPI-norm. If it can be obtained as a finite
linear combination of standard RPI-norms with positive
coefficients, then $\|S\|\le \|\Lambda\|$.
\end{prop}

\begin{proof}
Let us assume that some functions $\s_1,\ldots,\s_k\in AS^1(\R)$
and a $k$-tuple $(a_1,\ldots,a_k)$ of positive numbers exist such
that for every $\p\in AS^1(\R)$ it holds that
$\|\p\|=\sum_{i=1}^{k}a_i\cdot\|\p\|_{[\s_i]}$. Let us consider
the functions $\tilde\s_1,\ldots,\tilde\s_k$, such that, for
$i=1,\ldots,k$, if $\max|\s_i|=\max\s_i$ then $\tilde\s_i=\s_i$,
otherwise  $\tilde\s_i=-\s_i$. By Theorem~\ref{thmRPInorms}, it
holds that $\|\p\|=\sum_{i=1}^{k}a_i\cdot\|\p\|_{[\tilde\s_i]}$.
Therefore, by the exchange property and Proposition~\ref{propS},
we have that
$$\|S\|=\sum_{i=1}^{k}a_i\cdot\|S\|_{[\tilde\s_i]}=\sum_{i=1}^{k}a_i\cdot\|\tilde\s_i\|_{[S]}=\sum_{i=1}^{k}a_i\cdot\max|\tilde\s_i|=\sum_{i=1}^{k}a_i\cdot\max\tilde\s_i.$$
Analogously, by the exchange property and
Proposition~\ref{maxmin}, we obtain that
$\|\Lambda\|=\sum_{i=1}^{k}a_i\cdot(\max\tilde\s_i-\min\tilde\s_i)$.
The claim immediately follows since
$\sum_{i=1}^{k}a_i\cdot\min\tilde\s_i\le 0$  (recall that $\min
\s\le 0$ for every $\psi\in AS^1(\R)$).
\end{proof}

As a consequence of this property, we can furnish an example of
$RPI$-norm that cannot be represented as a linear combination with
positive coefficients of standard RPI-norms.

\begin{cor}
The RPI-norm $\|\p\|=\max|\p|+\lim_{t\to
+\infty}|\p(t)|$ cannot be represented as a finite linear
combination with
positive coefficients of standard RPI-norms. In particular, it is not a standard RPI-norm.
\end{cor}

\begin{proof}
It is sufficient to observe that $\|S\|=\max|S|+\lim_{t\to
+\infty}|S(t)|=2$, $\|\Lambda\|=\max|\Lambda|+\lim_{t\to
+\infty}|\Lambda(t)|=1$, and apply Proposition~\ref{comblin}.
\end{proof}

\begin{rem}
\label{domanda}
It could be interesting to know if the norm $\max|\p|+\lim_{t\to
+\infty}|\p(t)|$ can be represented either as a $\sup$ or as an $\inf$
of a
suitable set of standard RPI-norms.
\end{rem}

Another example of RPI-norm that cannot be expressed as a finite linear
combination with
positive coefficients of standard RPI-norms is the total variation.

\begin{prop}
\label{novariationRPI} The total variation cannot be represented
as a finite linear combination with positive coefficients of
standard RPI-norms. In particular, it is not a standard RPI-norm.
\end{prop}

\begin{proof}
If the total variation could be represented as a linear combination
with
positive coefficients of standard RPI-norms, the equality $V_\p=\sum_{i=1}^k
a_i\|\p\|_{[\s_i]}$ would hold for every $\p\in AS^1(\mathbb{R})$,
when a suitable set $\{a_1,\ldots,a_k\}$ of positive coefficients is
chosen.

By Theorem~\ref{thmRPInorms}  we would have that $V_\p\le
\sum_{i=1}^k a_i V_{\s_i}\cdot \max |\p|$. This inequality
contradicts the fact that we can easily find a function $\bar\p\in
AS^1(\R)$ such that $\max \bar\p\neq 0$ and the ratio
$\frac{V_{\bar\p}}{\max\bar\p}$ is arbitrarily large.
\end{proof}

 Nevertheless, the total variation can be seen
as the $\sup$ of a suitable set of standard RPI-norms, as shown in the following
Section~\ref{TVsup}.

\subsubsection{The total variation is a sup of standard
RPI-norms}\label{TVsup}

We have seen in Proposition~\ref{novariationRPI} that
 the total variation is not a standard RPI-norm. We now show that it
is the $\sup$ of a family of standard RPI-norms.

\begin{prop}\label{totalvariationsup}
 For every $n\ge 1$, let us set $L_n(t)=\sum_{i=0}^{n-1}(-1)^i\Lambda(t-4i)$.
 Then, for every $\p\in
AS^1(\R)$, we have that $V_\p=\sup_{n\in \N^+}\|\p\|_{[L_{n}]}=\lim_{n\to\infty}\|\p\|_{[L_n]}$.

\end{prop}

\begin{proof}
Let us prove the first equality.
By applying
Theorem~\ref{thmRPInorms} we obtain that
$\|\p\|_{[L_n]}\le\max|L_n|\cdot V_\p=V_\p$.
We just have to show
that for every $\eps>0$ an $n$ exists such that
$V_\p-\|\p\|_{[L_{n}]}\le\eps$. This is trivially true if
$\p=\mathbf{0}$ so let us assume $\p\ne\mathbf{0}$. Let $[a,b]$ be
a closed interval containing the support of
$\frac{\ d\p^*}{dt}$. Let us recall that $V_\p=V_{\p^*}=\sup_k\sup_{a=t_0<t_{1}<\ldots
<t_{k}=b}\sum_{i=0}^{k-1}|\p^*(t_{i+1})-\p^*(t_{i})|$. Hence, there
exist  $n\ge 1$, and a partition  $a=\bar\tau_{0}<\bar\tau_{1}<\ldots <\bar\tau_{n}=b$
of $[a,b]$,  such that
$$0\le
V_\p-\sum_{i=0}^{n-1}|\p^*(\bar\tau_{i+1})-\p^*(-\bar\tau_i)|\le \frac{\eps}{2}.$$
Possibly by substituting our partition with a simpler one, we can assume that
$\p^*(\bar\tau_{i+1})-\p^*(\bar\tau_i)\ne 0$ for $i=0,\ldots ,n-1$
and, for $n\ge 2$, ${\sign}(\p^*(\bar\tau_{i+1})-\p^*(\bar\tau_i))\ne
{\sign}(\p^*(\bar\tau_{i+2})-\p^*(\bar\tau_{i+1}))$ for $i=0,\ldots
,n-2$.

\begin {figure}
 \begin {center}

  \includegraphics [height=5cm] {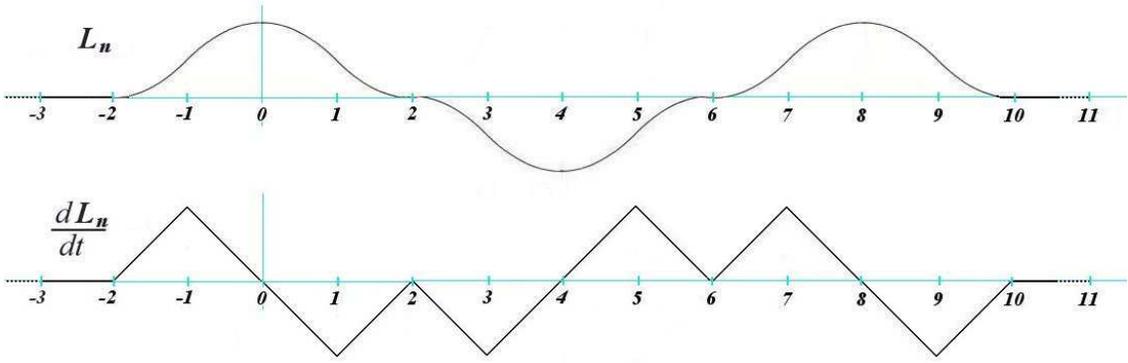}
  \caption {The function $L_n$  and its derivative (case $n=3$), used in Proposition~\ref{totalvariationsup}.}

\label{Ln3}
 \end {center}
\end {figure}

For every sufficiently small $\eta>0$, let us consider an
orientation-preserving $C^1$-diffeomorphism $h_\eta$ that takes the interval
$[\bar\tau_{i}-\eta,\bar\tau_{i}+\eta]$ onto the interval
$[4(i-1)+\eta,4i-\eta]$,  for every
integer $i$ with $0\le i\le n$. Let us observe that the function $L_n$ is monotone on every interval $[4(i-1)+\eta,4i-\eta]$ (cf. Figure \ref{Ln3}).

For $\eta$ small enough, by recalling that $\p^*(\bar\tau_{n})=\p^*(b)=0$, it is easy to prove that
 $$\left|\ \left|\p^*(\bar \tau_0)+\sum_{i=1}^{n-1}(-1)^i\cdot 2\cdot \p^*(\bar\tau_i)\right|-
\left|\int_{-\infty}^{+\infty}\p^*(\tau)\cdot\frac{d \left(L_n\circ h_\eta\right)}{d\tau}(\tau)\ d\tau\right|\ \right|\le \frac{\eps}{2}.$$

The assumptions about the differences $\p^*(\bar\tau_{i+1})-\p^*(\bar\tau_i)$ and, once again the condition
$\p^*(\bar\tau_{n})=\p^*(b)=0$, imply that
 $$\sum_{i=0}^{n-1}|\p^*(\bar\tau_{i+1})-\p^*(\bar\tau_{i})|=\left|\p^*(\bar \tau_0)+\sum_{i=1}^{n-1}(-1)^i\cdot 2\cdot \p^*(\bar\tau_i)\right|.$$

It follows that $$\left|V_\p -\left|\int_{-\infty}^{+\infty}\p^*(\tau)\cdot\frac{d \left(L_n\circ h_\eta\right)}{d\tau}(\tau)\ d\tau\right|\ \right|\le \eps.$$
Hence $V_\p-\|\p\|_{[L_{n}]}\le\eps$, and the first equality in our statement is proved.

The second equality follows from the first one, by observing that the sequence $(\|\p\|_{[L_n]})$ is increasing.
\end{proof}

\begin{rem}
In plain words, our proof of Proposition~\ref{totalvariationsup} is
based on recognizing that
$\sum_{i=0}^{n-1}|\p^*(\bar\tau_{i+1})-\p^*(\bar\tau_{i})|$ can be seen
as the value taken by the absolute value of the linear functional
$\delta_{\bar\tau_{0}}+\sum_{i=1}^{n-1}(-1)^i2\cdot
\delta_{\bar\tau_{i}}$ computed at $\p^*$, where $\delta_{t}$ is the
usual Dirac delta at point $t$. The reparametrization $h_\eta$
allows us to approximate
$\left|\delta_{\bar\tau_{0}}(\p^*)+\sum_{i=1}^{n-1}(-1)^i2\cdot
\delta_{\bar\tau_{i}}(\p^*)\right|$ by
$\left|\int_{-\infty}^{+\infty}\p^*(\tau)\cdot\frac{d \left(L_n\circ
h_\eta\right)}{d\tau}(\tau)\ d\tau\right|$, by ``concentrating'' at
$\bar\tau_0$ and at the other $\bar\tau_i$'s a signed variation of
$L_n\circ h_\eta$ approximately equal to $1$ and $\pm 2$,
respectively. This last idea will be developed in next
Section~\ref{DROSRPIN} by using general weights (not just $1$ and
$\pm 2$) for the Dirac deltas, and its generalization will lead to
the Representation Theorem~\ref{teoremadirappresentazione}. This
result assures that every standard RPI-norm can be seen as the
absolute value of a suitable linear combination of Dirac deltas,
maximized with respect to the movements that preserve the deltas'
position order. We could obtain Proposition~\ref{totalvariationsup}
as a consequence of the Representation Theorem, but we have
preferred to anticipate this result for the sake of clarity of
exposition. Moreover, this choice allows us to illustrate the ideas
that we are going to develop.
\end{rem}

\section{Discrete representation of standard RPI-norms}
\label{DROSRPIN}

In this section we show how to compute the standard RPI-norms in a simpler,
discrete way.

As we have seen in Section~\ref{twoex}, the standard RPI-norms $\max
|\p|$ and $\max\p-\min\p$ can be expressed as
$\sup_{t}\left|\delta_{t}(\p^*)\right|$ and $\sup_{t_1\le
t_2}\left|\delta_{t_1}(\p^*)-\delta_{t_2}(\p^*)\right|$ respectively,
where $\delta_{t}$ is the usual Dirac delta at point $t$.

This is a general property of standard RPI-norms, that they can be
seen as $\sup$ of the absolute value of linear combinations of Dirac
deltas. The basic idea underlying this fact is that the $\sup$ of $\left|\int_{-\infty}^{+\infty}\p^*(t)\cdot\frac{d(\s\circ
h)}{dt}(t)\ dt\right|$ for $h\in D^1_+(\R)$ is obtained  by  considering a sequence of reparametrizations more and more
concentrating the variation of $\s$ at suitable points. Passing from the integral definition of
standard RPI-norm to linear combinations of Dirac deltas, the sup with
respect to orientation-preserving reparametrizations is
substituted by a sup with respect to movements that shift  the Dirac
deltas' centers without changing their ordering (Theorem \ref{teoremadirappresentazionegenerico}). We shall see that
the best choice is to place these Dirac deltas' centers at critical
points of $\p$ (Representation Theorem \ref{teoremadirappresentazione}).

\begin {figure}
 \begin {center}
\begin{picture}(0,0)%
\includegraphics{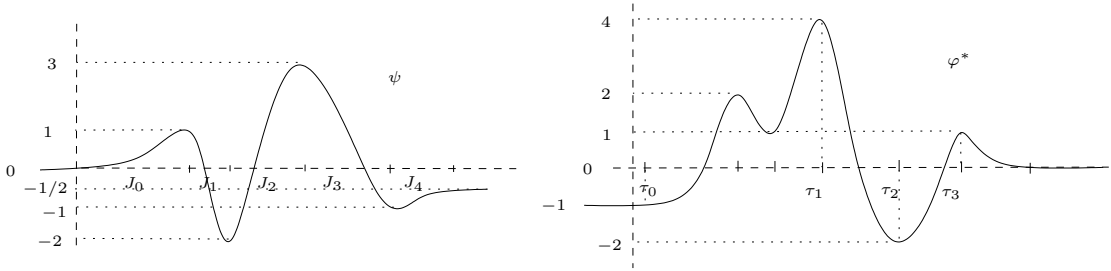}%
\end{picture}%
\setlength{\unitlength}{3158sp}%
\begingroup\makeatletter\ifx\SetFigFont\undefined%
\gdef\SetFigFont#1#2#3#4#5{%
  \reset@font\fontsize{#1}{#2pt}%
  \fontfamily{#3}\fontseries{#4}\fontshape{#5}%
  \selectfont}%
\fi\endgroup%
\begin{picture}(8685,2102)(24,-2695)
\put(963,-2056){\makebox(0,0)[lb]{\smash{{\SetFigFont{6}{7.2}{\rmdefault}{\mddefault}{\updefault}{\color[rgb]{0,0,0}$J_0$}%
}}}}
\put(1546,-2056){\makebox(0,0)[lb]{\smash{{\SetFigFont{6}{7.2}{\rmdefault}{\mddefault}{\updefault}{\color[rgb]{0,0,0}$J_1$}%
}}}}
\put(2519,-2045){\makebox(0,0)[lb]{\smash{{\SetFigFont{6}{7.2}{\rmdefault}{\mddefault}{\updefault}{\color[rgb]{0,0,0}$J_3$}%
}}}}
\put(3038,-1220){\makebox(0,0)[lb]{\smash{{\SetFigFont{6}{7.2}{\rmdefault}{\mddefault}{\updefault}{\color[rgb]{0,0,0}$\s$}%
}}}}
\put( 39,-1957){\makebox(0,0)[lb]{\smash{{\SetFigFont{6}{7.2}{\rmdefault}{\mddefault}{\updefault}{\color[rgb]{0,0,0}$0$}%
}}}}
\put(320,-2283){\makebox(0,0)[lb]{\smash{{\SetFigFont{6}{7.2}{\rmdefault}{\mddefault}{\updefault}{\color[rgb]{0,0,0}$-1$}%
}}}}
\put(283,-2511){\makebox(0,0)[lb]{\smash{{\SetFigFont{6}{7.2}{\rmdefault}{\mddefault}{\updefault}{\color[rgb]{0,0,0}$-2$}%
}}}}
\put(360,-1123){\makebox(0,0)[lb]{\smash{{\SetFigFont{6}{7.2}{\rmdefault}{\mddefault}{\updefault}{\color[rgb]{0,0,0}$3$}%
}}}}
\put(330,-1652){\makebox(0,0)[lb]{\smash{{\SetFigFont{6}{7.2}{\rmdefault}{\mddefault}{\updefault}{\color[rgb]{0,0,0}$1$}%
}}}}
\put(171,-2093){\makebox(0,0)[lb]{\smash{{\SetFigFont{6}{7.2}{\rmdefault}{\mddefault}{\updefault}{\color[rgb]{0,0,0}$-1/2$}%
}}}}
\put(2005,-2055){\makebox(0,0)[lb]{\smash{{\SetFigFont{6}{7.2}{\rmdefault}{\mddefault}{\updefault}{\color[rgb]{0,0,0}$J_2$}%
}}}}
\put(3155,-2050){\makebox(0,0)[lb]{\smash{{\SetFigFont{6}{7.2}{\rmdefault}{\mddefault}{\updefault}{\color[rgb]{0,0,0}$J_4$}%
}}}}
\put(4678,-2541){\makebox(0,0)[lb]{\smash{{\SetFigFont{6}{7.2}{\rmdefault}{\mddefault}{\updefault}{\color[rgb]{0,0,0}$-2$}%
}}}}
\put(4565,-1949){\makebox(0,0)[lb]{\smash{{\SetFigFont{6}{7.2}{\rmdefault}{\mddefault}{\updefault}{\color[rgb]{0,0,0}$0$}%
}}}}
\put(4718,-1677){\makebox(0,0)[lb]{\smash{{\SetFigFont{6}{7.2}{\rmdefault}{\mddefault}{\updefault}{\color[rgb]{0,0,0}$1$}%
}}}}
\put(7426,-1089){\makebox(0,0)[lb]{\smash{{\SetFigFont{6}{7.2}{\rmdefault}{\mddefault}{\updefault}{\color[rgb]{0,0,0}$\p^*$}%
}}}}
\put(4709,-1357){\makebox(0,0)[lb]{\smash{{\SetFigFont{6}{7.2}{\rmdefault}{\mddefault}{\updefault}{\color[rgb]{0,0,0}$2$}%
}}}}
\put(4709,-788){\makebox(0,0)[lb]{\smash{{\SetFigFont{6}{7.2}{\rmdefault}{\mddefault}{\updefault}{\color[rgb]{0,0,0}$4$}%
}}}}
\put(4239,-2230){\makebox(0,0)[lb]{\smash{{\SetFigFont{6}{7.2}{\rmdefault}{\mddefault}{\updefault}{\color[rgb]{0,0,0}$-1$}%
}}}}
\put(5004,-2089){\makebox(0,0)[lb]{\smash{{\SetFigFont{6}{7.2}{\rmdefault}{\mddefault}{\updefault}{\color[rgb]{0,0,0}$\tau_0$}%
}}}}
\put(6905,-2106){\makebox(0,0)[lb]{\smash{{\SetFigFont{6}{7.2}{\rmdefault}{\mddefault}{\updefault}{\color[rgb]{0,0,0}$\tau_2$}%
}}}}
\put(7382,-2118){\makebox(0,0)[lb]{\smash{{\SetFigFont{6}{7.2}{\rmdefault}{\mddefault}{\updefault}{\color[rgb]{0,0,0}$\tau_3$}%
}}}}
\put(6311,-2110){\makebox(0,0)[lb]{\smash{{\SetFigFont{6}{7.2}{\rmdefault}{\mddefault}{\updefault}{\color[rgb]{0,0,0}$\tau_1$}%
}}}}
\end{picture}%
\caption{The functions used in the example described in
Section~\ref{DROSRPIN}.} \label{esdiscr}
\end {center}
\end {figure}

By way of exemplification, let us consider the norm
$\|\p\|_{[\s]}$ where $\s(t)$ and $\p^*(t)$ are the functions
illustrated in Figure~\ref{esdiscr}. Let us denote by
$J_i=(a_i,b_i)$ the maximal open intervals where the derivative of
$\s$ has constant non-zero sign. It holds that
$\s(b_0)-\s(a_0)=1$, $\s(b_1)-\s(a_1)=-3$, $\s(b_2)-\s(a_2)=5$,
$\s(b_3)-\s(a_3)=-4$, $\s(b_4)-\s(a_4)=1/2$. In order to increase
the value
$\left|\int_{-\infty}^{\infty}\p^*(t)\cdot\frac{d(\s\circ
h)}{dt}(t)\ dt\right|$, it is convenient to take
reparametrizations $h_\eta$ that transform   smaller and smaller
neighbourhoods of suitable  points $t_0\le t_1\le t_2\le t_3\le
t_4$ ordinately to the intervals
$(a_0+\eta,b_0-\eta),(a_1+\eta,b_1-\eta),\ldots,(a_4+\eta,b_4-\eta)$,
with a smaller and smaller $\eta>0$ (or to $(a_0+\eta,b_1-\eta)$
if, say, $t_0=t_1$, and so on). By passing to the limit one
obtains that
$$\sup_h\left|\int_{-\infty}^{\infty}\p^*(t)\cdot\frac{d(\s\circ h)}{dt}(t)\ dt\right|=\sup_{t_0\le t_1\le t_2\le t_3\le t_4}\left|1\cdot\p^*(t_0)-3\cdot\p^*(t_1)+5\cdot\p^*(t_2)-4\cdot\p^*(t_3)+\frac{1}{2}\cdot\p^*(t_4)\right|.$$ In other words,
$$\sup_h\left|\int_{-\infty}^{\infty}\p^*(t)\cdot\frac{d(\s\circ h)}{dt}(t)\ dt\right|=\sup_{t_0\le t_1\le t_2\le t_3\le t_4}\left|\left(1\cdot\delta_{t_0}-3\cdot\delta_{t_1}+5\cdot\delta_{t_2}-4\cdot\delta_{t_3}+\frac{1}{2}\cdot\delta_{t_4}\right)(\p^*)\right|.$$
Now, one easily sees that, in order to get the greatest value, the $t_i$'s must be critical points of $\p^*$. In particular, in this case the $\sup$ is attained when $t_0=t_1=\tau_0$,  $t_2=\tau_1$, $t_3=\tau_2$, $t_4=\tau_3$, so that $\|\p\|_{[\s]}=30.5$.

We point out that the considered linear combinations can involve
infinitely many terms. This is the main difficulty to manage in
this section, and will require some computations.

The key result obtained in this section (Representation Theorem \ref{teoremadirappresentazione}) will be
fundamental in the next section, where we shall use it to prove that all the
standard RPI-norms of a piecewise monotone $C^1$-function $\p$ with compact support are sufficient to reconstruct $\p$, up to reparametrization and an arbitrarily small error
with respect to the total variation norm.

 The first step to get these results is defining a bilinear function $F$
that will be useful in the sequel.

\subsection{The functional $F$}
\label{functF}

\begin {figure}
 \begin {center}

  \includegraphics [height=4cm] {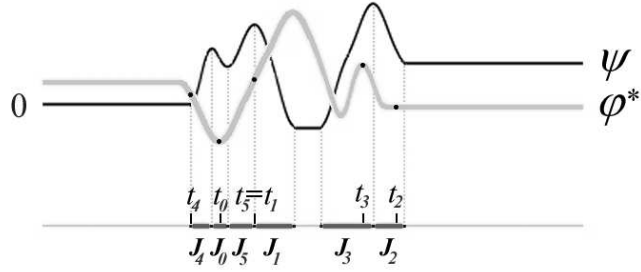}
  \caption {The set $\mathcal{J}(\s)=\left\{J_i\right\}_{i\in I}$ for the displayed function $\s$. A possible choice of the points $t_i$ is shown, with respect to $\p$ (the corresponding points on the graph of $\p^*$ are marked). Observe that the order $\preceq$ in $\mathcal{J}(\s)$ does not need to be given by the index $i$, and that some $t_i$'s belong to the boundary of the corresponding interval $J_i$.}

\label{tconindici}
 \end {center}
\end {figure}

All along the remainder of Section~\ref{DROSRPIN} we shall assume that two functions $\p,\s\in AS^1(\R)$ are given, with $\s\ne {\mathbf 0}$. Let us consider the countable set $\mathcal{J}(\s)$ of
all maximal open intervals of $\R$ where $\frac{d\s}{dt}$ does not
vanish. We shall set $\mathcal{J}(\s)=\left\{J_i\right\}_{i\in I}$, where $I$ is either the finite set $\{0,\ldots,n-1\}$ or the set $\N$.
For each $J_i=(a_i,b_i)\in \mathcal{J}(\s)$ we shall assume that a point $t_i\in [a_i,b_i]$ is chosen, such that
\begin{enumerate}
  \item  $t_i$ is a point where the restriction of $\p^*$ to the closed interval $[a_i,b_i]$ takes its maximum value, if $\s|_{a_i}^{b_i}>0$;
  \item  $t_i$ is a point where the restriction of $\p^*$ to the closed interval $[a_i,b_i]$ takes its minimum value, if $\s|_{a_i}^{b_i}<0$.
\end{enumerate}
Here we set $\s|_{a_i}^{b_i}=\s(b_i)-\s(a_i)$ (see example in Figure~\ref{tconindici}).

\begin{defn}\label{bp}
The set $\{t_i\}_{i\in I}$ is said to be a
\emph{set of basepoints for the pair $(\p,\s)$}.\end{defn}

On the set $\mathcal{J}(\s)$ we shall consider the order $\preceq$ induced by the $a_i$'s. In other words we shall set $J_i\preceq J_j$ if and only if $a_i\le a_j$. This order will not need to coincide with the order induced by the index $i$. The symbol $\sigma_i$  will denote
the sign of $\frac{d\s}{dt}$ in the interval $J_i$, i.e. $\sigma_i=\sign \ \s|_{a_i}^{b_i}$.

\begin{defn}\label{F}
We define the
bilinear functional $F: AS^1(\R)\times \left(AS^1(\R)-{\mathbf 0}\right)\rightarrow\R$ by setting $$F(\p,\s)=\sum_{i\in I} \s|_{a_i}^{b_i}\cdot\p^*(t_i)$$
where $\{t_i\}_{i\in I}$ is a set of
basepoints  for $(\p,\s)$.
\end{defn}

In other words, $F(\p,\s)= \sum_{i\in I} \s|_{a_i}^{b_i}\cdot
\delta_{t_i}(\p^*)$, where
$\delta_{t_i}$ is the usual Dirac delta at point $t_i$. Note
that the definition of $F(\p,\s)$ does not depend on the particular choice
of the set of basepoints for the pair $(\p,\s)$. The idea underlying
the definition of basepoint is to maximize each addend
$\s|_{a_i}^{b_i}\cdot\p^*(t)$ in the definition of $F$, when $t$ varies in $[a_i,b_i]$.

We observe that $\left| F(\p,\s)\right|\le \max
|\p|\cdot V_\s <+\infty$, since $\s$ has  bounded variation.

\begin{rem}
It is easy to verify that $\{t_i\}_{i\in I}$ is a
set of basepoints for $(\p,\s)$ if and only if for every $i\in I$
$$\p^*(t_i)=\sigma_i\cdot\max_{[a_i,b_i]}\left\{\sigma_i\cdot
\p^*\right\}.$$
Therefore, we
get this equivalent definition for $F$:
$$F(\p,\s)=\sum_{i\in I} \sigma_i\cdot\s|_{a_i}^{b_i}\cdot\max_{[a_i,b_i]}\left\{\sigma_i\cdot
\p^*\right\}.$$
Moreover, we observe that each set of basepoints for $(\p,\s)$ is
contained in the compact support of $\frac{d\s}{dt}$.
\end{rem}

\subsection{Some useful properties of the functional $F$}
Let us consider the set $H_{\s}$  of all orientation-preserving $C^1$-diffeomorphism of
the real line that take each interval $J_i\in \mathcal{J}(\s)$ to
itself.
The following lemma shows the key property of the
functional $F$.

\begin{lem}\label{interpretationlemma}
$\sup_{h\in H_\s}\left|\int_{-\infty}^{+\infty}\p^*(t)\cdot\frac{d
(\s\circ h)}{dt}(t)\ dt\right|=\max \left\{F(\p,\s),
F(\p,-\s)\right\}$.
\end{lem}

\begin{proof}
\textbf{PART 1: } $\sup_{h\in
H_\s}\left|\int_{-\infty}^{+\infty}\p^*(t)\cdot\frac{d (\s\circ
h)}{dt}(t)\ dt\right|\le \max \left\{F(\p,\s), F(\p,-\s)\right\}$.

Let $\bar \s=\s\circ h$ with $h\in H_\s$. Consider the countable
set $\mathcal{J}(\bar \s)$ of all maximal open intervals of $\R$
where $\frac{d\bar\s}{dt}$ does not vanish. Obviously,
$\mathcal{J}(\bar \s)=\mathcal{J}(\s)$ and ${\sign} \frac{d\bar
\s}{dt}={\sign} \frac{d \s}{dt}$.
Furthermore,
$\bar\s|_{a_i}^{b_i}=\s|_{a_i}^{b_i}$ for every index $i\in I$,
and $\{t_i\}_{i\in I}$ is a set of basepoints for $(\p,\bar \s)$,
i.e.
$$\p^*(t_i)=\sigma_i\cdot\max_{[a_i,b_i]}\left\{\sigma_i\cdot
\p^*\right\},$$ where $\sigma_i$ denotes the sign
taken by both $\frac{d\s}{dt}$ and $\frac{d\bar \s}{dt}$ on the open
interval $J_i=(a_i,b_i)$. Therefore
\begin{eqnarray*}
& &\int_{a_i}^{b_i}\p^*(t)\cdot\frac{d\bar \s}{dt}(t)\ dt=
\int_{a_i}^{b_i}\p^*(t)\cdot\sigma_i\cdot\left|\frac{d\bar
\s}{dt}(t)\right|\ dt \le
\\
&\le& \int_{a_i}^{b_i} \max_{[a_i,b_i]}\left\{\p^*
\cdot\sigma_i\right\} \cdot\frac{d\bar \s}{dt}(t)\cdot\sigma_i\ dt
=
\\
&=&\int_{a_i}^{b_i}\p^*(t_i)\cdot\frac{d\bar \s}{dt}(t)\
dt =\p^*(t_i)\cdot\bar\s|_{a_i}^{b_i}=\p^*(t_i)\cdot\s|_{a_i}^{b_i}.
\end{eqnarray*}
Hence,
\begin{eqnarray}\label{inequality}
& &\int_{-\infty}^{+\infty}\p^*(t)\cdot\frac{d(\s\circ h)}{dt}(t)\
dt \le \sum_{i\in I} \s|_{a_i}^{b_i}\cdot\p^*(t_i)=F(\p,\s).
\end{eqnarray}
By substituting $\s$ with $-\s$ in the previous inequality
(observe  that $\mathcal{J}(-\s)=\mathcal{J}(\s)$ and
$H_{-\s}=H_\s$), we get
\begin{eqnarray}\label{inequality2}
& &-\int_{-\infty}^{+\infty}\p^*(t)\cdot\frac{d(\s\circ
h)}{dt}(t)\ dt \le \sum_{i\in I}
-\s|_{a_i}^{b_i}\cdot\p^*(\tilde t_i)=F(\p,-\s)
\end{eqnarray}
where $\{\tilde t_i\}_{i\in I}$ is a set of basepoints for $(\p,-\bar \s)$.

It follows that
\begin{eqnarray}\label{inequality3}
& &\left|\int_{-\infty}^{+\infty}\p^*(t)\cdot\frac{d(\s\circ
h)}{dt}(t)\ dt \right|\le \max \left\{F(\p,\s), F(\p,-\s)\right\}
\end{eqnarray}
for every $h\in H_\s$. Therefore, the inequality
$$\sup_{h\in H_\s}\left|\int_{-\infty}^{+\infty}\p^*(t)\cdot\frac{d (\s\circ h)}{dt}(t)\ dt\right|\le
\max \left\{F(\p,\s),F(\p,-\s)\right\}$$ holds.
\medskip

\textbf{PART 2: }$\sup_{h\in
H_\s}\left|\int_{-\infty}^{+\infty}\p^*(t)\cdot\frac{d (\s\circ
h)}{dt}(t)\ dt\right|\ge \max \left\{F(\p,\s), F(\p,-\s)\right\}$.

Assume $\bar \s=\s\circ h$ with $h\in H_\s$. Once more,
$\mathcal{J}(\bar \s)=\mathcal{J}(\s)$, ${\sign} \frac{d\bar
\s}{ds}={\sign} \frac{d \s}{ds}$ and $\{t_i\}_{i\in I}$ is a set
of basepoints for $(\p,\bar \s)$, i.e.
$$\p^*(t_{i})=\sigma_i\cdot\max_{[a_i,b_i]}\left\{
\sigma_i \cdot\p^*\right\},$$ where $\sigma_i$ denotes the sign
taken by $\frac{d\s}{ds}$ and $\frac{d\bar \s}{ds}$ on the
interval $J_i$.




Let us choose an $\eps>0$. In order to avoid the problem of some
$t_i$'s possibly belonging to the boundary of $J_i$, we define a new set
$\{t'_i\}_{i\in I}$: for each interval
$J_i=(a_i,b_i)\in \mathcal{J}(\bar \s)$ we choose a
$t'_{i}\in (a_i,b_i)$ such that
$|\p^*(t_{i})-\p^*(t'_{i})|<\frac{\eps}{\ 2^i}$.

For each positive integer $k\le |I|$ let us choose a positive real number
$\eta$ such that
$\eta<\min\left\{\frac{t'_{i}-a_i}{2},\frac{b_i-t'_{i}}{2},\eps\right\}$,
for every $i\le k-1$. Then we can consider an orientation-preserving
diffeomorphism $h_{(\eta,k)}\in H_\s$ such that
\begin{description}
  \item[i)] for every $i\le k-1$, $h_{(\eta,k)}$ maps $(t'_{i}-\frac{\eta}{\ 2^i},t'_{i}+\frac{\eta}{\ 2^i})$ onto
$(a_i+\frac{\eta}{\ 2^i},b_i-\frac{\eta}{\ 2^i})$;
  \item[ii)] the restriction of $h_{(\eta,k)}$ to the set $\R-\bigcup_{i=0}^{k-1} J_i$
  is the identity.
\end{description}
Recall also that, because $h_{(\eta,k)}\in H_\s$,  for every index $i$ the map $h_{(\eta,k)}$ takes the
interval $J_i$ onto itself.

Since $\p$ is  continuous, we can also assume to choose
$\eta$ so small that the inequality
$|\p^*(t'_{i})-\p^*(h^{-1}_{(\eta,k)}(t))|<\frac{\eps}{2^i}$ holds
for any $t\in (a_i+\frac{\eta}{\ 2^i},b_i-\frac{\eta}{\
2^i})$ and any index $i\le k-1$.

Therefore, for any index $i\le k-1$, by setting $t=h_{(\eta,k)}(s)$
and recalling that the sign of $\frac{d{ \s}} {dt}$ is constant in
$(a_i,b_i)$ (we shall use this fact in several following passages),
\begin{eqnarray*}
&   & \left| \s|_{a_i}^{b_i} \cdot\p^*(t'_{i}) -\int_{a_i}^{b_i}\p^*
(s)\cdot\frac{d({ \s}\circ h_{(\eta,k)})} {ds}(s)\ ds \right|=
\\    
&=& \left| \s|_{a_i}^{b_i}\cdot\p^*(t'_{i}) -\int_{a_i}^{b_i}\p^*
(s)\cdot\frac{d{ \s}} {ds}(h_{(\eta,k)}(s))\cdot\frac{dh_{(\eta,k)}}{ds}(s)\
ds \right|=
\\    
&=& \left| \s|_{a_i}^{b_i}\cdot\p^*(t'_{i}) -\int_{a_i}^{b_i}\p^*
(h^{-1}_{(\eta,k)}(t))\cdot\frac{d{ \s}} {dt}(t)\ dt \right|\le
\\    
&\le& \left| \s|_{a_i}^{b_i}\cdot\p^*(t'_{i})
-\int_{a_i+\frac{\eta}{\ 2^i}}^{b_i-\frac{\eta}{\ 2^i}}\p^*
(h^{-1}_{(\eta,k)}(t))\cdot\frac{d{ \s}} {dt}(t)\ dt \right|+
\\    
& & + \ \ \ \left| \int_{a_i}^{a_i+\frac{\eta}{\ 2^i}}\p^*
(h^{-1}_{(\eta,k)}(t))\cdot\frac{d{ \s}} {dt}(t)\ dt \right|+ \left|
\int_{b_i-\frac{\eta}{\ 2^i}}^{b_i}\p^*
(h^{-1}_{(\eta,k)}(t))\cdot\frac{d{ \s}} {dt}(t)\ dt \right| \le
\\    
&\le& \left|
\s|_{a_i}^{b_i}\cdot\p^*(t'_{i})-\p^*(t'_{i})\cdot\int_{a_i+\frac{\eta}{\
2^i}}^{b_i-\frac{\eta}{\ 2^i}} \frac{d{ \s}} {dt}(t)\ dt
\right|+ \left|\int_{a_i+\frac{\eta}{\
2^i}}^{b_i-\frac{\eta}{\ 2^i}}
\left(\p^*(t'_{i})-\p^*(h^{-1}_{(\eta,k)}(t))\right)\cdot \frac{d{
\s}} {dt}(t)\ dt \right|+
\\    
& & + \ \ \ \max |\p|\cdot \int_{a_i}^{a_i+\frac{\eta}{\
2^i}}\left|\frac{d{ \s}} {dt}(t)\right|\ dt + \max
|\p|\cdot\int_{b_i-\frac{\eta}{\ 2^i}}^{b_i}\left|\frac{d{
\s}} {dt}(t)\right|\ dt \le
\\    
&\le& \left| \left(\s|_{a_i}^{b_i}-\s|_{a_i+\eta/
2^i}^{b_i-\eta/2^i}\right) \cdot\p^*(t'_{i})\right| +
\int_{a_i+\frac{\eta}{2^i}}^{b_i-\frac{\eta}{2^i}}
\left|\p^*(t'_{i})-\p^*(h^{-1}_{(\eta,k)}(t))\right|\cdot
\left|\frac{d{ \s}} {dt}(t)\right|\ dt +
\\    
& & + \ \ \ \left|\s|_{a_i}^{a_i+\eta/2^i}+\s|_{b_i-\eta/2^i}^{b_i}\right|\cdot\max |\p|
\le
\\    
&\le& \left|\s|_{a_i}^{a_i+\eta/2^i}+\s|_{b_i-\eta/2^i}^{b_i}\right|\cdot\max |\p|+\frac{\eps}{2^i}\cdot\int_{a_i+\frac{\eta}{\
2^i}}^{b_i-\frac{\eta}{\ 2^i}} \left|\frac{d{ \s}}
{dt}(t)\right|\ dt +
\left|\s|_{a_i}^{a_i+\eta/2^i}+\s|_{b_i-\eta/2^i}^{b_i}\right|\cdot\max |\p| =
\\    
&=& \left|\s|_{a_i}^{a_i+\eta/2^i}+\s|_{b_i-\eta/2^i}^{b_i}\right|\cdot\max |\p|+\frac{\eps}{2^i}\cdot\left|\s|_{a_i+\eta/2^i}^{b_i-\eta/2^i} \right| +\left|\s|_{a_i}^{a_i+\eta/2^i}+\s|_{b_i-\eta/2^i}^{b_i}\right|\cdot\max |\p| =
\\    
&=& 2\left(\left|\s|_{a_i}^{a_i+\eta/2^i}\right|+\left|\s|_{b_i-\eta/2^i}^{b_i}\right|\right)\cdot\max |\p|+\frac{\eps}{2^i}\cdot\left|\s|_{a_i+\eta/2^i}^{b_i-\eta/2^i} \right| \le
\\    
&\le& 4\cdot \frac{\eta}{\ 2^i}\cdot \max \left|\frac{d{
\s}}{dt}\right|\cdot \max |\p|+\frac{\eps}{2^i}\cdot\left(\max {
\s}-\min { \s}\right).
\end{eqnarray*}

It follows that, for every positive integer $k\le |I|$ and every $\eps>0$, a
small enough positive $\eta=\eta(k,\eps)\le\eps$ exists such that,
denoting by $A_k$  the set $\bigcup_{i=0}^{k-1} (a_i,b_i)$, we
have
\begin{eqnarray*}
& &\left| \sum_{i=0}^{k-1} \left(\s|_{a_i}^{b_i}\cdot\p^*(t_{i}) -\int_{A_k}\p^*
(s)\cdot\frac{d({ \s}\circ h_{(\eta,k)})} {ds}(s)\ ds
\right)\right|\le\\
&\le &\left| \sum_{i=0}^{k-1}
\s|_{a_i}^{b_i}\cdot\left(\p^*(t_{i})-\p^*(t'_{i})\right)\right|+\left|
\sum_{i=0}^{k-1} \left(\s|_{a_i}^{b_i}\cdot\p^*(t'_{i}) -\int_{A_k}\p^*(s)\cdot\frac{d({
\s}\circ h_{(\eta,k)})} {ds}(s)\ ds
\right)\right|\le\\
&\le & V_\s\cdot\sum_{i=0}^{k-1} \frac{\eps}{\ 2^i}+
\sum_{i=0}^{k-1} \left|\s|_{a_i}^{b_i}\cdot\p^*(t'_{i}) -\int_{A_k}\p^*(s)\cdot\frac{d({
\s}\circ h_{(\eta,k)})} {ds}(s)\ ds
\right|\le\\
& &\ \ \ \le 2\eps\cdot V_\s+
8\eta\cdot\max \left|\frac{d{ \s}}{dt}\right|\cdot \max |\p|+2\eps\cdot(\max { \s}-\min { \s}).
\end{eqnarray*}
Hence
\begin{eqnarray*}
& &\left| \sum_{i=0}^{k-1} \s|_{a_i}^{b_i}\cdot\p^*(t_{i})
-\int_{-\infty}^{+\infty}\p^*(s)\cdot\frac{d({ \s}\circ
h_{(\eta,k)})}{ds}(s)\ ds
\right|\le\\
& &\ \ \ \le 2\eps\cdot\left(V_\s+\max { \s}-\min { \s}\right)+
8\eta\cdot \max \left|\frac{d{
\s}}{dt}\right|\cdot\max |\p|+\left|\int_{\R-A_k}\p^*(s)\cdot\frac{d( \s\circ
h_{(\eta,k)})}{ds}(s)\ ds \right|.
\end{eqnarray*}
By definition of $h_{(\eta,k)}$, the function $ \s\circ
h_{(\eta,k)}$ equals $ \s$ on ${\R-A_k}$ and hence, when $I=\N$,
$$\lim_{k\to|I|}\int_{\R-A_k}\p^*(s)\cdot\frac{d( \s\circ
h_{(\eta,k)})}{ds}(s)\ ds=\lim_{k\to|I|}\int_{\R-A_k}\p^*
(s)\cdot\frac{d \s}{ds}(s)\ ds=0,$$
and analogously when $k=|I|$,
$$\int_{\R-A_k}\p^*(s)\cdot\frac{d( \s\circ
h_{(\eta,k)})}{ds}(s)\ ds=\int_{\R-A_k}\p^*
(s)\cdot\frac{d \s}{ds}(s)\ ds=0.$$
Therefore, recalling that
$\eta\le\eps$, the following inequality holds for every large enough $k$ (if $|I|=\infty$), and for $k=|I|$ (if $|I|<\infty$):
\begin{eqnarray*}\label{limit1}
& &\left| \sum_{i=0}^{k-1} \s|_{a_i}^{b_i}\cdot\p^*(t_{i})
-\int_{-\infty}^{+\infty}\p^*(s)\cdot\frac{d({ \s}\circ h_{(\eta,k)})}
{ds}(s)\ ds \right|\le \\
& & \le 2\eps\cdot\left(V_\s+\max { \s}-\min { \s}+4\cdot\max \left|\frac{d{ \s}}{dt}\right|\cdot\max |\p|\right)+2\eps.
\end{eqnarray*}
Hence
\begin{eqnarray*}
& &\left| \sum_{i=0}^{|I|-1} \s|_{a_i}^{b_i}\cdot\p^*(t_{i})
-\int_{-\infty}^{+\infty}\p^*(s)\cdot\frac{d({ \s}\circ h_{(\eta,k)})}
{ds}(s)\ ds \right|\le
\end{eqnarray*}
\begin{eqnarray*}\label{limit2}
& &\ \ \ \le 2\eps\cdot\left(V_\s+\max { \s}-\min { \s}+4\cdot\max \left|\frac{d{ \s}}{dt}\right|\cdot\max |\p|+1\right)+\left|\sum_{i=k}^{|I|-1}
\s|_{a_i}^{b_i}\cdot\p^*(t_{i})\right|
\end{eqnarray*}
(here and in the sequel, when $I=\N$, we set $|I|-1=\infty$).

Since $F(\p,\s)=\sum_{i=0}^{|I|-1} \s|_{a_i}^{b_i}\cdot\p^*(t_{i})$ and it is finite, if $k$ is large enough (in case $|I|=\infty$), or $k=|I|$ (in case $|I|<\infty$) we get $\left|\sum_{i=k}^{|I|-1}
\s|_{a_i}^{b_i}\cdot\p^*(t_{i})\right|\le 2\eps$. Hence
\begin{eqnarray*}
& &\left| F(\p,\s) -\int_{-\infty}^{+\infty}\p^*(s)\cdot\frac{d({
\s}\circ h_{(\eta,k)})} {ds}(s)\ ds \right|\le
\end{eqnarray*}
\begin{eqnarray}\label{limit3}
& &\ \ \ \le 2\eps\cdot\left(V_\s+\max { \s}-\min { \s}+4\cdot\max \left|\frac{d{ \s}}{dt}\right|\cdot\max |\p|+2\right).
\end{eqnarray}
Inequality (\ref{limit3}) proves that for every $ \bar\eps>0$ an
orientation-preserving diffeomorphism $h^+\in H_\s$ exists such that
$$\left|\int_{-\infty}^{+\infty}\p^*(t)\cdot\frac{d (\s\circ h^+)}{dt}(t)\ dt\right|\ge\int_{-\infty}^{+\infty}\p^*(t)\cdot\frac{d (\s\circ h^+)}{dt}(t)\ dt\ge F(\p,\s)- \bar\eps.$$

By substituting $\s$ with $-\s$ and observing that
$H_{-\s}=H_\s$ we get for every $ \bar\eps>0$ an
orientation-preserving diffeomorphism $h^-\in H_\s$ exists such that
$$\left|\int_{-\infty}^{+\infty}\p^*(t)\cdot\frac{d (\s\circ h^-)}{dt}(t)\ dt\right|\ge -\int_{-\infty}^{+\infty}\p^*(t)\cdot\frac{d (\s\circ h^-)}{dt}(t)\ dt\ge F(\p,-\s)- \bar\eps.$$
Therefore the inequality
$$\sup_{h\in H_\s}\left|\int_{-\infty}^{+\infty}\p^*(t)\cdot\frac{d (\s\circ h)}{dt}(t)\ dt\right|\ge
\max \left\{F(\p,\s),F(\p,-\s)\right\}$$holds. This implies our
claim.
\end{proof}

The next result, proved by applying the previous lemma,  motivates the introduction of the functional $F$.

\begin{prop}\label{trivialequalitylemma}
$\|\p\|_{[\s]}=\sup_{\hat \s\in [\s]}\max \left\{F(\p,\hat
\s),F(\p,-\hat \s)\right\}$.
\end{prop}

\begin{proof}
Lemma~\ref{interpretationlemma} implies that
\begin{eqnarray*}
& &\|\p\|_{[\s]}=\sup_{\hat \s\in
[\s]}\left|\int_{-\infty}^{+\infty}\p^*(t)\cdot\frac{d \hat \s}{dt}(t)\
dt\right|=
\\
&=& \sup_{\hat \s\in [\s]}\sup_{h\in H_{\hat
\s}}\left|\int_{-\infty}^{+\infty}\p^*(t)\cdot\frac{d (\hat \s\circ
h)}{dt}(t)\ dt\right|= \sup_{\hat \s\in [\s]}\max \left\{F(\p,\hat
\s),F(\p,-\hat \s)\right\}.
\end{eqnarray*}
\end{proof}

\subsection{Standard RPI-norms as absolute values of linear combinations of Dirac deltas}

The next theorem simplifies the computation of the standard
RPI-norms, bypassing the concept of basepoints for the pair
$(\p,\s)$. First we define  a new set $\mathcal{T}(\p,\s)$ based on the natural ordering $\preceq$, previously introduced on the set $\mathcal{J}(\s)$. We recall that  $J_i\preceq J_j$ if and only if $a_i\le a_j$. This order does not need to coincide with the order induced by the index $i$. We also recall that the set $I$ indexing $\mathcal{J}(\s)$ is assumed to be either the finite set $\{0,1\ldots,n-1\}$ or the set $\N$.

\begin{defn}\label{natord}
Let $\p,\s\in AS^1(\R)-\{{\mathbf 0}\}$, and let $[a,b]$ be the smallest closed interval containing the support of $\frac{\ d\p^*}{dt}$. We shall denote  by $\mathcal{T}(\p,\s)$
the set of all the  sequences $(\tau_i)$ of points of $[a,b]$ such that
$J_{i_1}\preceq J_{i_2}$ implies $\tau_{i_1}\le \tau_{i_2}$, for every $i_1,i_2\in I$, and $\tau_i=b$ if $i\not\in I$. The sequences in $\mathcal{T}(\p,\s)$ will be said to be \emph{compatible with $(\p,\s)$}. The terms $\tau_i=b$ with $i\not\in I$ will be called \emph{dummy terms} of the sequence $(\tau_i)\in \mathcal{T}(\p,\s)$.
\end{defn}

\begin{rem}\label{completamento}
When $I=\mathbb{N}$, an ordered set of basepoints is a compatible sequence itself. When $I$ is finite, starting from a set of basepoints, we can obtain a compatible sequence by adding infinitely many dummy terms $\tau_i=b$ to our finite sequence. As a matter of fact, the dummy terms will  not be used in our computations.
\end{rem}

\begin{rem}\label{serieconvergenti}
For every $(\tau_i)\in\mathcal{T}(\p,\s)$, the series $\sum_{i\in I} \s|_{a_i}^{b_i}\cdot\p^*(\tau_i)$ converges, since $\left|\sum_{i\in I} \s|_{a_i}^{b_i}\cdot\p^*(\tau_i)\right|\le \max|\p|\cdot V_\s$.
\end{rem}

\begin{thm}\label{teoremadirappresentazionegenerico}
 If $\p,\s\neq \mathbf{0}$ then $\|\p\|_{[\s]}=\sup_{(\tau_i)\in
\mathcal{T}(\p,\s)}\left|\sum_{i\in I} \s|_{a_i}^{b_i}\cdot\p^*(\tau_i)\right|$.
\end{thm}


\begin{proof}
If $\hat \s\in [\s]$, an orientation-preserving diffeomorphism
$h\in D^1_+(\R)$ exists such that $\s=\hat \s\circ h$, and
$\mathcal{J}(\hat \s)=\{h(J_i)\}_{i\in I}$. For each index $i$, we define
the open interval $(\alpha_i,\beta_i)$ by setting $(\alpha_i,\beta_i)=h(J_i)$. Let us choose
a set $\{\hat t^+_i\}$ of basepoints for the pair $(\p,\hat \s)$
by taking each $\hat t^+_i$ in the closure of the interval $h(J_i)$. Analogously,
let us choose a set $\{\hat t^-_i\}$ of basepoints for the pair
$(\p,-\hat \s)$ by taking each $\hat t^-_i$ in the closure of the interval
$h(J_i)$.

Because of the definition of $F$,
we obtain that
$$F(\p,\hat \s)=\sum_{i\in I}
\hat\s|_{\alpha_i}^{\beta_i}\cdot\p^*(\hat t^+_i)=\sum_{i\in I} \s|_{a_i}^{b_i}\cdot\p^*(\hat t^+_i)$$
$$F(\p,-\hat \s)=-\sum_{i\in I}
\hat\s|_{\alpha_i}^{\beta_i}\cdot\p^*(\hat t^-_i)=-\sum_{i\in I} \s|_{a_i}^{b_i}\cdot\p^*(\hat t^-_i).$$
Therefore,
 $$\max \left\{F(\p,\hat \s),F(\p,-\hat \s)\right\}\le
\max \left\{\left|
\sum_{i\in I} \s|_{a_i}^{b_i}\cdot\p^*(\hat t^+_i)
\right|,\left|
\sum_{i\in I} \s|_{a_i}^{b_i}\cdot\p^*(\hat t^-_i)
\right|\right\}.$$
Let $\left(\hat t^+_i\right)$ be a compatible sequence for $(\p,\hat\s)$ obtained from the set of basepoints  $\{\hat t^+_i\}$, and analogously, let
$\left(\hat t^-_i\right)$ be a compatible sequence for   $(\p,-\hat\s)$  obtained from the set of basepoints  $\{\hat t^-_i\}$ (recall Remark \ref{completamento}). It is easy to see that $\left(\hat t^+_i\right)$ and $\left(\hat t^-_i\right)$ are compatible sequences also for $(\p,\s)$. Therefore, for any $\hat \s\in [\s]$ the
inequality
$$\max
\left\{F(\p,\hat \s),F(\p,-\hat \s)\right\}\le \sup_{(\tau_i)\in
\mathcal{T}(\p,\s)}\left|
\sum_{i\in I} \s|_{a_i}^{b_i}\cdot\p^*(\tau_i)
\right|$$ holds. From
Prop.~\ref{trivialequalitylemma} the inequality
$$\|\p\|_{[\s]}\le\sup_{(\tau_i)\in
\mathcal{T}(\p,\s)}\left|
\sum_{i\in I} \s|_{a_i}^{b_i}\cdot\p^*(\tau_i)
\right|$$ follows.

On the other hand, because of the continuity of $\p$, for
every $(\tau_i)\in \mathcal{T}(\p,\s)$, every positive integer $k\le |I|$ and every
$\eps>0$,  an orientation-preserving diffeomorphism $h_{k,\eps}\in
D^1_+(\R)$ exists, such that the distance between the number
$\p^*(\tau_i)$ and each value in the set $\p^*(h_{k,\eps}(J_i))$ is
not greater than $\eps$, for $i\le k-1$ (it is sufficient to choose a diffeomorphism taking each $J_i$ into an interval contained in a small neighborhood of $\tau_i$). We point out that here we are using the hypothesis that $(\tau_i)$ is a sequence compatible with $(\p,\s)$.

Let us consider  $\hat
\s_{k,\eps}=\s\circ h_{k,\eps}^{-1}$, and choose a set $\{\hat
t^+_i\}$ of basepoints for the pair $(\p,\hat \s_{k,\eps})$ and a
set $\{\hat t^-_i\}$ of basepoints for the pair $(\p,-\hat
\s_{k,\eps})$.
For each index $i\le k-1$, we define
the open interval $(\alpha'_i,\beta'_i)$ by setting $(\alpha'_i,\beta'_i)=h_{k,\eps}(J_i)$.
As before, we observe that $\mathcal{J}(\hat
\s_{k,\eps})=\{h_{k,\eps}(J_i)\}_{i\in I}$ and that
$\hat t^+_i,\hat t^-_i$ belong to the closure of  $h_{k,\eps}(J_i)$ for every index $i\le k-1$.
Since $\left|\p^*\left(\hat
t^+_i\right)-\p^*\left(\tau_i\right)\right|\le \eps$
and $\left|\p^*\left(\hat
t^-_i\right)-\p^*\left(\tau_i\right)\right|\le
\eps$, from
$\hat\s_{k,\eps}|_{\alpha'_i}^{\beta'_i}=\s|_{a_i}^{b_i}$
we get

$$\left|\sum_{i=0}^{k-1} \hat \s_{k,\eps}|_{\alpha'_i}^{\beta'_i}\cdot\p^*\left(\hat t^+_i\right)-
\sum_{i=0}^{k-1} \s|_{a_i}^{b_i}\cdot\p^*\left(\tau_i\right)\right|\le \eps\cdot V_\s$$

and

$$\left|\sum_{i=0}^{k-1} \hat \s_{k,\eps}|_{\alpha'_i}^{\beta'_i}\cdot\p^*\left(\hat t^-_i\right)-
\sum_{i=0}^{k-1} \s|_{a_i}^{b_i}\cdot\p^*\left(\tau_i\right)\right|\le \eps\cdot V_\s.$$

Furthermore,
$$\left|\sum_{i=k}^{|I|-1} \hat \s_{k,\eps}|_{\alpha'_i}^{\beta'_i}\cdot\p^*\left(\hat t^+_i\right)\right|
\le \max |\p|\cdot\sum_{i=k}^{|I|-1}
\left|\hat \s_{k,\eps}|_{\alpha'_i}^{\beta'_i}\right|=
\max |\p|\cdot\sum_{i=k}^{|I|-1} \left|\s|_{a_i}^{b_i}\right|.$$

Since $\s$ is a function of bounded variation, if $k$ is large
enough we get

$$\left|\sum_{i=k}^{|I|-1} \hat \s_{k,\eps}|_{\alpha'_i}^{\beta'_i}\cdot\p^*\left(\hat t^+_i\right)\right|
\le \eps.$$

(Here and in the following, \emph{$k$ large enough} means $k=|I|$ if $|I|$ is finite, and in this case every empty summation is assumed to take the value $0$.)

Analogously, if $k$ is large enough we get

$$\left|\sum_{i=k}^{|I|-1} \hat \s_{k,\eps}|_{\alpha'_i}^{\beta'_i}\cdot\p^*\left(\hat t^-_i\right)\right|
\le \eps, \ \ \left|\sum_{i=k}^{|I|-1} \s|_{a_i}^{b_i}\cdot\p^*\left(\tau_i\right)\right|
\le \eps.$$

Therefore for every $\eps>0$ we can find a large enough index $k$
such that

$$\left|\sum_{i\in I} \hat \s_{k,\eps}|_{\alpha'_i}^{\beta'_i}\cdot\p^*\left(\hat t^+_i\right)-
\sum_{i\in I} \s|_{a_i}^{b_i}\cdot\p^*\left(\tau_i\right)\right|\le \eps\cdot V_\s+2\eps$$

and

$$\left|\sum_{i\in I} \hat \s_{k,\eps}|_{\alpha'_i}^{\beta'_i}\cdot\p^*\left(\hat t^-_i\right)-
\sum_{i\in I} \s|_{a_i}^{b_i}\cdot\p^*\left(\tau_i\right)\right|\le \eps\cdot V_\s+2\eps.$$

By recalling the definition of $F$ we obtain

 $$F(\p,\hat \s_{k,\eps})+\eps\cdot V_\s+ 2\eps\ge
\sum_{i\in I} \s|_{a_i}^{b_i}\cdot\p^*\left(\tau_i\right),$$

$$F(\p,-\hat \s_{k,\eps})+\eps\cdot V_\s+ 2\eps\ge
-\sum_{i\in I} \s|_{a_i}^{b_i}\cdot\p^*\left(\tau_i\right).$$

These two last inequalities and the arbitrariness of $\eps$ imply that

$$\sup_{\hat \s\in [\s]}\max \left\{F(\p,\hat \s),F(\p,-\hat \s)\right\}\ge \left|\sum_{i\in I} \s|_{a_i}^{b_i}\cdot\p^*\left(\tau_i\right)\right|$$ for any
$(\tau_i)\in \mathcal{T}(\p,\s)$.
It follows that

$$\sup_{\hat \s\in [\s]}\max \left\{F(\p,\hat \s),F(\p,-\hat \s)\right\}\ge \sup_{(\tau_i)\in \mathcal{T}(\p,\s)}\left|\sum_{i\in I} \s|_{a_i}^{b_i}\cdot\p^*\left(\tau_i\right)\right|.$$

From
Prop.~\ref{trivialequalitylemma} the inequality
$$\|\p\|_{[\s]}\ge\sup_{(\tau_i)\in
\mathcal{T}(\p,\s)}\left|\sum_{i\in I} \s|_{a_i}^{b_i}\cdot\p^*\left(\tau_i\right)\right|$$ follows. Hence our statement
is proved.
\end{proof}

\subsection{Optimal sequences in $\mathcal{T}(\p,\s)$}

The previous Theorem~\ref{teoremadirappresentazionegenerico} raises an interesting
issue: is the $\sup$ equaling $\|\p\|_{[\s]}$ actually a $\max$? In Prop. \ref{propoptimal} we shall give an affirmative answer to this question.

We consider the following definition:

\begin{defn}
  Every sequence $(\bar \tau_i)\in \mathcal{T}(\p,\s)$ such that $\|\p\|_{[\s]}=\left|\sum_{i\in I}
  \s|_{a_i}^{b_i}\cdot\p^*(\bar\tau_i)\right|$  is said to be \emph{optimal for $(\p,\s)$}.
  The set of all optimal sequences for $(\p,\s)$  will be denoted by $\mathcal{O}(\p,\s)$.
  \end{defn}

In the sequel, optimal sequences will be obtained as convergent subsequences of sequences in $\mathcal{T}(\p,\s)$. Therefore we shall need the following lemma.

\begin{lem}\label{sottosucc}
From each sequence $(T^n)$ of sequences belonging to $\mathcal{T}(\p,\s)$ it is possible to extract a subsequence that pointwise converges to a sequence $\bar T\in \mathcal{T}(\p,\s)$.
\end{lem}

\begin{proof}
Consider the smallest interval
$[a,b]$ containing   the compact support of
$\frac{d\p^*}{dt}$.
Let $(T^n)$ be a sequence of sequences belonging to $\mathcal{T}(\p,\s)$. For every fixed $n$, $T^n=(t^n_i)$ with $t^n_i\in [a,b]$. Hence, the sequence $(T^n)$ admits a  subsequence $\left(T^{n^0_r}\right)$ (varying $r$) such that $\left(t^{n^0_r}_0\right)$ converges to some $\bar \tau_0\in [a,b]$. The sequence $\left(T^{n^0_r}\right)$ admits a  subsequence $\left(T^{n^1_r}\right)$ such that $\left(t^{n^1_r}_1\right)$ converges to some $\bar \tau_1\in [a,b]$. Since $\left(T^{n^1_r}\right)$ is a subsequence of $\left(T^{n^0_r}\right)$, it still holds that $\left(t^{n^1_r}_0\right)$ converges to  $\bar \tau_0$. By iterating this argument, for every $s\in \mathbb{N}$, we can extract a subsequence $\left(T^{n^{s+1}_r}\right)$ from $\left(T^{n^{s}_r}\right)$ such  that the sequence $\left(t_i^{n^{s+1}_r}\right)$ converges to $\bar \tau_i\in [a,b]$ (varying $r$) for every  fixed $i$ with $0\le i\le s+1$. By construction, the diagonal sequence  $\left(T^{n^r_r}\right)$ is a subsequence of $\left(T^{n}\right)$. Furthermore, for every $s\in \N$  $\left(T^{n^r_r}\right)$ is a subsequence of $\left(T^{n^s_r}\right)$, if we ignore the first $s$ terms of both these sequences.  Hence $\lim_{r\to +\infty}t_i^{n^r_r}=\bar \tau_i$ for every fixed $i\in \mathbb{N}$. Finally, $(\bar\tau_i)\in \mathcal{T}(\p,\s)$ since each $\bar\tau_i\in [a,b]$, and for every $i_1,i_2\in I$, if $J_{i_1}\preceq J_{i_2}$ then $t_{i_1}^{n^{r}_{r}}\le t_{i_2}^{n^{r}_{r}}$ for every $r$, implying that  $\lim_{r\to +\infty}t_{i_1}^{n^{r}_{r}}\le \lim_{r\to +\infty}t_{i_2}^{n^{r}_{r}}$. Hence the sequence $\left(T^{n^r_r}\right)$ proves our statement.
\end{proof}

In the next pages, each sequence obtained by the method described in the proof of the previous lemma will be said to be  \emph{``obtained by a diagonalization process''}.

\begin{prop} \label{propoptimal}
  If $\p,\s\neq \mathbf{0}$ then  $\mathcal{O}(\p,\s)$ is not empty.
  \end{prop}

\begin{proof} On the basis of Theorem~\ref{teoremadirappresentazionegenerico}, for each  $n\in \N$ we can take a  sequence $T^n=(t^n_i)\in \mathcal{T}(\p,\s)$
in such a way that $\|\p\|_{[\s]}=\lim_{n\to
\infty}\left|\sum_{i\in I}
  \s|_{a_i}^{b_i}\cdot\p^*(t^n_i)\right|$. By Lemma \ref{sottosucc}, the sequence of sequences $(T^n)$ admits a subsequence $(T^{n_r})$ that pointwise converges to a sequence $\bar T=(\bar\tau_i)$ compatible with $(\p,\s)$.
If we denote each sequence $T^{n_r}$ by $(t_i^{n_r})$ (varying $i$), the following equalities hold:

\begin{equation}
\label{limsum}
\left|\sum_{i\in I} \s|_{a_i}^{b_i}\cdot\p^*(\bar\tau_i)\right|=\left|\sum_{i\in I} \s|_{a_i}^{b_i}\cdot\lim_{r\to\infty}\p^*\left(t^{n_r}_i\right)\right|=\lim_{r\to\infty}\left|\sum_{i\in I} \s|_{a_i}^{b_i}\cdot\p^*\left(t^{n_r}_i\right)\right|=\|\p\|_{[\s]}
\end{equation}
where the second equality follows from the fact that $\s$ has  bounded variation.
\end{proof}


However, a stronger result holds, stating that there exist optimal sequences for $(\p,\s)$ containing only critical points for $\p^*$.

\begin{lem}\label{sottosuccopt}
From each sequence $(T^n)$ of sequences in $\mathcal{O}(\p,\s)$ it is possible to extract a subsequence that pointwise converges to a sequence $\bar T\in \mathcal{O}(\p,\s)$.
\end{lem}

\begin{proof}
By Lemma \ref{sottosucc}, $(T^n)$ admits a subsequence
$\left(T^{n_r}\right)$  that pointwise converges to a sequence $\bar T\in \mathcal{T}(\p,\s)$.
By recalling that $\s$ has  bounded variation, it is easy to verify that $\bar T\in
\mathcal{O}(\p,\s)$ (cf. the equalities (\ref{limsum}) in the previous proof of Proposition~\ref{propoptimal}).
\end{proof}

Now we can prove the following result, improving Proposition \ref{propoptimal}.

\begin{prop} \label{propoptimalcritical}
   If $\p,\s\neq \mathbf{0}$ then a sequence $(\hat \tau_i)\in\mathcal{O}(\p,\s)$ exists,
  where each $\hat \tau_i$ is a critical point of $\p^*$.
  \end{prop}

\begin{proof}
Proposition~\ref{propoptimal} shows that the set
$\mathcal{O}(\p,\s)$ of all optimal sequences in $\mathcal{T}(\p,\s)$
 is not empty. Let $K_{\p^*}$ be the set of all critical points
of $\p^*$. If $T=(\tau_i)\in \mathcal{O}(\p,\s)$ we define the weight
$w(T)=\sum_{i\in I} \gamma_i\cdot \left|\s|_{a_i}^{b_i}\right|$, where
$\gamma_i=\min \{\tau_i-x|x\in K_{\p^*}, x\le \tau_i \}$ (in other
words $\gamma_i$ is the distance between $\tau_i$ and the first
critical point of $\p^*$  on its left). This positive terms series
converges since it is smaller than $(b-a)\cdot\sum_{i\in I}
\left|\s|_{a_i}^{b_i}\right|\le (b-a)\cdot V_\s$.

For any  $n\in \N$ we can take a  sequence
$T^n=(\tau^n_i)\in \mathcal{O}(\p,\s)$ in such a way that
$\lim_{n\to \infty}w(T^n)=\inf_{T\in \mathcal{O}(\p,\s)}w(T)$.
By Lemma \ref{sottosuccopt},  we can extract from $(T^n)$ a subsequence pointwise converging to an optimal sequence  $\widehat{T}=(\hat\tau_i)$.

Let us set $\gamma^n_i=\min \{\tau^n_i-x|x\in K_{\p^*}, x\le \tau^n_i \}$ for every $n,i\in\N$, and
$\widehat{\gamma}_i=\min \{\hat\tau_i-x|x\in K_{\p^*}, x\le \hat\tau_i \}$ for every $i\in\N$. Since the set $K_{\p^*}$ is closed, we can easily prove that, for every $i\in \mathbb{N}$, either $\widehat{\gamma}_i=\lim_{n\to\infty}\gamma^n_i$ or $\widehat{\gamma}_i=0$ although $\lim_{n\to\infty}\gamma^n_i\ne 0$. By recalling once again that $\s$ has  bounded variation, it follows that
$\lim_{n\to\infty}w(T^n)\ge w(\widehat{T})$. Since $\lim_{n\to \infty}w(T^n)=\inf_{T\in \mathcal{O}(\p,\s)}w(T)$, we get $w(\widehat{T})=\inf_{T\in \mathcal{O}(\p,\s)}w(T)$.

Now we prove by contradiction that $w(\widehat{T})=0$. Assume
$w(\widehat{T})>0$. Then an index $j\in I$ exists such that
$\hat\tau_j\not \in K_{\p^*}$ and, since $K_{\p^*}$ is closed, we can find an $\eta>0$ for
which the closure of the open interval $U=(\hat\tau_j-\eta,\hat\tau_j+\eta)$ does
not meet $K_{\p^*}$. We want to show that we can move all points of $\widehat{T}$ in $U$ leftwards, and get an optimal sequence with a weight that is strictly less than $w(\widehat{T})$. This will generate our contradiction.

In order to do that, let us consider a $C^1$ function
$\rho:\mathbb{R}\to \mathbb{R}$ such that
\begin{itemize}
\item $\rho=0$ outside $U$;
\item $\rho>0$ in $U$;
\item $\max_{\overline{U}}\left|\frac{d\rho}{dt}\right|< \min_{\overline{U}} \left|\frac{\ d\p^*}{dt}\right|$.
\end{itemize}

The last hypothesis guarantees that the function $\p^+=\p^*+\rho$
is a (possibly orientation-reversing) diffeomorphism from $U$ onto
its image $\p^+(U)$. Since  $\p^+(U)=\p^*(U)$, we can consider the
function from $U$ to $U$ that takes each point $t$ to the unique
point $t'$ such that $\p^*(t')=\p^+(t)$ (observe that either both $\p^+$ and $\p^*$ are strictly increasing in $U$ or both $\p^+$ and $\p^*$ are strictly decreasing in $U$). We can
extend this function to a function $h^+:\R\to \R$ by defining it
to equal the identity outside $U$. It is immediate to verify that
$h^+$ is an orientation-preserving diffeomorphism, since $\frac{\
d\p^*}{dt}$ and $\frac{\ d\p^+}{dt}$ take the same sign in $U$.
Analogously, the function $\p^-=\p^*-\rho$ is a diffeomorphism
from $U$ onto its image $\p^-(U)=\p^*(U)$. Hence we can consider
the function from $U$ to $U$ that takes each point $t$ to the
unique point $t'$ such that $\p^*(t')=\p^-(t)$. We can extend this
function to a function $h^-:\R\to \R$ by defining it to equal the
identity outside $U$, and $h^-$ is an orientation-preserving
diffeomorphism.

Now, let us define two new sequences $(\tau_i^+)$ and $(\tau_i^-)$.
For every $i\in \N$ we set $\tau_i^+=h^+(\hat\tau_i)$, $\tau_i^-=h^-(\hat\tau_i)$,
so that $\p^*(\tau_i^+)=\p^+(\hat\tau_i)$ and $\p^*(\tau_i^-)=\p^-(\hat\tau_i)$.

Since $h^+$ is an orientation-preserving diffeomorphism, $\hat\tau_i\le \hat\tau_j$ if and only if $\tau_i^+\le \tau_j^+$. That means that also $(\tau_i^+)$ is a sequence compatible with $(\p,\s)$. Analogously, also $(\tau_i^-)$ results to be a sequence compatible with $(\p,\s)$.
Since the sequence $(\hat\tau_i)$ is optimal, if
$\sum_{i\in I}
  \s|_{a_i}^{b_i}\cdot\p^*(\hat\tau_i)\ge 0$
 the following
statements hold:
$$\sum_{i\in I}
  \s|_{a_i}^{b_i}\cdot\p^*(\hat\tau_i) - \sum_{i\in I}
  \s|_{a_i}^{b_i}\cdot\p^*(\tau^+_i)=-\sum_{i\in I}
  \s|_{a_i}^{b_i}\cdot\rho(\hat\tau_i)\ge 0,$$

$$\sum_{i\in I}
  \s|_{a_i}^{b_i}\cdot\p^*(\hat\tau_i) - \sum_{i\in I}
  \s|_{a_i}^{b_i}\cdot\p^*(\tau^-_i)=\sum_{i\in I}
  \s|_{a_i}^{b_i}\cdot\rho(\hat\tau_i)\ge 0,$$

  and hence $\sum_{i\in I}
  \s|_{a_i}^{b_i}\cdot\rho(\hat\tau_i)=0$.

On the other hand, if $\sum_{i\in I}
  \s|_{a_i}^{b_i}\cdot\p^*(\hat\tau_i)< 0$
 the optimality of $(\hat\tau_i)$ implies the following
statements:

$$-\sum_{i\in I}
  \s|_{a_i}^{b_i}\cdot\p^*(\hat\tau_i) + \sum_{i\in I}
  \s|_{a_i}^{b_i}\cdot\p^*(\tau^+_i)=\sum_{i\in I}
  \s|_{a_i}^{b_i}\cdot\rho(\hat\tau_i)\ge 0,$$

$$-\sum_{i\in I}
  \s|_{a_i}^{b_i}\cdot\p^*(\hat\tau_i) + \sum_{i\in I}
  \s|_{a_i}^{b_i}\cdot\p^*(\tau^-_i)=-\sum_{i\in I}
  \s|_{a_i}^{b_i}\cdot\rho(\hat\tau_i)\ge 0,$$

and hence $\sum_{i\in I}
  \s|_{a_i}^{b_i}\cdot\rho(\hat\tau_i)=0$.

  Therefore, in any case $\sum_{i\in I}
  \s|_{a_i}^{b_i}\cdot\rho(\hat\tau_i)=0$, implying
  that
  $$\sum_{i\in I}
  \s|_{a_i}^{b_i}\cdot\p^*(\hat\tau_i)=\sum_{i\in I}
  \s|_{a_i}^{b_i}\cdot\p^*(\tau^+_i)=\sum_{i\in I}
  \s|_{a_i}^{b_i}\cdot\p^*(\tau^-_i).$$
  It follows that also the
  sequences $T^+=(\tau^+_i)$ and $T^-=(\tau^-_i)$ belong to $\mathcal{O}(\p,\s)$. Moreover, since $\rho(t)> 0$ if $t\in U$, it holds that either $h^+$ or $h^-$ moves every point in $U$ leftwards (according to whether $\frac{\ d\p^*}{dt}$ is negative or positive in $U$, respectively), while both of them do not move the points outside $U$. Therefore either $w(T^+)<w(\widehat{T})$ or
  $w(T^-)<w(\widehat{T})$ must hold, against our hypotheses. Hence the
equality $w(\widehat{T})=0$ is proved. It follows that
$\widehat{\gamma}_i=0$ for every index $i\in \N$, i.e. every $\hat\tau_i$ is a
critical point for $\p^*$.
\end{proof}

Proposition~\ref{propoptimalcritical} allows us to obtain immediately the next
useful result, strengthening Theorem~\ref{teoremadirappresentazionegenerico}. We state first a new definition.

\begin{defn}
Let  $\p,\s\in AS^1(\R)-\{\mathbf{0}\}$. We denote by
$\mathcal{C}(\p,\s)$  the set of all sequences $(\hat
\tau_i)\in\mathcal{O}(\p,\s)$   such that $\hat\tau_i$ is a critical
point of $\p^*$ for every index $i\in I$. We shall say that these sequences
are the \emph{optimal critical sequences for $(\p,\s)$}.
\end{defn}

\begin{thm}[\textbf{Representation Theorem}]\label{teoremadirappresentazione}
Let $\p,\s\in AS^1(\R)-\{\mathbf{0}\}$. Then $$\|\p\|_{[\s]}=\max_{(\hat
\tau_i)\in \mathcal{C}(\p,\s)}\left|\sum_{i\in I} \s|_{a_i}^{b_i}\cdot\p^*(\hat
\tau_i)\right|.$$
\end{thm}

\begin{rem}\label{remcordiscretizzazione}
Another way to express Theorem~\ref{teoremadirappresentazione} is stating
that the standard RPI-norm $\|\p\|_{[\s]}$ equals the value $\max_{(\hat \tau_i)\in
\mathcal{C}(\p,\s)}\left|\sum_{i\in I}
\s|_{a_i}^{b_i}\cdot\delta_{\hat \tau_i}(\p^*)\right|$. In other
words, previous Theorem~\ref{teoremadirappresentazione} makes available an
equivalent discrete definition for standard RPI-norms. This
definition allows for easier computations.
\end{rem}

We conclude this section with a remark.

\begin{rem}
\label{tuttenormestandard} If a $k$-tuple $(m_0,\ldots,m_{k-1})\in
\mathbb{R}^k$ with $m_i\neq 0$ for at least one index $i$ is
given, then a  function $\s\in AS^1(\R)-\{\mathbf{0}\}$   exists
such that, for every $\p\in AS^1(\R)$, the value
$\max_{\tau_0\le\tau_1\le\ldots\le\tau_{k-1}}\left|\sum_{i=0}^{k-1}
m_i\cdot\delta_{\tau_i}(\p^*)\right|$ equals the standard RPI-norm
$\|\p\|_{[\s]}$. In fact, the Representation Theorem shows that it
is sufficient to choose some points $a_0<b_0<a_1<b_1<\ldots
<a_{k-1}<b_{k-1}$ and an almost sigmoidal $C^1$-function $\s$ such
that $\s$ is monotone in $(a_i,b_i)$ and $\s|_{a_i}^{b_i}=m_i$ for
every $0\le i\le k-1$, while $\frac{d\s}{dt}=0$ outside the set
$\bigcup_{i=0}^{k-1} [a_i,b_i]$.
\end{rem}

In other words, this means that any finite linear combination of
Dirac deltas corresponds to a standard RPI-norm. This is a partial
converse of Remark~\ref{remcordiscretizzazione}, stating that any
standard RPI-norm corresponds to a (not necessarily finite) linear
combination of Dirac deltas. It might be interesting to know under
which hypotheses the statement seen in
Remark~\ref{tuttenormestandard} is true for series of Dirac
deltas.

\section{Relationship between RPI-norms and standard RPI-norms}
\label{DEP}

A piecewise monotone almost sigmoidal function is an almost sigmoidal function that is monotone in each connected component of the complement of a finite set. In this section we shall prove a key result in this paper,
showing that all the RPI-norms of piecewise monotone $C^1_c$-functions are
determined by standard RPI-norms (Theorem~\ref{RT}).

Before dealing with the technical details of our proofs, it may be
useful to sketch the underlying ideas. The basic question to be
answered could be formulated in this way: \emph{``How can we use
the information contained in the standard RPI-norms in order to
reconstruct the function $\p$?''} In order to make this point
clear, let us consider for  instance a function $\s$ that   is
associated with the linear combination of Dirac deltas
$\Sigma_3=\delta_{t_0}-\delta_{t_1}+\delta_{t_2}$, where the
values $t_0,t_1,t_2$ are set equal to three suitable critical
points of $\p^*$, according to the Representation Theorem
\ref{teoremadirappresentazione}. In order to get some more
information about $\p$, we have to change $\s$ (and consequently
$\Sigma_3$). The simplest way to change $\s$ and $\Sigma_3$ is to
slightly perturb one of the three weights $1,-1,1$ in our linear
combination of deltas. E.g., we can consider the linear
combination
$\Sigma_3^\eps=(1+\eps)\delta_{t_0}-\delta_{t_1}+\delta_{t_2}$,
associated with a suitable function $\s_\eps$. For $\eps$ small
enough, the choice of $t_0,t_1,t_2$ for which
$\left|\Sigma_3^\eps(\p^*)\right|$ is maximum  allows also $\left|\Sigma_3(\p^*)\right|$ to
attain its maximum value. This ``invariance of the basepoints
$t_0,t_1,t_2$ with respect to small changes of the weights'' and
the fact that $\Sigma_3^\eps(\p^*)$ and $\Sigma_3(\p^*)$ take the
same sign will allow us to write the following equalities:

$$\|\p\|_{[\s_\eps]}-\|\p\|_{[\s]}=\left|\Sigma_3^\eps(\p^*)\right|-\left|\Sigma_3(\p^*)\right|=
\eps\cdot
\delta_{t_0}(\p^*)\cdot\sign
\left(\Sigma_3(\p^*)\right)=\eps\cdot
\p^*(t_0)\cdot\sign
\left(\Sigma_3(\p^*)\right).$$

It follows that the function $\|\p\|_{[\s_\eps]}$
is differentiable with respect to $\eps$ and that $\frac{d\|\p\|_{[\s_\eps]}}{d\eps}(0)$ equals
$\p^*(t_0)\cdot\sign
\left(\Sigma_3(\p^*)\right)$. So,  we get that
the value $v_0$ taken by $\p^*$ at the critical
point $t_0$  equals
$\frac{d\|\p\|_{[\s_\eps]}}{d\eps}(0)\cdot\sign
\left(\Sigma_3(\p^*)\right)$.

We can repeat the above procedure to obtain the values taken by
$\p^*$ at the other two critical points $t_1$ and $t_2$. Hence, so
far, we know that $\p^*$ is a function that takes the values $v_0$,
$v_1$ and $v_2$ in this order, when $t$ varies from $-\infty$ to
$+\infty$. By taking functions $\s$ with an increasing number of
oscillations (i.e. $\s$ corresponding to
$\sum_{i=0}^{n-1}(-1)^i\delta_{t_i}$), we obtain more and more
information about the values of $\p^*$ at critical points. Since $\p$ is piecewise monotone, if the
number of oscillations of $\psi$ is large enough then $\p^*$ is
monotone between two suitable critical points $t_i$ and
$t_{i+1}$ of $\p^*$.
Obviously, we have not enough
information to locate the points $t_i$'s, but we are able to
reconstruct the oscillations of $\p^*\cdot\sign
\left(\Sigma_3(\p^*)\right)$, up to reparametrization and an
arbitrarily small error with respect to the variation norm.
Roughly speaking, these are the ideas we are going to use.

First of all we recall the formal definition of piecewise monotone
function.

\begin{defn}
We say that $f:\R\to\R$ is \emph{piecewise monotone} if a finite
set $W\subset\R$ exists such that $f$ is monotone in each
connected component of the complement of $W$.  Each such a set $W$
will be said to be a \emph{separating set for $f$}. If $W$ is also
minimal with respect to inclusion, it will be said to be a
\emph{minimal separating set for $f$}. We define $l(f)$ as the
minimum of the cardinalities of the separating sets for $f$.
\end{defn}

Obviously, $l(S)=0$ and $l(\Lambda)=1$. Note that if $f\in
C^1_c(\R)$ and $f\neq \mathbf{0}$ then $l(f)\ge 1$. It is easy to
show that all minimal separating sets for $f$ take the same
cardinality $l(f)$. This follows from the next simple proposition
(we omit the immediate proof):

\begin{prop}
\label{W} Let $f:\R\to\R$ be a piecewise monotone function. Let
$W=\{t_0,\ldots,t_{m}\}$ and $W'=\{t'_0,\ldots,t'_{n}\}$ be two
separating sets for $f$, with $t_0<t_1<\ldots <t_m$ and $W'$
minimal.  If for some $i$ and $j$ it holds that
$t_j<t'_i<t_{j+1}$, then $f$ is constant either in $[t_j,t'_i]$ or
in $[t'_i,t_{j+1}]$.
\end{prop}

We also observe that the concept of piecewise monotone almost
sigmoidal $C^1$-function is invariant under reparametrization, and
that $l(\p)=l(\p^*)$. Moreover, the points of a minimal separating
set for $\p^*$ are necessarily critical points for $\p^*$.

In the rest of this section, when $\p$ is a piecewise monotone
almost sigmoidal $C^1$-function with non-empty compact support
(i.e. $\p\neq \mathbf{0}$), we let $[a,b]$ denote the minimal
interval containing the support of $\p^*$. Moreover, if
$\left\{t_0,\ldots,t_{l(\p)-1}\right\}$ is a minimal separating
set for $\p^*$, we assume it is increasingly ordered and we
define  $c=\min_{0\le i\le
l(\p)}\left|\p(t_i)-\p(t_{i-1})\right|$, where we set $t_{-1}=a$
and $t_{l(\p)}=b$. This meaning of the symbols $t_{-1}$ and
$t_{l(\p)}$ will be maintained in the following pages.

Before proceeding, we need to introduce a new family of functions.

\begin{defn}
Let $n\ge 1$. For every  vector $e=(\eps_0,\ldots,\eps_{n-1})\in \mathbb{R}^n$
  we define the functions $S_n,S_n^e:\R\to \R$
  by setting $S_n(t)=\sum_{i=0}^{n-1}(-1)^iS(t-2i)$ and $S_n^e(t)=\sum_{i=0}^{n-1}\left((-1)^i+\eps_i\right)S(t-2i)$.
  \end{defn}

In plain words, the function $S_n^e$ is a perturbation of the
function $S_n$. In particular, for $e=(0,\ldots,0)$, the function
$S_n^e$ equals the function $S_n$. We observe that the functions
$S_n^e$ are piecewise polynomial and belong to
$AS^1(\R)$. We also note that $\|\p\|_{[S_n]}\le \|\p\|_{[S_{n+1}]}$ for every $n\ge 1$ and every $\p\in AS^1(\R)$.

The following lemma is a key passage towards the proof of the Reconstruction Theorem \ref{RT} for piecewise monotone functions in $C^1_c(\R)$.

\begin{lem}\label{ma}
Let $\p\ne \bf 0$ be a piecewise monotone $C^1$-function with compact support. Let $[a,b]$ denote the minimal interval containing the support of $\p^*$. If $\left\{t_0,\ldots,t_{l(\p)-1}\right\}$ is an (increasingly ordered) minimal separating set for $\p^*$ and $c=\min_{0\le i\le l(\p)}\left|\p(t_i)-\p(t_{i-1})\right|$ (where we set $t_{-1}=a$ and $t_{l(\p)}=b$),  then
\begin{enumerate}
\item $\|\p\|_{[S_n]}=\left|\sum_{i=0}^{l(\p)-1}(-1)^{i}\cdot\p^*(t_i)\right|=\frac{1}{2}\cdot V_\p$ for any $n\ge l(\p)$;
\item $\|\p\|_{[S_{n}]}\le\|\p\|_{[S_{l(\p)}]}-c$, for $1\le n< l(\p)$.
\end{enumerate}
\end{lem}

\begin{proof}
We start proving {\em 1.} Let us consider the finite set $\mathcal{J}(S_n)$ of
all maximal open intervals of $\R$ where $\frac{dS_n}{dt}$ does not
vanish. We note that $\mathcal{J}(S_n)=\left\{J_i\right\}_{i\in I}$, where $J_i=(-1+2i,1+2i)$, $I=\{0,\ldots,n-1\}$,  and $S_n\left|_{-1+2i}^{1+2i}\right.=(-1)^{i}$.

From the Bounding Lemma \ref{boundsC1} for functions in $C^1_c(\R)$, we obtain that $\|\p\|_{[S_{n}]}\le V_\p/2$ (we point out that Remark \ref{remcordiscretizzazione} implies $\|\Lambda\|_{[S_n]}=1$). Observing that $\left|\sum_{i=0}^{l(\p)-1}(-1)^{i}\cdot\p^*(t_i)\right|=V_\p/2$, we easily conclude that $(t_0,\ldots,t_{l(\p)-1},b,b,b,\ldots )$ is an optimal sequence for $(\p, S_n)$ and that $\|\p\|_{[S_{n}]}= V_\p/2$ for any $n\ge l(\p)$.

We now prove {\em 2.} Let $1\le n<l(\p)$. Since $\|\p\|_{[S_n]}\le \|\p\|_{[S_{n+1}]}$, it is sufficient to prove that $\|\p\|_{[S_{l(\p)-1}]}\le\|\p\|_{[S_{l(\p)}]}-c$.

If $l(\p)=2$ then $\max\p\ge c$, $-\min\p\ge c$ so that
$$\|\p\|_{[S_{l(\p)}]}-\|\p\|_{[S_{l(\p)-1}]} = \|\p\|_{[\Lambda]}-\|\p\|_{[S]}=\max\p-\min\p-\max|\p|\ge c.$$

Let us now assume that $l(\p)\ge 3$.
Let $T=(\tau_i)$ be an optimal sequence for $(\p,S_{l(\p)-1})$, increasingly ordered  so that $\|\p\|_{[S_{l(\p)-1}]}=\left|\sum_{i=0}^{l(\p)-2}(-1)^{i}\cdot\p^*(\tau_i)\right|$.
The key point of the proof relies in understanding where to place $\tau_0, \tau_1,\ldots, \tau_{l(\p)-2}$,  in order to achieve optimality.

First of all, we can assume that
\begin{description}
\item[i)] \emph{the only repeated point in the sequence $T$ is  $b$}
\end{description}
since consecutive points appear with opposite weights, and that
\begin{description}
\item[ii)]  \emph{each $\tau_i$ belongs to $\{t_0,t_1,\ldots, t_{l(\p)-1}\}\cup \{a,b\}$.}
\end{description}
Indeed, if  statement \textbf{ii)} were false we could easily find a better sequence than $T$ by using the monotonicity of $\p^*$ outside the set
$\{t_0,t_1,\ldots, t_{l(\p)-1}\}$, and hence $T$ would not be optimal.

Now, let us note that $\tau_{l(\p)-1}=b$, since $\tau_{l(\p)-1}$ is the first dummy point of $T \in \mathcal{T}\left(\p,S_{l(\p)-1}\right)$ (recall Definition \ref{natord}). Therefore, if $\tau_0=a$  we obtain
$$\sum_{i=0}^{l(\p)-2}(-1)^{i}\cdot\p^*(\tau_i)=-\sum_{i=0}^{l(\p)-2}(-1)^{i}\cdot\p^*(\tau_{i+1})$$
because $\p^*(a)= \p^*(b)=0$. Hence, possibly by replacing the optimal sequence $(\tau_i)$ with the new optimal sequence $(\tau_{i+1})$, we can substitute {\bf ii)} with
\begin{description}
\item[iii)]  \emph{each $\tau_i$ belongs to $\{t_0,t_1,\ldots, t_{l(\p)-1}\}\cup \{b\}$.}
\end{description}

%
%

Possibly changing $\p$ with $-\p$, we can assume that
\begin{description}
\item[iv)]  $\sum_{i=0}^{l(\p)-2}(-1)^{i}\cdot\p^*(\tau_i)\ge 0$.
\end{description}

Now, two mutually exclusive cases are possible:
\begin{description}
\item[A)]  \emph{$\p^*$ is increasing in $[a,t_0]$; }
\item[B)]  \emph{$\p^*$ is decreasing in $[a,t_0]$.}
\end{description}

We can prove the following property (in the following we shall set $t_{l(\p)}=b$):
\begin{description}
\item[v-A)]  \emph{In case \textbf{A)} for each $i$ with $0\le i\le l(\p)-2$, exactly one index $j$ exists such that $0\le j\le l(\p)$ and  $\tau_i=t_j$. Furthermore, $i-j$ is even.}
\end{description}
The new statement is that $i-j$ is even. In order to prove \textbf{v-A)} by contradiction, let us assume  that $i$ is even and $j$ is odd. Since \textbf{A)} holds and $j$ is odd, $\p^*$ is decreasing in $[t_{j-1},t_j]$ (this follows from the definition of minimal separating set). Although $\p^*$ may be not strictly decreasing in $[t_{j-1},t_j]$, the definition of minimal separating set for $\p^*$ implies that $\p^*(t_j)<\p^*(t_{j-1})$. Hence $(-1)^i\cdot \p^*(\tau_i)=\p^*(t_j)<\p^*(t_{j-1})=(-1)^i\cdot \p^*(t_{j-1})$. Now, let us assume  that $i$ is odd and $j$ is even. Since \textbf{A)} holds and $j$ is even, $\p^*$ is increasing in $[t_{j-1},t_j]$. Hence $(-1)^i\cdot \p^*(\tau_i)=-\p^*(t_j)<-\p^*(t_{j-1})=(-1)^i\cdot \p^*(t_{j-1})$.  Therefore, if $i$ and $j$ did not have the same parity the sequence $T$ would not be optimal, since $(-1)^i\cdot \p^*(\tau_i)<(-1)^i\cdot \p^*(t_{j-1})$ and we could obtain a better sequence by redefining $\tau_i=t_{j-1}$. This proves that $i-j$ is even, in case \textbf{A)}.

Moreover, we can prove the following property (once again, we shall set $t_{l(\p)}=b$):
\begin{description}
\item[v-B)] \emph{In case \textbf{B)} for each $i$ with $0\le i\le l(\p)-2$, exactly one index $j$ exists such that $0\le j\le l(\p)$ and  $\tau_i=t_j$. Furthermore, $i-j$ is odd.}
\end{description}
The new statement is that $i-j$ is odd. In order to prove \textbf{v-B)} by contradiction, let us assume  that both $i$ and $j$ are even. Since \textbf{B)} holds and $j$ is even, $\p^*$ is decreasing in $[t_{j-1},t_j]$ (this follows from the definition of minimal separating set). Although $\p^*$ may be not strictly decreasing in $[t_{j-1},t_j]$, the definition of minimal separating set for $\p^*$ implies that $\p^*(t_j)<\p^*(t_{j-1})$. Hence $(-1)^i\cdot \p^*(\tau_i)=\p^*(t_j)<\p^*(t_{j-1})=(-1)^i\cdot \p^*(t_{j-1})$. Now, let us assume that both $i$ and $j$ are odd. Since \textbf{B)} holds and $j$ is odd, $\p^*$ is increasing in $[t_{j-1},t_j]$. Hence $(-1)^i\cdot \p^*(\tau_i)=-\p^*(t_j)<-\p^*(t_{j-1})=(-1)^i\cdot \p^*(t_{j-1})$.  Therefore, if $i$ and $j$ had the same parity the sequence $T$ would not be optimal, since $(-1)^i\cdot \p^*(\tau_i)<(-1)^i\cdot \p^*(t_{j-1})$ and we could obtain a better sequence by redefining $\tau_i=t_{j-1}$. This proves that $i-j$ is odd, in case \textbf{B)}.

Properties \textbf{v-A)} and \textbf{v-B)} imply that both in case \textbf{A)} and in case \textbf{B)}
\begin{description}
\item[v)]  \emph{for each $i$ with $0\le i\le l(\p)-2$, exactly one index $j$ exists such that $0\le j\le l(\p)$ and $\tau_i=t_j$. Furthermore, if $i$ is even then $\p^*(t_j)>\p^*(t_{j-1})$, while if $i$ is odd then $\p^*(t_j)<\p^*(t_{j-1})$ (here we set $t_{-1}=a$).}
\end{description}

Indeed, if $i$ is even and \textbf{A)} holds property \textbf{v-A)} implies that $j$ is even, so that $\p^*$ is increasing in $[t_{j-1},t_j]$.  If $i$ is even and \textbf{B)} holds property \textbf{v-B)} implies that $j$ is odd, so that $\p^*$ is increasing in $[t_{j-1},t_j]$. If $i$ is odd and \textbf{A)} holds property \textbf{v-A)} implies that $j$ is odd, so that $\p^*$ is decreasing in $[t_{j-1},t_j]$. If $i$ is odd and \textbf{B)} holds property \textbf{v-B)} implies that $j$ is even, so that $\p^*$ is decreasing in $[t_{j-1},t_j]$. In summary, if $i$ is even then $\p^*$ is increasing in $[t_{j-1},t_j]$, while if $i$ is odd then $\p^*$ is decreasing in $[t_{j-1},t_j]$. This proves property \textbf{v)}.

Now we can prove that
\begin{description}
\item[vi)]  \emph{$\tau_{l(\p)-3}\neq b$ (and hence $\tau_i\neq b$ for every $i\le l(\p)-3$).}
\end{description}

In order to check {\bf vi)}, let us proceed by contradiction and assume that $\tau_{l(\p)-3}=b$. Note that necessarily  $\tau_{l(\p)-2}=\tau_{l(\p)-1}=b$. Since the number of points $\tau_i$ that are different from $b$ is at most $l(\p)-3$, at least three points of the separating set $\{t_0,t_1,\ldots, t_{l(\p)-1}\}$  do not coincide with any $\tau_i$.
Let us consider the largest index $k\le l(\p)-1$ such that $t_{k}$ does not belong to the set $\left\{\tau_i\right\}$.
Since, by definition, $t_{k+1}\in \left\{\tau_i\right\}$, {\bf v-A)} and {\bf v-B)} imply that $t_{k-1}$ does not belong to the set $\left\{\tau_i\right\}$. Indeed, properties {\bf v-A)} and {\bf v-B)} guarantee that for every index $i$ the points $\tau_i$ and $\tau_{i+1}$ are separated by an even number of points $t_j$ (possibly $0$). Let us define $\bar\iota=\min\left\{i:\tau_i= t_{k+1}\right\}$ (we set $t_{l(\p)}=b$) and consider the sequence
 $T'=(\tau'_i)$ obtained by setting $\tau'_i=\tau_i$ if $i< \bar\iota$, $\tau'_{\bar \iota}=t_{k-1}$, $\tau'_{\bar \iota +1}=t_{k}$ and $\tau'_{i}=\tau_{i-2}$ if $i\ge \bar\iota+2$.
In other words $T'$ differs from $T$ for the insertion of the pair of points $t_{k-1}$ and $t_{k}$. Since $\tau_{\bar \iota}=t_{k+1}$ with $\bar\iota\le l(\p)-3$, by applying {\bf v)} we have $\p^*(t_{k+1})>\p^*(t_k)$ if $\bar\iota$ is even, and $\p^*(t_{k+1})<\p^*(t_k)$ if $\bar\iota$ is odd. As a consequence, $\p^*(t_{k-1})>\p^*(t_{k})$ if $\bar\iota$ is even, and $\p^*(t_{k-1})<\p^*(t_{k})$ if $\bar\iota$ is odd. Therefore, both for $\bar\iota$ even and for $\bar\iota$ odd we get
\begin{eqnarray*}
(-1)^{\bar \iota}\p^*(\tau'_{\bar \iota})+(-1)^{\bar \iota +1}\p^*(\tau'_{\bar \iota +1})=(-1)^{\bar \iota}\big(\p^*(t_{k-1})
-\p^*(t_{k})\big)>0.
\end{eqnarray*}
Therefore, we have that
\begin{eqnarray*}
\sum_{i=0}^{l(\p)-2}(-1)^{i}\cdot\p^*(\tau_i)& =& \sum_{i=0}^{\bar \iota -1}(-1)^{i}\cdot\p^*(\tau_i)+ \sum_{i=\bar \iota}^{l(\p)-2}(-1)^{i}\cdot\p^*(\tau_i)=\\
&  =& \sum_{i=0}^{\bar \iota -1}(-1)^{i}\cdot\p^*(\tau'_i)+\sum_{i=\bar \iota +2}^{l(\p)}(-1)^{i}\cdot\p^*(\tau'_i)<\\
&  < &\sum_{i=0}^{l(\p)}(-1)^{i}\cdot\p^*(\tau'_i),
\end{eqnarray*}
contradicting the optimality of $T$. Hence property \textbf{vi)} is proved.
\bigskip

The previous property \textbf{vi)} shows that $\left\{\tau_0, \tau_1,\ldots, \tau_{l(\p)-3}\right\}\subseteq\left\{t_0, t_1,\ldots, t_{l(\p)-1}\right\}$, and hence exactly $l(\p)-2$ points in the separating set must belong to the set $\left\{\tau_0, \tau_1,\ldots, \tau_{l(\p)-3}\right\}$. From all this we deduce that the only dispositions allowed for $\tau_0, \tau_1,\ldots, \tau_{l(\p)-2}$ are  the following ones:
\begin{enumerate}
\item $\tau_0=t_1,\tau_1=t_2,\ldots ,\tau_{l(\p)-2}=t_{l(\p)-1}$;
\item $\tau_0=t_0,\tau_1=t_1,\ldots ,\tau_{l(\p)-2}=t_{l(\p)-2}$;
\item an index $k$ with $1\le k\le l(\p)-1$ exists such that $\left\{t_0,\ldots,t_{l(\p)-1}\right\}=\left\{\tau_0,\ldots,\tau_{l(\p)-3}\right\}\cup \left\{t_{k-1},t_k\right\}$.
\end{enumerate}

The first two cases happen when $\tau_{l(\p)-2}\neq b$, the last case when $\tau_{l(\p)-2}= b$.
Note that if $\tau_{l(\p)-2}\neq b$, exactly $l(\p)-1$ points in the separating set must belong to the set $\left\{\tau_0, \tau_1,\ldots, \tau_{l(\p)-2}\right\}$, and recall once again that for every index $i$ the points $\tau_i$ and $\tau_{i+1}$ are separated by an even number of points $t_j$ (possibly $0$). Hence the only $t_j$'s that can miss in the set $\left\{\tau_0, \tau_1,\ldots, \tau_{l(\p)-2}\right\}$ are $t_0$ and $t_{l(\p)-1}$. This observation produces the cases 1 and 2. If $\tau_{l(\p)-2}= b$ a gap of two consecutive $t_j$'s is possible, implying the case 3.

In the first case,
since $\sign(\p^*(t_0))=\sign \left(\sum_{i=0}^{l(\p)-1}(-1)^{i}\cdot\p^*(t_i)\right)$ (we can easily verify this equality both in the cases \textbf{A)} and \textbf{B)}) and $c\le \left|\p^*(t_0)\right|\le \|\p\|_{[S_1]}\le \|\p\|_{[S_{l(\p)}]}$,
we have that
\begin{eqnarray*}
& &\|\p\|_{[S_{l(\p)-1}]}=\left|\sum_{i=0}^{l(\p)-2}(-1)^{i}\cdot\p^*(\tau_i)\right|=\left|\sum_{i=1}^{l(\p)-1}(-1)^{i}\cdot\p^*(t_i)\right|=\left|\sum_{i=0}^{l(\p)-1}(-1)^{i}\cdot\p^*(t_i)-\p^*(t_0)\right|=\\
& & =\left|\ \left|\sum_{i=0}^{l(\p)-1}(-1)^{i}\cdot\p^*(t_i)\right|-\left|\p^*(t_0)\right|\right|
=\left|\ \|\p\|_{[S_{l(\p)}]}-\left|\p^*(t_0)\right|\right|=\|\p\|_{[S_{l(\p)}]}-\left|\p^*(t_0)\right|\le \|\p\|_{[S_{l(\p)}]}-c.
\end{eqnarray*}
In the second case, the claim is proved analogously, since $c\le \left|\p^*(t_{l(\p)-1})\right|\le \|\p\|_{[S_1]}\le \|\p\|_{[S_{l(\p)}]}$ and
$$\sign\left((-1)^{l(\p)-1}\p^*(t_{l(\p)-1})\right)=\sign \left(\sum_{i=0}^{l(\p)-1}(-1)^{i}\cdot\p^*(t_i)\right).$$
This last equality immediately follows from the equality $\sign(\p^*(t_0))=\sign \left(\sum_{i=0}^{l(\p)-1}(-1)^{i}\cdot\p^*(t_i)\right)$, by substituting the function $\p^*(t)$ with the function $(-1)^{l(\p)-1}\p^*(-t)$ and inverting the order of the points in the separating set.

Therefore, in case $2$ we have that
\begin{eqnarray*}
& &\|\p\|_{[S_{l(\p)-1}]}=\left|\sum_{i=0}^{l(\p)-2}(-1)^{i}\cdot\p^*(\tau_i)\right|=\left|\sum_{i=0}^{l(\p)-2}(-1)^{i}\cdot\p^*(t_i)\right|=\\
& & =\left|\sum_{i=0}^{l(\p)-1}(-1)^{i}\cdot\p^*(t_i)-(-1)^{l(\p)-1}\p^*(t_{l(\p)-1})\right|=\\
& & =\left|\ \left|\sum_{i=0}^{l(\p)-1}(-1)^{i}\cdot\p^*(t_i)\right|-\left|\p^*(t_{l(\p)-1})\right|\right|\le \|\p\|_{[S_{l(\p)}]}-c.
\end{eqnarray*}
Finally, in case $3$,
\begin{eqnarray*}
 & &\|\p\|_{[S_{l(\p)-1}]}=\left|\sum_{i=0}^{l(\p)-2}(-1)^{i}\cdot\p^*(\tau_i)\right|=\left|\sum_{i=0}^{l(\p)-1}(-1)^{i}\cdot\p^*(t_i)-\left((-1)^{k-1}\p^*(t_{k-1})+ (-1)^{k}\p^*(t_{k})\right)\right|=\\
& &=   \left|\  \left|\sum_{i=0}^{l(\p)-1}(-1)^{i}\cdot\p^*(t_i)\right|-\left|(-1)^{k-1}\p^*(t_{k-1})+ (-1)^{k}\p^*(t_{k})\right|\right|\le \|\p\|_{[S_{l(\p)}]}-c,
\end{eqnarray*}
because
$$\sign\left((-1)^{k-1}\p^*(t_{k-1})+ (-1)^{k}\p^*(t_{k})\right)=\sign \left(\sum_{i=0}^{l(\p)-1}(-1)^{i}\cdot\p^*(t_i)\right)$$
(as can be easily verified both in case \textbf{A)} and in case \textbf{B)}) and
$$c\le \left| (-1)^{k-1}\p^*(t_{k-1})+ (-1)^{k}\p^*(t_{k}) \right|\le \|\p\|_{[S_2]}\le \|\p\|_{[S_{l(\p)}]}.$$
\end{proof}

An immediate consequence of the previous lemma is the following result, allowing us to deduce the value of $l(\p)$ from the knowledge of the standard RPI-norms $\|\p\|_{[S_{n}]}$, $n\ge 1$, when $\p$ is piecewise monotone and belongs to $C^1_c(\R)$.

\begin{cor}\label{corellediphi}
Let $\p\ne \bf0$ be a piecewise monotone $C^1$-function with compact support. The value $l(\p)$ is equal to the smallest integer $N$ such that $\|\p\|_{[S_{N}]}=\|\p\|_{[S_{n}]}$ for every $n\ge N$.
\end{cor}

In the following lemma, we consider the pairs $(\p, S_{l(\p)}^e)$ varying $e$,  and we show that for $e$ small enough they all admit the same optimal sequence.

\begin{lem}\label{bo}
Let $\p\ne \bf0$ be a piecewise monotone $C^1$-function with compact support. Let
 $W=\left\{t_0,\ldots,t_{l(\p)-1}\right\}$ be an (increasingly ordered) minimal separating set for $\p^*$.
For every $e=(\eps_0,\ldots,\eps_{l(\p)-1})\in \R^{l(\p)}$ with
$\max_i|\eps_i|<\min\left\{ \frac{c}{2\cdot l(\p)\cdot
\max|\p|},1\right\}$, it holds that
$(t_0,\ldots,t_{l(\p)-1},b,b,b,\ldots)$ is an optimal sequence for
$\left(\p,S_{l(\p)}^e\right)$, i.e.
$$\left|\sum_{i=0}^{l(\p)-1}\left((-1)^{i}+\eps_i\right)\cdot\p^*(t_i)\right|=\left\|\p\right\|_{\left[S_{l(\p)}^e\right]}.$$
\end{lem}

\begin{proof}
First of all we note that $c>0$ by definition.
  Let us consider the finite set $\mathcal{J}(S_{l(\p)}^e)$ of
all maximal open intervals of $\R$ where $\frac{dS_{l(\p)}^e}{dt}$ does not
vanish. Since $\max|\eps_i|< 1$, we have that $\mathcal{J}(S_{l(\p)}^e)=\left\{J_i\right\}_{i\in I}$, where $J_i=(-1+2i,1+2i)$, $I=\{0,\ldots,l(\p)-1\}$,  and $S_{l(\p)}^e\left|_{-1+2i}^{1+2i}\right.=(-1)^{i}+\eps_i$.

Consider  an optimal sequence for $\left(\p,S_{l(\p)}^e\right)$,
$T=(\tau_0,\ldots,\tau_{l(\p)-1},b,b,b,\ldots)$ (increasingly
ordered). As we have already done in the proof of Lemma \ref{ma},
we can assume that $\left\{\tau_i\right\}\subseteq W\cup
\left\{a,b\right\}$,
%
and that the only repeated point  in $T$ is $b$. We claim that
necessarily
$(\tau_0,\ldots,\tau_{l(\p)-1})=(t_0,\ldots,t_{l(\p)-1})$.
Otherwise, either at least one point at the beginning of $\left(\tau_0,\ldots,\tau_{l(\p)-1}\right)$ is equal to $a$, or at least one point at the end of the same $l(\p)$-tuple is equal to $b$. In this case
$\left\{\tau_0,\ldots,\tau_{l(\p)-1}\right\}-\left\{a,b\right\}$
is properly included in $\left\{t_0,\ldots,t_{l(\p)-1}\right\}$,
so that
$\left|\sum_{i=0}^{l(\p)-1}(-1)^{i}\p^*(\tau_i)\right|\le\left\|\p\right\|_{\left[S_{l(\p)-1}\right]}$, since $\p^*(a)=\p^*(b)=0$.
Hence Lemma \ref{ma} implies that
$$\left|\sum_{i=0}^{l(\p)-1}(-1)^{i}\p^*(t_i)\right|-\left|\sum_{i=0}^{l(\p)-1}(-1)^{i}\p^*(\tau_i)\right|=
\left\|\p\right\|_{\left[S_{l(\p)}\right]}-
\left|\sum_{i=0}^{l(\p)-1}(-1)^{i}\p^*(\tau_i)\right|\ge \left\|\p\right\|_{\left[S_{l(\p)}\right]}-
\left\|\p\right\|_{\left[S_{l(\p-1)}\right]}
\ge c$$
and thus
$$\left|\sum_{i=0}^{l(\p)-1}\left((-1)^{i}+\eps_i\right)\cdot\p^*(t_i)\right|-\left|\sum_{i=0}^{l(\p)-1}\left((-1)^{i}+\eps_i\right)\cdot\p^*(\tau_i)\right|\ge c-2\cdot l(\p)\cdot \max|\p|\cdot\max_i|\eps_i|>0 ,$$
contradicting the optimality of the sequence $T$.
\end{proof}

We now show that the standard RPI-norms of $\p$ with respect to $S_{l(\p)}^e$ allow us to obtain the values of $\p$ at the points of each minimal separating set.

\begin{lem}\label{derivabile}
Let $\p\ne \textbf{0}$ be a piecewise monotone $C^1$-function with compact support. Let us consider an (increasingly ordered)  minimal separating set $W=\left\{t_0,\ldots,t_{l(\p)-1}\right\}$ for $\p^*$.
For every $e=(\eps_0,\ldots,\eps_{l(\p)-1})\in \R^{l(\p)}$ with $\max_i|\eps_i|<\min\left\{ \frac{c}{2\cdot l(\p)\cdot \max|\p|},1\right\}$ it holds that
\begin{enumerate}
\item the function $\|\p\|_{\left[S_{l(\p)}^e\right]}$ is differentiable in
the variables $\eps_0,\ldots,\eps_{l(\p)-1}$;
\item  $\p^*(t_{i})=s\cdot\frac{\partial}{\partial\eps_i
}\|\p\|_{\left[S_{l(\p)}^e\right]}(0)$ for $0\le i\le l(\p)-1$,
 where $s$ is the sign of
$\sum_{i=0}^{l(\p)-1} (-1)^{i}\cdot\p^*(t_i)$.
\end{enumerate}
\end{lem}

\begin{proof}
From Lemma \ref{bo} it follows that for every $e=(\eps_0,\ldots,\eps_{l(\p)-1})\in \R^{l(\p)}$ with $|\eps_0|,\ldots,|\eps_{l(\p)-1}|< \min\left\{\frac{c}{2\cdot l(\p)\cdot \max|\p|},1\right\}$ it holds that $\|\p\|_{\left[S_{l(\p)}^e\right]}=\left|\sum_{i=0}^{l(\p)-1}
\left((-1)^{i}+\eps_{i}\right)\cdot\p^*(t_i)\right|$.
As a consequence,  the function $\|\p\|_{\left[S_{l(\p)}^e\right]}$ is differentiable in each
variable $\eps_i$ and
$$\frac{\partial}{\partial\eps_i
}\|\p\|_{\left[S_{l(\p)}^e\right]}(0)=s\cdot\p^*(t_{i})$$ for $0\le i\le l(\p)-1$,
 where $s$ is the sign of
$\sum_{i=0}^{l(\p)-1} (-1)^i\cdot\p^*(t_i)$ (note that
$s$ does not depend on $e$).
\end{proof}

We are now ready to prove that given  a piecewise monotone
$C^1$-function $\p$ with compact support, it is possible to construct a
piecewise polynomial almost sigmoidal $C^1$-function $\hat\p$  such that
$\hat\p$ approximates $\pm\p$ in the total variation norm up to
reparametrization, and $\|\p\|=\|\hat\p\|$ for any RPI-norm on
$AS^1(\R)$. The key point here is that this construction is based
just on the knowledge of the values taken by the standard
RPI-norms  at $\p$. In other words,  the RPI-norms of piecewise
monotone  functions in $C^1_c(\R)$  are determined by the standard
RPI-norms.

\begin{thm}[\textbf{Reconstruction Theorem for piecewise monotone functions in $C^1_c(\R)$}]\label{RT}
Assume that $\p$ is a piecewise monotone $C^1$-function with compact support.  Let $e=(\eps_0,\ldots,\eps_{l(\p)-1})\in \R^{l(\p)}$, where $|\eps_0|,\ldots,|\eps_{l(\p)-1}|\le \eps$. If  $0<\eps < \min\left\{ \frac{c}{2\cdot l(\p)\cdot \max|\p|},1\right\}$  we can  define the piecewise polynomial $C^1$-function
$\widehat{\p}\in AS^1(\R)$ by setting
$$\widehat{\p}(t)= \frac{\partial\|\p\|_{\left[S_{l(\p)}^e\right]}}{\partial\eps_{l(\p)-1}
}(0)\cdot S(t)+
\sum_{j=1}^{l(\p)-1}\left(\frac{\partial\|\p\|_{\left[S_{l(\p)}^e\right]}}{\partial\eps_{l(\p)-j-1}
}(0)-\frac{\partial\|\p\|_{\left[S_{l(\p)}^e\right]}}{\partial\eps_{{l(\p)-j}}
}(0)\right)\cdot S(t-2j),$$ so that the following statements hold:
\begin{enumerate}
\item there exists an orientation-preserving $C^1$-diffeomorphism $h:\R\to \R$ such that either $V_{\p\circ h-\widehat \p}\le \eps$ or  $V_{-\p\circ h-\widehat \p}\le \eps$;
\item $\|\p\|=\|\widehat\p\|$ for any RPI-norm on $AS^1(\R)$.
\end{enumerate}
\end{thm}

\begin{proof}
Let us assume $\p\ne \bf0$, otherwise the claims are trivial. Let
$\left\{t_0,t_1,\ldots,t_{l(\p)-1}\right\}$ be a minimal
separating set for $\p^*$, with $ t_0< t_1<\ldots < t_{l(\p)-1}$.
By applying Lemma \ref{derivabile}, for $i=0,\ldots, l(\p)-1$,
$\p^*(t_i)=s\cdot
\frac{\partial\|\p\|_{[S_{l(\p)}^e]}}{\partial\eps_i\ }(0)$, with
$s=\sign\left(\sum_{i=0}^{l(\p)-1}(-1)^i\cdot\p^*(t_i)\right)$.
Thus, the norms $\|\p\|_{[S_{l(\p)}^e]}$ (varying $e$) allow us to
determine, up to the sign, the value of $\p^*$ at  the points
$t_0,t_1, \ldots ,t_{l(\p)-1}$. Furthermore, we have that
$$s\cdot\widehat{\p}(t)=\p(-t_{l(\p)-1})\cdot S(t)+
\sum_{j=1}^{l(\p)-1}\left(\p(-t_{l(\p)-j-1})-\p(-
t_{l(\p)-j})\right)\cdot S(t-2j).$$

We observe that $\p$ is monotone  in each of the intervals
$(-\infty,-t_{l(\p)-1}]$, $[ -t_{l(\p)-1}, -t_{l(\p)-2}], \ldots ,
[ -t_{1}, -t_{0}]$, $[-t_{0},+\infty)$.
 From Proposition
\ref{generalizzazione} it follows that for any
$\eps>0$  there exist an $AS^1(\R)$ function $\p_\eps$ and an
orientation-preserving $C^1$-diffeomorphism $ h:\R\to\R$ such that
$V_{\p_\eps-\p}\le \eps$ and $\p_\eps\circ h=s\cdot\widehat{\p}$.
Therefore,
$$V_{s\cdot (\p\circ h)-\widehat{\p}}=V_{ \p\circ h-s\cdot\widehat{\p}}=V_{ \p\circ h-\p_\eps\circ h}=V_{ \p-\p_\eps}\le
\eps,$$
proving the first claim. Furthermore, by the Bounding Lemma,
$$\big|\|\p\|-\|\widehat{\p}\|\big|=\big|\|s\cdot (\p\circ h)\|-\|\widehat{\p}\|\big|\le \|s\cdot (\p\circ h)-\widehat{\p}\|\le V_{s\cdot (\p\circ h)-\widehat{\p}}\cdot\|S\|\le
\eps
\cdot\|S\|.$$
Thus, passing to the limit for $\eps$ tending to $0$, also the second claim is proved.
\end{proof}

\section{Open problems and conclusions}
\label{OP}

In this paper we have studied the main properties of the
reparametrization invariant norms on $AS^1(\R)$, focusing the
key role of standard RPI-norms. We have proved that these norms
allow for the reconstruction of any  piecewise monotone $C^1$-function with compact support up to reparametrization, so determining the
value of any other RPI-norm on the same function. 

However, many
problems remain open. First of all the theory has been developed
just for the space $AS^1(\R)$ and our last results require also to
assume that the considered functions are piecewise monotone and have compact support. The
extensions to less regular spaces could be desirable. 

Another
interesting extension could be the passage from spaces of
functions defined on $\R$ to spaces of functions defined on $\R^n$
or a closed manifold. We notice that in the case of a closed
manifold  $\M$ the extension could be intertwined with the study
of the topology of $\M$, requiring to substitute the role of the Dirac
delta with other distributions. 

Moreover, we have left open the question  about the existence of reparametrization invariant norms not obtainable as sup of standard RPI-norms.  

We postpone the research on these issues to other papers.

\subsection*{Acknowledgments}
We wish to thank L. Ambrosio, A. Cohen, D. Guidetti, S. Masnou  and B. Volzone for their helpful suggestions. Anyway, the
authors are solely responsible for any possible errors. Thanks to M. Ferri for his indispensable support and friendship.

This paper is dedicated to Matteo, Marta and Tommaso.

\bibliographystyle{amsplain}
\providecommand{\bysame}{\leavevmode\hbox
to3em{\hrulefill}\thinspace}

\end{document}